\documentclass[11pt]{amsart}

\usepackage{psfig}
\usepackage{color}
\usepackage{array,amscd,amsmath,amsthm}
\usepackage{amssymb}
\usepackage{graphicx}
\usepackage[all]{xy}
\usepackage{rotating}
\xyoption{2cell}
\xyoption{rotate}
\xyoption{line}
\usepackage{lscape}
\allowdisplaybreaks

\newtheorem{mytheorem}{Theorem}
\newtheorem{mycorollary}{Corollary}

\newtheorem{theorem}{Theorem}[section]
\newtheorem{proposition}[theorem]{Proposition}
\newtheorem{lemma}[theorem]{Lemma}

\newtheorem*{conjecture}{Conjecture}

\theoremstyle{definition}
\newtheorem{definition}[theorem]{Definition}
\newtheorem*{example}{Example}
\newtheorem*{remark}{Remark}
\newtheorem*{proposition*}{Proposition}
\newtheorem*{theorem*}{Theorem}
\newtheorem*{definition*}{Definition}
\newtheorem*{notation}{Notation}
\newtheorem*{strategy}{Strategy}
\newtheorem*{warning}{Warning}

\newcommand{\C}{\mathbb{C}}
\newcommand{\Z}{\mathbb{Z}}
\newcommand{\Q}{\mathbb{Q}}
\newcommand{\R}{\mathbb{R}}
\newcommand{\N}{\mathbb{N}}

\def\M{\mathcal{M}}
\def\Mcomb{\mathcal{M}^{comb}}
\def\Mbar{\overline{\M}}
\def\Mbarcomb{\Mbar^{comb}}
\def\Mbartri{\Mbar^{\triangle}}
\def\Mhat{\widehat{\M}^{comb}}

\def\Mfun{\mathfrak{M}}
\def\Mfunbar{\overline{\Mfun}}

\def\Cc{\mathcal{C}}

\def\T{\mathcal{T}}
\def\Tfun{\mathfrak{T}}

\def\pa{\partial}
\def\S{\Sigma}
\def\Scomb{\Sigma^{comb}}
\def\What{\widehat{W}}
\def\Wbar{\overline{W}}

\def\vol{\text{vol}}
\def\k{\kappa}
\def\rar{\rightarrow}
\def\arr#1#2{\stackrel{#1}{#2}}
\def\hra{\hookrightarrow}
\def\a{\alpha}
\def\b{\beta}
\def\ua{\mathbf {\a}}

\def\ub{\mathbf {\b}}
\def\Ao{A^\circ}
\def\A8{A_{\infty}}
\def\Ahat{\widehat{A}}

\def\l{\lambda}
\def\La{\Lambda}
\def\H{\mathcal{H}}
\def\i{\iota}
\def\g{\gamma}
\def\G{\Gamma}
\def\D{\Delta}
\def\Dc{\Delta^{\circ}}
\def\d{\delta}
\def\lra{\longrightarrow}
\def\lmt{\longmapsto}
\def\s{\sigma}
\def\e{\varepsilon}
\def\xb{\bar{x}}
\def\la{\langle}
\def\ra{\rangle}
\def\ol#1{\overline{#1}}
\def\wh#1{\widehat{#1}}

\begin{document}

\title{Combinatorial classes on $\Mbar_{g,n}$
are tautological}

\author{Gabriele Mondello}

\address{Dipartimento di Matematica \\
Universit\`a degli Studi di Bologna \\
Piazza di Porta S. Donato, 5 \\
40127 Bologna (Italy)}

\email{mondello@dm.unibo.it}


\begin{abstract}
The combinatorial description via ribbon graphs of the moduli space
of Riemann surfaces makes it possible
to define combinatorial cycles in a natural
way. Witten and Kontsevich first conjectured that these classes are polynomials
in the tautological classes. We answer affirmatively to this conjecture
and find recursively all the polynomials.
\end{abstract}


\maketitle


\begin{section}*{Introduction}
\begin{subsection}{History}
Let $g$ and $n$ be nonnegative integers such that
$2g-2+n>0$ and set $P=\{p_1,\dots,p_n\}$.
Denote by $\M_{g,P}$ the moduli space
of compact Riemann surfaces $S$ of genus $g$ with
an injection $P\hra S$.

In the early '80s a combinatorial description
of $\M_{g,P}$ was discovered.
Thanks to ideas
mainly due to Mumford and Thurston,
it was proven that, if $S$ is a compact oriented
surface of genus $g$ and $P$ is nonempty, then
there is a homeomorphism
between the Teichm\"uller space $\T(S,P)$
and a certain subset $|\Ao(S,P)|$ of the
realization of the arc complex $A(S,P)$,
which is equivariant with respect to the action of
the modular group $\G_{S,P}$.
This gives a homeomorphism between $\M_{g,P}\times\D_P$
and the orbicellular complex $|\Ao(S,P)|/\G_{S,P}$,
which can be equivalently expressed in terms
of the ribbon graph complex $\Mcomb_{g,P}$.

This combinatorial description led to
impressive results about the topology of $\M_{g,P}$:
among the others,
the computation of the virtual cohomological dimension
of $\G_{S,P}$ by Harer (see \cite{harer:virtual}),
the computation of the virtual Euler characteristic
of $\G_{S,P}$ by Harer-Zagier and Penner (see \cite{harer:euler} and
\cite{penner:euler}), computations and estimates of Weil-Petersson
volumes (see for instance
\cite{penner:volumes}, \cite{kmz}
and \cite{grushevsky:wp})
and Witten's conjecture by Kontsevich
(see \cite{kontsevich:intersection}).

There are (at least) two different ways to define
the above homeomorphism. While combinatorially equivalent,
they are geometrically very different.
The first way uses the hyperbolic geometry of noncompact
Riemann surfaces with finite volume
and is due to Penner (see \cite{penner:hyperbolic}) and
Bowditch-Epstein (see \cite{bowditch-epstein:triangulations}).

The second way uses results of Jenkins (see \cite{jenkins:57})
and Strebel (see \cite{strebel:67})
on meromorphic quadratic differentials, in particular
results of existence and uniqueness for differentials with
closed trajectories, and is due to Harer, Mumford
and Thurston (see \cite{harer:virtual}).
We refer to \cite{strebel:quadratic} for a detailed treatment
of this subject, but see also \cite{hubbard-masur:foliations}.

The conjecture of Witten and Kontsevich
(see \cite{kontsevich:intersection}),
which can be formulated in both the
combinatorial descriptions equivalently, 
says that the $W$ cycles, which are defined
as the locus of ribbon graphs with
assigned valencies of their vertices,
are Poincar\'e dual to
the tautological classes $\k$
on the moduli space.
In fact, these subcomplexes
$W$ determine homology classes with noncompact support
on $\M_{g,P}$;
so we naturally obtain cohomology classes
with rational coefficients by Poincar\'e duality,
because $\M_{g,P}$ is an orbifold.
More precisely, the conjecture says that
a multiple of $\k_r$ is Poincar\'e dual to
$W_{2r+3}$, whose support is the locus of ribbon graphs
with one vertex of valency at least $2r+3$.

More generally, let $m_*=(m_{-1},m_0,m_1,\dots)$ be a sequence of nonnegative
integers such that $\sum_{i\geq -1}(2i+1)m_i=4g-4+2n$
and let $\Mcomb_{m_*,P}\subset\Mcomb_{g,P}$ be the orbicellular
subcomplex of ribbon graphs whose top-dimensional orbicells are
parametrized by ribbon graphs with $m_i$ vertices of valency
$2i+3$. Notice, by the way, that $\Mcomb_{m_*,P}\cong\Mcomb_{g,P}$ if
$m_*=(0,4g-4+2n,0,0,\dots)$.
For every $l=(l_{p_1},\dots,l_{p_n}) \in \R_+^P$
denote by $\Mcomb_{m_*,P}(l)$
the subset of graphs in $\Mcomb_{m_*,P}$ such that
the $p_i$-th hole has perimeter $2 l_{p_i}$. Remark that
$\Mcomb_{m_*,P}(\R_+^P)$ is homeomorphic to
$\Mcomb_{m_*,P}(l)\times\R_+^P$ for every $l\in\R_+^P$.

Kontsevich \cite{kontsevich:intersection}
and Penner \cite{penner:poincare}
proved that for every $l\in\R_+^P$ the orbicomplex
$\Mcomb_{m_*,P}(l)$ has an orientation, and so
the classifying map
$\Mcomb_{m_*,P}(l) \rar \M_{g,P}$ defines
a homology class with noncompact support
$W_{m_*,P}$ on $\M_{g,P}$, that does not
depend on the choice of $l\in \R_+^P$, and
which will be called {\it combinatorial class}.
Moreover, Kontsevich introduced combinatorial ``compactifications''
$\Mbarcomb_{m_*,P}$ which still have orientations
and (in the case $m_{-1}=0$)
embed as subcomplexes $\Mbarcomb_{m_*,P}\hra\Mbarcomb_{g,P}$;
thus, they define cycles
$\Wbar_{m_*,P}(\R_+^P)\in H^{BM}_*(\Mbarcomb_{g,P}(\R_+^P);\Q)$
in Borel-Moore homology and ordinary cycles
$\Wbar_{m_*,P}(l)\in H_*(\Mbarcomb_{g,P}(l);\Q)$
for every $l\in\R_+^P$.
In fact, $\Mbarcomb_{g,P}(l)$ is homeomorphic
to a quotient $\Mbar'_{g,P}$
of the Deligne-Mumford compactification $\Mbar_{g,P}$
for all $l\in\R_+^P$ and the class $\Wbar_{m_*,P}(l)$
on $\Mbar'_{g,P}$ does not depend on $l$.
However, while $\Mbar_{g,P}$ is an orbifold and
so virtual Poincar\'e duality holds and allows us
to write equalities that mix rational homology
and cohomology classes,
$\Mbar'_{g,P}$ might have ugly singularities,
so it is not so easy
to get a cohomology class on $\Mbar_{g,P}$.
In any case, we will always consider homology
and cohomology groups with rational coefficients
in what follows.
\begin{conjecture}[Witten-Kontsevich \cite{kontsevich:intersection}]
For every $m_*=(0,m_0,m_1,\dots)$ such that
$\sum_i (2i+1)m_i=4g-4+2|P|$, the class
$W_{m_*,P} \in H_*^{BM}(\M_{g,P})$ is Poincar\'e dual
to a polynomial in the $\k$ classes
$f_{m_*}(\k_1,\k_2,\dots) \in H^*(\M_{g,P})$.
\end{conjecture}

First results towards a proof of the conjecture
were obtained by Penner
\cite{penner:poincare}: using a result of Wolpert
\cite{wolpert:homology} on the Weil-Petersson metric,
he could deal with the case $W_5=12\k_1$.

The approach of Arbarello and Cornalba \cite{arbarello-cornalba:combinatorial} 
passes through matrix models,
Di Francesco-Itzykson-Zuber's theorem \cite{dfiz:polynomial}
and Kontsevich's compactification $\Mbarcomb_{g,P}$
and led to stronger results.
In fact,
they found a way to compute in principle all the
$W_{m_*,P}$ in terms of the kappa classes 
and reported their results in
lower codimensions, giving a strong evidence to the
conjecture.
For example, they discovered that on $\M_{g,P}$
the cycle $W_{(0,m_0,3,0,\dots),P}$ is dual
to $288\k_1^3-4176\k_1\k_2+20736\k_3$.
Looking at a number of results such as the previous one,
they refined the conjecture as follows.
\begin{conjecture}[\cite{arbarello-cornalba:combinatorial}]
Consider the algebra of polynomials $\Q[t]:=\Q[t_1,t_2,\dots]$ where
each $t_i$ has degree $1$. Then for every $m_*=(0,m_0,m_1,\dots)$
there exists a polynomial $f_{m_*}\in \Q[t]$
of degree $\sum_{i\geq 1}m_i$ such that
\[
W^*_{m_*,P}=f_{m_*}(\k_1,\k_2,\dots)\in H^*(\M_{g,P})
\]
where $W^*_{m_*,P}$ is the Poincar\'e
dual of $W_{m_*,P}$. Moreover $f_{m_*}$ looks like
\[
f_{m_*}(t)=
\prod_{i\geq 1}\frac{(2^{i+1}(2i+1)!!)^{m_i}}{m_i!} t_i^{m_i}
+ \text{(terms of lower degree).}
\]
\end{conjecture}
In any event, the meaning of the other coefficients of $f_{m_*}$
was still obscure.

Really, they compared the combinatorial classes
and the kappa classes as functionals
on the algebra generated by the $\psi$ classes, which are defined both on
$\Mbar_{g,P}$ and on $\Mbarcomb_{g,P}(l)$.
In this way, they were able to compute the difference
$W^*_{m_*,P}-f_{m_*}(\k)$ in some concrete cases
up to some minor uncertainty.

In this paper we give an affirmative answer
to the previous conjecture and we exhibit a formula
that permits one to compute all the polynomials $f_{m_*}$
inductively on their degree.

Quite recently, K. Igusa \cite{igusa:mmm} \cite{igusa:kontsevich}
and K. Igusa-M. Kleber \cite{igusa:trees} have proven very similar
results by different methods.
\end{subsection}
\begin{subsection}{Contents of the paper}
We begin with some introductory material on the Teichm\"uller space
and the moduli space of curves, the combinatorial description
and the tautological classes.

Next, we introduce the combinatorial classes $\Wbar_{m_*,P}$
on $\Mbarcomb_{g,P}$ and some generalized classes
$\Wbar_{m_*,\rho,P}\subset \Wbar_{m_*,P\cup Q}$
on $\Mbarcomb_{g,P\cup Q}$ depending on a
$\rho:Q\lra \Z_{\geq -1}$, which are
defined prescribing that every $q\in Q$ marks a vertex of
valency at least $2\rho(q)+3$.
 
After that, we lift the cycles $\Wbar_{m_*,P}(l)$
on $\Mbarcomb_{g,P}$
to cycles on $\Mbar_{g,P}$
(resp. the cycles $\Wbar_{m_*,\rho,P}$
from $\Mbarcomb_{g,P\cup Q}$ to $\Mbar_{g,P\cup Q}$).
However,
in order to do that, we have to pay
a price. In fact, we do not obtain
absolute homology classes of $\Mbar_{g,P}$
(resp. $\Mbar_{g,P\cup Q}$)
but only homology classes relative to
a certain algebraic closed subset
$\S_{g,P}\subset\Mbar_{g,P}$
(resp. $\S^Q_{g,P}\subset\Mbar_{g,P\cup Q}$),
which in any case is contained in the boundary.

We remark that the $\Wbar_{m_*,P}$ classes can be obtained
by pushing $\Wbar_{m_*,\rho,P}$ forward
via a combinatorial analogue of the forgetful map
$\pi_Q:\Mbar_{g,P\cup Q} \rar \Mbar_{g,P}$.
Hence, we prove first that
$\Wbar_{m_*,\rho,P}$ is
Poincar\'e dual to a polynomial in the $\psi$
classes, if all nontrivalent vertices of the general
ribbon graph of $\Wbar_{m_*,\rho,P}$ are marked
by $Q$.
Then we obtain our result for any combinatorial class,
because we know how $\psi$ classes behave under forgetful
morphisms (Faber's formula).

The simplest case is the class $W^q_{2r+3}$ supported on the
subcomplex of $P\cup\{q\}$-marked ribbon graphs in which
$q$ marks a vertex of valency at least $2r+3$. We prove that
$W^q_{2r+3}$ is Poincar\'e dual to $2^{r+1}(2r+1)!!\psi_q^{r+1}$
on $\M_{g,P\cup\{q\}}$; with a little more care we obtain also
that $W_{2r+3}$ is Poincar\'e dual to $2^{r+1}(2r+1)!!\k_r$
on $\M_{g,P}$, essentially because
$(\pi_q)_*(\psi_q^{r+1})=\k_r$.

Very roughly, we want to verify the above equality
$W^q_{2r+3}=2^{r+1}(2r+1)!!(\psi_q^{r+1})^*$
(where $(\psi_q^{r+1})^*$ is the Poincar\'e dual of $\psi_q^{r+1}$)
by coupling both sides with a closed PL differential
form $\eta$ on $\Mbarcomb_{g,P\cup\{q\}}$.

In this computation, we exploit a nice and explicit
PL differential form $\ol{\omega}_q$ on $\Mbarcomb_{g,P\cup\{q\}}$
(found by Kontsevich), whose class pulls back to $\psi_q$
on $\Mbar_{g,P\cup \{q\}}$. 

Hence, the problem translates to showing that,
for some perimeter lengths $\tilde{l}$ and $\tilde{l}'$,
the equality
\[
\int_{\Mbarcomb_{g,P\cup\{q\}}(\tilde{l})}\ol{\omega}_q^{r+1}\wedge\eta=
2^{r+1}(2r+1)!!\int_{\Wbar^q_{2r+3}(\tilde{l}')}\eta+
\int_{\ol{N}^q(\tilde{l}')}\eta
\]
holds, where $\ol{N}^q(\tilde{l}')$ is a certain cycle sitting in the boundary.
The main problem is that the combinatorial class $\Wbar^q_{2r+3}$
is defined in the slices $\Mbarcomb_{g,P\cup\{q\}}(\tilde{l}')$
with $\tilde{l}'_q=0$, while the differential form $\ol{\omega}_q$
is defined only on the slices $\Mbarcomb_{g,P\cup\{q\}}(\tilde{l})$
with $\tilde{l}_q>0$.

The key idea to overcome this difficulty is to define
a deformation retraction $\H^q_0$ that
shrinks the $q$-th hole and so provides a bridge between
the region $\{\tilde{l}_q>0\}$ and the slice $\{\tilde{l}_q=0\}$.
In other words, $\H^q_0$ makes it possible to
recover the combinatorial class $\Wbar^q_{2r+3}$
as push-forward
of $\ol{\omega}_q^{r+1}$.
We remark that this $\H^q_0$ retracts some cells
sitting in the interior of $\Mbarcomb_{g,P\cup\{q\}}$
to the boundary. Sometimes the curious behaviour of this map
gives rise to some technical problems; nevertheless,
it is a key tool in the proofs and it has the merit
of being easily visualizable.

Once we have our retraction $\H^q_0$,
we discover that, by simple reasons of degree,
the restriction of
the differential form $\ol{\omega}_q^{r+1}\wedge(\H^q_0)^*\eta$
on $\Mbarcomb_{g,P\cup\{q\}}(\tilde{l})$
is supported on the smallest
subcomplex $\ol{Y}_{2r+3}(\tilde{l})$
that contains all the cells parametrized by
ordinary ribbon graphs whose $q$-th hole is bordered by $2r+3$ edges.

Next, we dissect $\ol{Y}_{2r+3}$ into subcomplexes $\ol{Y}^i_{2r+3}$
according to the topology of the $q$-th hole.
In this way, the restriction
of $\H^q_0$ to each $\ol{Y}^i_{2r+3}$ is 
generically a fibration whose fibers $F^i$
are simplicial complexes of dimension $2r+3$.
Hence
\[
\int_{\Mbarcomb_{g,P\cup\{q\}}(\tilde{l})}\ol{\omega}_q^{r+1}
\wedge(\H^q_0)^*\eta=
\sum_i
\int_{\H^q_0(\ol{Y}^i_{2r+3}(\tilde{l}'))}\eta
\int_{F^i\cap\{\tilde{l}_q=\e\}} \ol{\omega}_q^{r+1}
\]
and we get the result analyzing $\H^q_0(\ol{Y}^i_{2r+3})$ and
computing the integral on the fibers.
We underline that,
once the spaces and the machinery are set up, the actual calculation
is really straightforward.

For example, the class $\Wbar^q_{2r+3}$ arises as image of top-dimensional
simplices when the hole $q$ is contractible,
i.e. no edge borders the hole $q$ from both sides.
In this case, the fiber is just one simplex (provided $r\geq 1$) and
the integral on the fiber is exactly $\frac{(r+1)!}{(2r+2)!}$.
Other top-dimensional simplices give rise to boundary terms.

Hence, in the case of combinatorial
classes with only one nontrivalent vertex,
the main result (in a slightly simplified version)
is the following.
\begin{mytheorem} \label{th:first}
For any $g$ and $n\geq 1$, the equality
\[
\Wbar^q_{2r+3}
=\frac{(2r+2)!}{(r+1)!}(\psi_q^{r+1})^*
\]
holds in $H_{6g-6+2n-2r}(\Mbar_{g,P\cup\{q\}},
\pa\M_{g,P\cup\{q\}})$ for every $r\geq -1$.
Moreover, for $r\geq 1$ the equality
\[
\Wbar_{2r+3}
=2^{r+1}(2r+1)!!\,\k^*_r
\]
holds in $H_{6g-6+2n-2r}(\Mbar_{g,P},\pa\M_{g,P})$.
\end{mytheorem}
As an example, we have the following corollary which was
already proven by Arbarello and Cornalba in a very different manner
(see \cite{arbarello-cornalba:combinatorial}).
\begin{mycorollary}
For every $g$ and $|P|=n\geq 1$ such that $2g-2+n>0$
the following equalities hold
\[
\begin{array}{ll}
\displaystyle{
\Wbar^q_5+\d^q_{irr}+
\sum_{\substack{g',I\neq \emptyset,P}}  \d^q_{g',I}
=12(\psi_q^2)^* }
& \quad\text{in $H_{6g-8+2n}(\Mbar_{g,P\cup\{q\}},\S^q_{g,P})$} \\
\displaystyle{
\Wbar_5+\d_{irr}+
\sum_{\substack{g',I\neq\emptyset,P}}  \d_{g',I}
=12 \k^*_1 }
& \quad\text{in $H_{6g-8+2n}(\Mbar_{g,P},\S_{g,P})$}
\end{array}
\]
where $\d^q_{g',I}$ is the image of the morphism
\[
\Mbar_{g',I\cup\{p'\}} \times\Mbar_{0,\{q,q',q''\}}
\times\Mbar_{g-g',I^c\cup\{p''\}} \rar \Mbar_{g,P\cup\{q\}}
\]
that glues $p'$ with $q'$ and $p''$ with $q''$
(analogously for $\d^q_{irr}$).
\end{mycorollary}

Next, we pass to examine the case of
a general combinatorial class $\Wbar_{m_*,P}$.
As explained before, we recover them by pushing some
$\Wbar_{m_*,\rho,P}$ forward through a combinatorial
version of the forgetful morphism
$\pi_Q:\Mbar_{g,P\cup Q}\lra\Mbar_{g,P}$.
The proof works more or less as before,
but a new phenomenon occurs when we shrink
via $\H_0^Q$ all the $Q$-marked holes. 
In fact, the cycles in the image of $\H_0^Q$
do not depend only on the shape of a single
hole $q_i\in Q$, but also on the configuration
of the $Q$-marked holes in the ribbon graph.
For instance, in the case $Q=\{q_1,q_2\}$,
the map $\H^{q_1,q_2}_0$ produces two
different cycles if the two holes $q_1$ and $q_2$
``touch'' each other (i.e. they have at least
one vertex in common) or if they are detached.
Hence, a careful combinatorial analysis is needed
in order to prove the theorem for a general $\Wbar_{m_*,\rho,P}$.

The notations and the
results about the classes $\Wbar_{m_*,\rho,P}$ are quite
heavy, so here we content ourselves
to state the theorem in the
simpler case of $\Wbar_{m_*,P}$ and to defer to
Section \ref{sec:second} for more refined results.

Choose $\rho: Q \rar \N_+$ and $m_*=(0,m_0,m_1,\dots)$
such that $m_0\geq 0$, $m_i=\rho^{-1}(i)$ for $i>0$
and $\sum_{i\geq 0} (2i+1) m_i=4g-4+2n$.
Let $\mathfrak{P}_{Q}$ be the set of partitions of $Q$ and
for all $\mu \subset Q$ define
$\rho_{\mu}=\sum_{q_i\in \mu} \rho(q_i)$.
\begin{mytheorem}[simplified version] \label{th:second}
For any $g\geq 0$ and $P\neq\emptyset$ such that $2g-2+|P|>0$,
the following equality holds
\begin{equation*}
 \sum_{M \in \mathfrak{P}_{Q}}
 C(\rho,M) \Wbar_{m_*(M),P}
=(\pi_Q)_*\left[\prod_{q\in Q}
2^{\rho(q)+1} (2\rho(q)+1)!!\,
(\psi_q^{\rho(q)+1})^*\right]
\end{equation*}
in $H_*(\Mbar_{g,P},\pa\M_{g,P})$,
where
\begin{equation*}
m_i(M)=
\begin{cases}
m_0 & \text{if $i=0$} \\
|\{\mu\in M| \rho_{\mu}=i\}| & \text{if $i\neq 0$}
\end{cases}
\end{equation*}
and $C(\rho,M)$ is an explicit integer coefficient that depends
only on $\rho$ and $M$.
\end{mytheorem}
The theorem gives a recipe, which works inductively on $|Q|$, to calculate all
the coefficients of $f_{m_*}$. As an example we have the following.
\begin{mycorollary}
For every $g\geq 0$ and $P\neq\emptyset$ such that $2g-2+|P|>0$
and for every $a,b\geq 1$, the
following identity
\begin{multline*}
2^{\d_{a,b}}\Wbar_{2a+3,2b+3}=
2^{a+b+2}(2a+1)!!(2b+1)!!(\k^*_a\k^*_b+\k^*_{a+b}) \\
-2^{a+b+1}(2a+2b+3)!!\k^*_{a+b}
\end{multline*}
holds in $H_*(\Mbar_{g,P}\,\pa\M_{g,P})$.
\end{mycorollary}
\end{subsection}
\begin{subsection}{Plan of the paper}
In the first section we recall a few facts on the Teichm\"uller
space and the moduli space of curves and we introduce
Kontsevich's compactification
$\Mbartri_{g,P}$ as a quotient of the product
$\Mbar_{g,P}\times\D_P$ of Deligne-Mumford compactification
by $\D_P=\{l\in\R_{\geq 0}^P|\,\sum_{p\in P}l_p=1\}$.
Moreover, we describe the fiber of the quotient map
$\xi:\Mbar_{g,P}\times\D_P\rar\Mbartri_{g,P}$.

In the second section we introduce the arc complex $A(S,P)$
(following the presentation of Looijenga \cite{looijenga:cellular})
and the ribbon graph complex
$\Mbarcomb_{g,P}$ (following Kontsevich \cite{kontsevich:intersection})
and we recall the equivalence of these two constructions
(already present in \cite{looijenga:cellular})
and Kontsevich's isomorphism $\Mbarcomb_{g,P}\cong\Mbartri_{g,P}\times\R_+$.
We follow the approach by means of quadratic differentials.
However, since our tools and techniques are essentially combinatorial,
there is an identical argument in the parallel setting in hyperbolic
geometry.

In the third section we define the tautological classes
and we recall the string and dilaton equations and,
more generally, Faber's formula which govern the push-forward
of $\psi$ classes via the forgetful maps.

In the fourth section we introduce $\Mbarcomb_{m_*,P}$
and their generalizations $\Mbarcomb_{m_*,\rho,P}$.
Moreover, we describe Kontsevich's combinatorial
representatives $\ol{\omega}$ of the $\psi$ classes and we
recall how their (weighted) sum gives a sort of symplectic
form and so an orientation form
on these complexes.

In the fifth section we develop the main technical tools,
namely the retraction
$\H^q_0$ in the simplest case
and the combinatorial forgetful
map $\pi_q^{comb}$, which are used in
the proof of Theorem \ref{th:first}
contained in Section 6.

Parallely, the seventh section extends the retraction
and the combinatorial forgetful map to the case
of many marked vertices; while, in the eighth section,
we introduce combinatorial classes with rational tails and
we prove the full version of Theorem \ref{th:second}.

Finally, in the Appendix 
we describe Looijenga's modification $\Ahat(S,P)$
of the arc complex
(see \cite{looijenga:cellular}),
that maps to the Deligne-Mumford compactification
in a $\G_{S,P}$-equivariant way,
and we discuss how it might give a canonical way
to lift the combinatorial classes to $\Mbar_{g,P}$
without ambiguities.
\end{subsection}
\begin{subsection}{Acknowledgements}
I warmly thank my advisor Enrico Arbarello for having
introduced me to the subject and for his constant support
spanning many years.
I also wish to thank Robert Penner for useful discussions,
and Domenico Fiorenza and Riccardo Murri
for valuable suggestions and remarks.
Finally, I want to thank the referee for carefully
reading the paper and for several
comments and improvements.
\end{subsection}
\end{section}

\begin{section}{Moduli spaces of curves and compactifications}
Let $S$ be a compact connected oriented surface
of genus $g$ and let $P \hra S$
be an injection of $n$ points such that $\chi(S\setminus P)=2-2g-n<0$.

\begin{definition}
A family of $P$-pointed surfaces is
a couple $(\pi,s)$ where $\pi: \Cc \rar B$ is a proper
differentiable submersion whose fibers are oriented
connected surfaces and $\{s_p : B \rar \Cc | p \in P \}$
is a collection of disjoint sections.
An {\it $(S,P)$-marking} is an equivalence class
of oriented diffeomorphisms
$f: S \times B \arr{\sim}{\lra} \Cc$
that commute with
the projections onto $B$ and such that
$f(p,b)=s_p(b)$ for every $p \in P$.
Two markings $f \sim \tilde{f}$ are {\it equivalent} if and only if
\[
\tilde{f}^{-1} \circ f:(S,P)\times B \rar (S,P)\times B
\]
is a vertically (i.e. over $B$) isotopic to the identity
relatively to $P$.
\end{definition}

A {\it conformal structure} is an atlas such
that the changes of coordinates are differentiable and
preserve the angles. There is an obvious
bijection between conformal structures and complex
structures and between conformal structures and Riemannian
metrics up to multiplication by a positive function.

\begin{remark}
By uniformization, the universal covering of
a Riemann surface $S\setminus P$ with negative Euler characteristic
is isomorphic to the unit disk $\D=\{z\in \C\,|\,|z|<1\}$.
Thus, the Poincar\'e metric $\frac{4dz\, d\bar{z}}{(1-|z|^2)^2}$
on $\D$ descends to a complete hyperbolic metric
$S\setminus P$ of finite volume with cusps at $P$, which is unique
because the analytic automorphisms of $\D$ are isometries.
Vice versa, given a complete hyperbolic metric of finite volume
on $S\setminus P$ one can
associate a complex structure in a canonical way.
As we are taking a conformal approach, we will not
pursue the hyperbolic point of view in what follows.
\end{remark}

\begin{definition}
Let $(\pi,s)$ be a family of $P$-pointed surfaces.
A {\it conformal structure} on $(\pi,s)$ is
a differentiable atlas of $\Cc$ which endows $\Cc_b\setminus\cup s_p(b)$
with a conformal structure for all $b \in B$.
\end{definition}

We say that two marked families $(\Cc,f)$ and $(\Cc',f')$ of
$P$-pointed surfaces with conformal structure are
{\it isomorphic} if there is a diffeomorphism $t:\Cc \arr{\sim}{\lra} \Cc'$
such that $t \circ f'=f$ and the restriction to each fiber
$t_b:\Cc_b \arr{\sim}{\lra} \Cc'_b$ is conformal outside the sections.

The Teichm\"uller functor
$\mathfrak{T}_{S,P}: \mathrm{(Top. \,\, Spaces)} \lra \mathrm{(Sets)}$
associates to every manifold $B$ the set of isomorphism classes
of $(S,P)$-marked families of $P$-pointed surfaces
over $B$ with conformal structure.
It is represented by a complex
manifold $\T_{S,P}$ of (complex) dimension $3g-3+n$,
which is diffeomorphic to a ball.
Except in the case $(g,n)=(0,3)$ it is never compact.

The modular group
$\G_{S,P}:=\mathrm{Diff_+(S,P)}/\mathrm{Diff_0(S,P)}$
of connected components of the space
of oriented diffeomorphisms of $(S,P)$
acts on the $(S,P)$-markings and so on $\Tfun_{S,P}$.
Its quotient is denoted by $\Mfun_{g,P}$ and
classifies smooth families of $P$-pointed Riemann surfaces
of genus $g$ up to isomorphism.
However this functor is not representable,
so the topological quotient
$\T_{S,P}/\G_{S,P}$ is only a
{\it coarse moduli space}.

The functor $\Mfun_{g,P}$ admits a natural extension
$\Mfunbar_{g,P}$ that classifies flat families
of stable $P$-pointed Riemann surfaces of (arithmetic) genus $g$,
where {\it stable} means that the singularities look like
$\{xy=0\}\subset \C^2$ in local analytic coordinates
and that each connected
component of the smooth locus has negative topological Euler
characteristic.
The functor $\Mfunbar_{g,P}$
admits a coarse moduli space $\overline{M}_{g,P}$
which is a normal irreducible projective variety
with quotient singularities and which
contains $M_{g,P}$ as a Zariski-dense open
subset. It can be seen that $\Mfunbar_{g,P}$ is in fact
represented by an orbifold $\Mbar_{g,P}$ which is connected
and compact.

Now we want to introduce a different compactification
$\Mbartri_{g,P}$
due to Kontsevich \cite{kontsevich:intersection}
which will be very useful in what follows.
Really, we slightly modify Kontsevich's construction
as we realize our space as a quotient of $\Mbar_{g,P}\times \D_P$
by an equivalence relation, where
$\D_P:=\{l\in \R_{\geq 0}^P|\,\,\sum_{p\in P}l_p=1\}$.

If $(S,l)$ is an element of $\Mbar_{g,P}\times \D_P$, then
we say that an irreducible component of $S$
is {\it positive} (with respect to $l$) if it contains a point
$p \in P$ such that $l_p>0$.
So we declare
that $(S,l)$ is equivalent to $(S',l')$ if
$l=l'$ and there is a homeomorphism of pointed surfaces
$S \arr{\sim}{\lra} S'$ which is holomorphic on the positive
components of $S$.

\begin{figure}[h]
\resizebox{10cm}{!}{\input{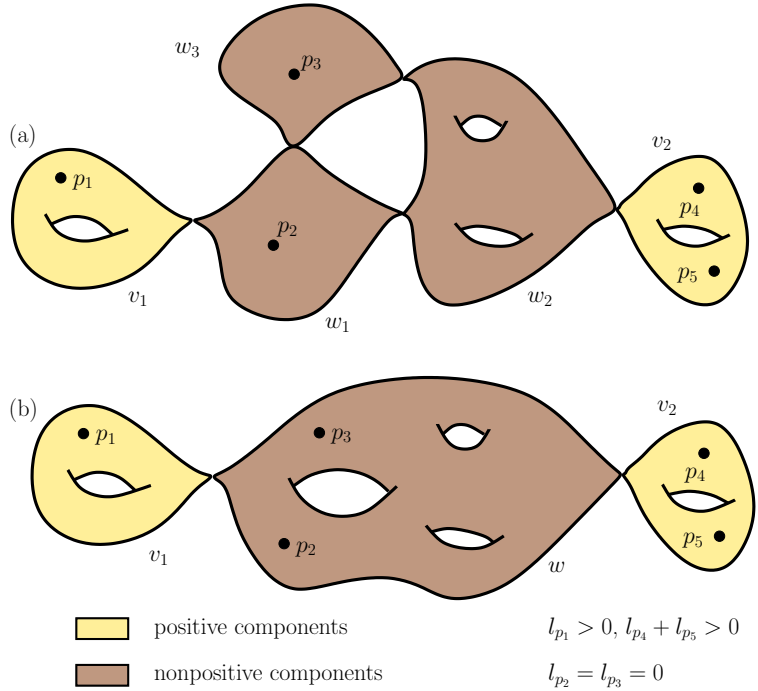}}
\caption{Two surfaces with different dual graphs}
\label{fig:curve}
\end{figure}
As this relation would not give back a Hausdorff space,
we consider its closure, which can be described as follows.

We attach to every $S$ its dual graph $\g_S$, whose vertices are
irreducibile components and whose edges are nodes of $S$.
Moreover, every vertex $v$ is labelled by a couple $(g_v,P_v)$,
where $g$ is the geometric genus of $v$ and $P_v \subset P$
is the set of marked points lying on $v$.
Given $(S,l)$ as before, consider the following two moves:
\begin{enumerate}
\item
if two nonpositive vertices $v$ and $v'$ are joined by
an edge $e$, then we can build a new graph discarding $e$,
melding $v$ and $v'$ together to obtain a new vertex $w$
which we label with $(g_w,P_w):=(g_v+g_{v'},P_v\cup P_{v'})$
\item
if a nonpositive vertex $v$ has a loop $e$, we can
make a new graph discarding $e$ and labelling $v$
with $(g_v+1,P_v)$.
\end{enumerate}
Applying these moves to $\g_S$ iteratively until the process ends,
we are given back a {\it reduced dual graph}
$\g_{S,l}^{red}$. \label{sec:reduced-dual-graph}
Call $V_0(S,l)$ the subset of vertices $v$ of $\g_{S,l}^{red}$
such that $l_p=0$ for every $p \in P_v$
and call $V_+(S,l)$ the subset of positive components of $S$.
\begin{figure}[h]
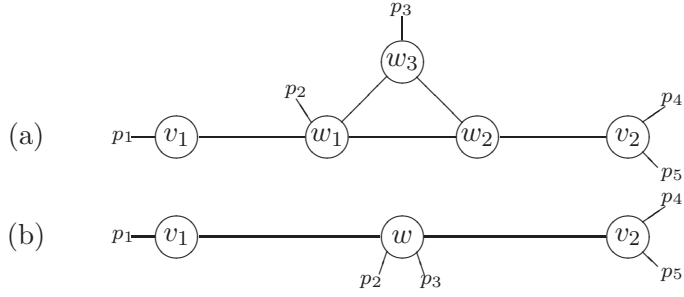

\begin{align*}
\xy
,(-20,0)*{\text{(a)}}
\ar@{-}(3,0);(17,0)
\ar@{-}(23,0);(37,0)
\ar@{-}(22,2);(28,8)
\ar@{-}(38,2);(32,8)
\ar@{-}(43,0);(57,0)
\ar@{-}(-3,0);(-6,0) 
\ar@{-}(18,2);(16,5) 
\ar@{-}(30,13);(30,16) 
\ar@{-}(62,2);(65,4) 
\ar@{-}(64,-4);(66,-5) 
,(-7,0)*{_{p_1}}
,(16,6)*{_{p_2}}
,(30,17)*{_{p_3}}
,(66,5)*{_{p_4}}
,(66,-5)*{_{p_5}}
,(0,0)*+{v_1}*\cir<9pt>{}
,(20,0)*+{w_1}*\cir<9pt>{}
,(30,10)*+{w_3}*\cir<9pt>{}
,(40,0)*+{w_2}*\cir<9pt>{}
,(60,0)*+{v_2}*\cir<9pt>{}
\endxy
\\
\xy
,(-20,0)*{\text{(b)}}
\ar@{-}(3,0);(27,0)
\ar@{-}(33,0);(57,0)
\ar@{-}(-3,0);(-6,0) 
\ar@{-}(28,-2);(27,-5) 
\ar@{-}(32,-2);(33,-5) 
\ar@{-}(62,2);(65,4) 
\ar@{-}(64,-4);(66,-5) 
,(-7,0)*{_{p_1}}
,(26,-6)*{_{p_2}}
,(34,-6)*{_{p_3}}
,(66,5)*{_{p_4}}
,(66,-5)*{_{p_5}}
,(0,0)*+{v_1}*\cir<9pt>{}
,(30,0)*+{w}*\cir<9pt>{}
,(60,0)*+{v_2}*\cir<9pt>{}
\endxy
\end{align*}
\caption{Dual graphs of surfaces in Fig.~\ref{fig:curve}}
\label{fig:curve_dual_graph}
\end{figure}
\begin{example}
Consider the nodal surfaces of genus $5$ as in Fig.~\ref{fig:curve}
and let $l$ be given in such a way that $l_{p_1}>0$ and
$l_{p_4}+l_{p_5}>0$, but $l_{p_2}=l_{p_3}=0$.
Hence, components $w_1$, $w_2$ and $w_3$ in case (a) and $w$ in case (b)
are nonpositive, while $v_1$ and $v_2$ are positive.
Their associated dual graphs are obviously different
(see Fig.~\ref{fig:curve_dual_graph}); but in case (b)
it is reduced, while in case (a) it is not.
In fact, it is easy to check that both surfaces have the same
reduced dual graph (which of course coincide with the dual
graph of the surface (b)).
\end{example}
For every couple $(S,l)$ denote by $\bar{S}$ the
quotient of $S$ obtained collapsing every nonpositive
component to a point.
We say that $(S,l)$ and $(S',l')$ are equivalent
if $l=l'$ and
there exist a homeomorphism
$\bar{f}:\bar{S} \arr{\sim}{\lra} \bar{S}'$
and an isomorphism
$f^{red}: \g_{S,l}^{red} \arr{\sim}{\lra} \g_{S',l'}^{red}$
of reduced dual graphs such that
$\bar{f}$ and $f^{red}$ are compatible
and the restriction of $\bar{f}$ to each component is
holomorphic.

Call $\Mbartri_{g,P}:=\Mbar_{g,P}\times\D_P/\sim$ the quotient
and $\xi: \Mbar_{g,P}\times \D_P \rar
\Mbartri_{g,P}$ the natural projection.
Remark that $\Mbartri_{g,P}$ is compact
and that $\xi$ commutes with the projection onto $\D_P$.
For every $l$ in $\D_P$ we will
denote by $\Mbartri_{g,P}(l)$ the subset of points
of the type $[S,l]$ and we will write $\Mbartri_{g,P}(L)$
for $\cup_{l\in L}\Mbartri_{g,P}(l)$ where $L\subset \D_P$.
Then it is easy to see that $\Mbartri_{g,P}(\Dc_P)$
is in fact homeomorphic to a product $\Mbartri_{g,P}(l)\times\Dc_P$
for any given $l \in \Dc_P$.

Finally notice that the fibers of $\xi$ are isomorphic to moduli spaces.
More precisely, consider a point $[S,l]$ of $\Mbartri_{g,P}$
and call $Q_v$ the subset of half-edges
of $\g_{S,l}^{red}$ incident on the vertex $v$
for every $v \in V_0(S,l)$.
Thus $\xi^{-1}([S,l])\cong \prod_{v\in V_0(S,l)}
\Mbar_{g_v,P_v\cup Q_v}$.
\end{section}
%
%
\begin{section}{Combinatorial description}
\label{sec:combinatorial}
Fix a compact connected oriented surface $S$ of genus $g$
and an injection $P:=\{p_1,\dots,p_n\} \hra S$ with $n>0$.
\begin{subsection}{The arc complex}\label{subsec:arc}
Let $\mathcal A$ be the set of isotopy classes
relative to $P$ of unoriented loops or arcs
embedded in $S$ that intersect $P$ exactly in the extremal point(s).
The {\it arc complex} is the abstract simplicial complex $A(S,P)$ whose
$k$-simplices are subsets $\ua=\{ \a_0, \dots, \a_k \}$ of $\mathcal A$
that are representable by a system of $k+1$ arcs
and loops intersecting only in $P$. We will denote
its geometric realization by $|A|$.

A simplex $\ua=\{\a_0,\dots,\a_k\}$ of $A$
is called {\it proper} if its star is finite,
or equivalently if $S \setminus \cup_{i=0}^k \a_i$ is a disjoint union of
open disks, each one containing at most one point of $P$.
The subset $\A8 \subset A$ of improper simplices is
a subcomplex; we denote $\Ao:=A\setminus \A8$ the
subset of proper ones and by $|\Ao|$ its ``geometric realization''
$|A|\setminus|\A8|$.

We will associate a marked ribbon graph $G_{\ua}$
to every proper simplex $\ua$ in a natural way
and a metric on $G_{\ua}$ to every internal point of $|\ua|$.
Let us fix some notation first.

\begin{definition}
An {\it (ordinary) ribbon graph} $G$ is a triple
$(X(G),\s_0,\s_1)$ such that
$X(G)$ is a nonempty finite set, $\s_0$ is a permutation of $X(G)$
and $\s_1$ is a fixed-point-free involution of $X(G)$. 
Let denote by $X_i(G)$ the set of orbits in $X(G)$ with respect
to the action of $\s_i$ for $i=0,1$.
\end{definition}
\begin{figure}[h]
\resizebox{8cm}{!}{\input{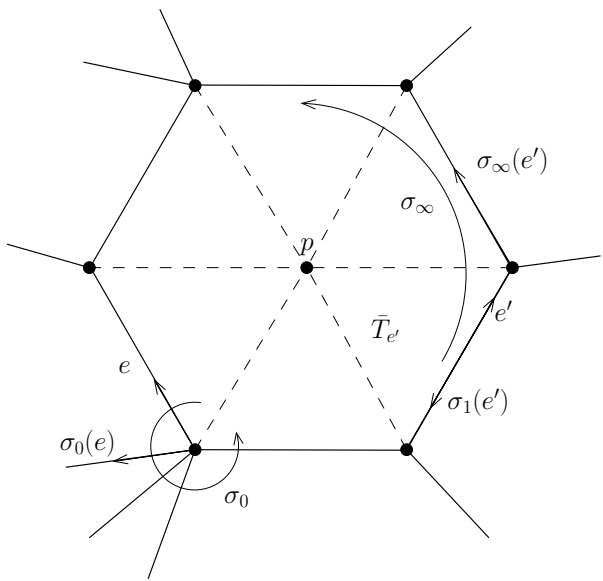}}
\caption{Geometric representation of a ribbon graph}
\label{fig:sigma}
\end{figure}
Remark that this definition is equivalent to the intuitive
one given in terms of a graph plus a cyclic ordering
of the half-edges incident on each vertex (see Fig.~\ref{fig:sigma}).
In fact we shall regard $X(G)$ as the set of oriented edges of $G$,
$X_0(G)$ as the set of vertices and $X_1(G)$ as the set
of unoriented edges. So we can identify $\s_0$ with
the operator that sends an oriented
edge outcoming from a vertex $v$ to the following oriented edge
outcoming from $v$ with respect to a given cyclic order,
and $\s_1$ with the operator that simply reverses the orientation
of the given oriented edge. Consequently, there is a natural
bijection between connected components of the ribbon graph $G$
and orbits in $X(G)$ under the action of the subgroup
$\la \s_0,\s_1 \ra \subset \mathfrak{S}(X(G))$.
Finally we can define $\s_\infty$ requiring that
$\s_\infty \s_1 \s_0 =1$, so that
$X_\infty(G)$ naturally corresponds to the set of holes of $G$
and $\s_\infty$ rotates the edges that border each hole.
To each ribbon graph $(X(G),\s_0,\s_1)$ we can associate
a dual one $G^*:=(X(G),\s_\infty^{-1},\s_1)$ such that $(G^*)^*=G$.
\begin{definition}
A {\it $P$-marking} of $G$ is an injection
$x: P \hra X_0(G) \cup X_\infty(G)$ such that $X_\infty(G)$ is
in the image of $x$.
A {\it realization} of the ribbon graph $G$
into an oriented surface $S'$ is a realization of the graph $|G|$
together with an embedding $|G|\hra S'$ which is
compatible with the cyclic ordering of the half-edges
incident on every vertex of $|G|$.
A realization is {\it proper} if $S'\setminus |G|$ is
a disjoint union of disks and of pointed disks.
\end{definition}
We call $(G,x)$ {\it reduced} if every unmarked vertex has valency
greater than two. In what follows ribbon graphs are intended to
be reduced.

To each proper simplex $\ua=\{\a_0,\dots,\a_k\}$
we can associate a connected ribbon graph $G_{\ua}^*$ simply taking
as $X(G_{\ua}^*)$ the set of oriented versions of $\a_i$'s,
as $\s_1$ the sense-reversing operator and making $\s_0$
rotate edges outcoming from a point $p$ counterclockwise
with respect to the given orientation of $S$.
It is easy to see that $G_{\ua}:=(G_{\ua}^*)^*$ inherits a $P$-marking:
we call it the ``dual'' ribbon graph associated to $\ua$.

Now we show that both $G^*_{\ua}$ and $G_{\ua}$ admit proper
realizations $|G^*_{\ua}|$ and $|G_{\ua}|$ in $S$, which
are canonically defined up to isotopy. This canonicity depends
on the fact that orientation-preserving embeddings of
a disk into an oriented surface are isotopic.

For $G^*_{\ua}$ it is sufficient to choose
explicit representatives for the arcs $\a_0,\dots,\a_k$ and to take
them as edges of the ribbon graph.
For $G_{\ua}$ one proceeds in the following way.
As $\ua$ is a proper simplex, the surface $S\setminus |G^*_{\ua}|$
is a disjoint union of disks and of pointed disks; then, one chooses
as vertices of $|G_{\ua}|$
one point in each unmarked disk and
the marked point in each pointed disk.
Then, for every $\a_i$ one draws on $S$ an edge which joins the vertices
corresponding to the disks separated by $\a_i$ and
which intersects the explicit representative of
$\a_i$ transversely in one point.
It is easy to see that this determines a proper realization
of the $P$-marked graph $G_{\ua}$ in $(S,x)$, which is
canonical up to isotopy.
\begin{remark}
Fix realizations of $G_{\ua}$ and $G^*_{\ua}$.
Both determine a cellular decomposition of $S$
and so a complex of cellular chains in a natural way:
call $C_{\bullet}(S)$ and $C^*_{\bullet}(S)$ these complexes.
The construction described above induces natural
isomorphisms $C_i(S)\cong C^*_{2-i}(S)$ for $i=0,1,2$,
and it is easy to see that they induce Poincar\'e
duality in homology.
\end{remark}
\begin{definition}
A {\it metrized} ribbon graph is a couple $(G,l)$ where
$G$ is a ribbon graph and $l$ is a {\it positive metric}
of total length $1$ on $G$,
i.e. a point of $\Dc_{X_1(G)}$.
\end{definition}
Actually it is clear that a point $a$ of $|\ua|^\circ \subset |\Ao|$
correspond to a positive metric
on $G_{\ua}^*$ and so on $G_{\ua}$.
Moreover, if $\l: |A| \rar \D_P$ is the simplicial map that
sends a vertex $\{\a\}$ of $|A|$ to the barycenter
of the extremal points of the arc $\a$, then the restriction of
$\l$ to a proper simplex
is the {\it circumference function} of the associated ribbon
graph, that is, it sends a metrized ribbon graph $(G,a)$
to the point whose $p$-th coordinate is half the perimeter
of the $p$-marked hole (it is zero in the case that $p$ marks
a vertex).

To each metrized ribbon graph $(G,a)$ we can
canonically associate a Riemann surface obtained
by gluing half-infinite strips
\begin{equation*}
S(G,a):= \left( \coprod_{e \in X(G)} T_e \right) \Big{/ \sim}
\end{equation*}
where $T_e=[0,e(a)] \times [0,\infty]/[0,e(a)]\times\{\infty\}$
and $\sim$ is the equivalence
relation generated by
$T_e \owns (t,0) \sim (e(a)-t,0) \in T_{\s_1(e)}$
and $T_e \owns (e(a),s) \sim (0,s) \in T_{\s_\infty(e)}$.
Call $\bar{T}_e$ the image of $T_e$ under the above identifications
and (if $G$ is $P$-marked)
$\bar{T}_p$ the union of the $\bar{T}_e$'s for all $e\in x(p)$
and notice that the conformal structures on
$\bar{T}_e\setminus(\{\infty\}\cup\{0\}\times\{0\}\cup\{e(a)\}\times\{0\})
\subset\R^2\cong\C$ glue
to give a conformal structure on $S(G,a)$ minus a finite set
and that it induces a well-defined unique complex structure on
the whole $S(G,a)$.
It is clear that a $P$-marking on $G$
descends to a $P$-marking $x': P \hra S(G,a)$
and that $G$ has a natural proper realization in $(S(G,a),x')$.
As $G$ is also properly realized in $(S,x)$ (canonically
up to isotopy), this determines
an isotopy class of diffeomorphisms
$(S,P) \arr{\sim}{\lra} (S(G,a),x'(P))$ and a classifying map
$\Psi:|\Ao(S,P)| \rar \T_{S,P}$ to the Teichm\"uller space.

\begin{theorem}[Harer-Mumford-Thurston]
The map
\[
(\Psi,\l):|\Ao(S,P)| \arr{\sim}{\lra} \T_{S,P}\times \D_P
\]
is a $\G_{S,P}$-equivariant homeomorphism.
Hence the quotient
\[
\Phi:|\Ao(S,P)|/\G_{S,P} \arr{\sim}{\lra} \M_{g,P}\times \D_P
\]
is a homeomorphism too.
\end{theorem}
\begin{proof}
One can construct a tautological family
of Riemann surfaces $\mathcal{C} \lra |\Ao(S,P)|$
whose restriction over a simplex $\ua$ is
real-analytic. So $\Psi$ is continuous by the universal
property of the Teichm\"uller space
and $\Psi|_{\ua}$ is real-analytic for every $\ua$.

Next, we recall that $|\Ao(S,P)|$ can be given a structure
of differentiable manifold compatible with the
piecewise linear one (see \cite{hubbard-masur:foliations}).
Hence, if we prove that $(\Psi,\l)$ is bijective,
then we can conclude that
it is open and so a homeomorphism by invariance of domain.

To prove that $(\Psi,\l)$ is bijective, we notice that
every $S(G,a)$ is canonically 
endowed with a meromorphic quadratic differential $\beta$
which restricts to $(dz)^2$ on each $\bar{T}_e\subset\C$.
This differential has three interesting properties (among others):
\begin{enumerate}
\item
it is holomorphic in $S(G,a)\setminus x'(P)$ and
almost all its horizontal trajectories
(namely, those defined by $\mathrm{Arg}(\beta)=0$) are closed
\item
for every $p\in P$: if $\l_p(a)>0$, then $\beta$ has a double pole
on $x'(p)$ with quadratic residue $-\left(\frac{\l_p(a)}{\pi}\right)^2$;
if $\l_p(a)=0$, then $x'(p)$ lies either a simple pole
of $\beta$ or on a critical trajectory
or is a zero of $\beta$
\item
the zeroes of $\beta$ and the marked points
$\{x'(p)\,|\, \l_p(a)=0\}$ are
the vertices of $|G_{\ua}|\subset S(G,a)$,
the horizontal trajectories
between these points (i.e. the critical trajectories)
are the edges of $|G_{\ua}|$;
hence the critical graph $\mathrm{Crit}(\beta)$ of $\beta$
coincides with $|G_{\ua}|$.
\end{enumerate}
Now we invoke the following celebrated result.
\begin{theorem}[Jenkins-Strebel,
\cite{jenkins:57} \cite{strebel:67}]\label{th:js}
Let $\tilde{S}$ be a compact Riemann surface and
$\tilde{x}:P\hra \tilde{S}$ a nonempty
subset such that $\chi(\tilde{S}\setminus\tilde{x}(P))<0$.
Then, for every nonzero function
$\tilde{h}:P \rar \R_{\geq 0}$, there exists a
unique meromorphic quadratic differential
$\beta(\tilde{S},P,\tilde{h})$
on $\tilde{S}$ with the following properties:
\begin{enumerate}
\item
$\beta(\tilde{S},P,\tilde{h})$
is holomorphic in $\tilde{S}\setminus\tilde{x}(P)$
and it has almost all
closed horizontal trajectories (and so at most double poles)
\item
every closed trajectory of $\beta(\tilde{S},P,\tilde{h})$
is isotopic
inside $\tilde{S}\setminus\tilde{x}(P)$ to a simple loop winding around
the point $\tilde{x}(p)$ for some $p\in P\setminus \tilde{h}^{-1}(0)$
\item
for every $p \in P\setminus \tilde{h}^{-1}(0)$,
the differential $\beta(\tilde{S},P,\tilde{h})$ has
a double pole in $\tilde{x}(p)$ with
quadratic residue $-\left(\frac{\tilde{h}(p)}{\pi}\right)^2$
\item
for every $p \in \tilde{h}^{-1}(0)$, the point $\tilde{x}(p)$
is contained inside the
critical graph of $\beta(\tilde{S},P,\tilde{h})$
(and so $\tilde{x}(p)$ is at worst a simple pole).
\end{enumerate}
\end{theorem}

A consequence of the previous theorem is
that the critical graph of $\beta(\tilde{S},P,\tilde{h})$
is a $P$-marked ribbon graph, properly embedded
in $(\tilde{S},\tilde{x})$. Furthermore it inherits a metric
from the quadratic differential. Hence, we
have produced a well-defined map
\[
\begin{array}{ccc}
\T_{S,P}\times \D_P & \lra & |\Ao(S,P)| \\
([f:(S,P)\rar(\tilde{S},\tilde{x})],h) & \mapsto &
\mathrm{Crit}(\beta(\tilde{S},P,h))
\end{array}
\]
and one can see that it
is the (set-theoretical) inverse of $(\Psi,\l)$.
Hence $(\Psi,\l)$ is bijective and so a homeomorphism.
\end{proof}
\end{subsection}
\begin{subsection}{More on the arc complex}
We have just seen that there is a correspondence
between points belonging to proper simplices of $|A(S,P)|$
and $P$-marked Riemann surfaces $\tilde{S}$ of genus $g$
together with an oriented diffeomorphism $f:S\arr{\sim}{\lra} \tilde{S}$
and a perimeter length $2l$, with $l\in\D_P$.

Now we recall how to describe points of $|\A8(S,P)|\subset |A(S,P)|$.
This has already been done in \cite{bowditch-epstein:triangulations}
and \cite{penner:simplicial} in the hyperbolic setting and it can
be adapted to the conformal case with minor modifications.\\

Keep the notation as in the previous subsection and
fix an improper simplex $\ua=\{\a_0,\dots,\a_k\}$
of the arc complex $A(S,P)$.
As before, $\ua$ determines a (possibly disconnected)
ribbon graph $G^*_{\ua}=(X(G^*_{\ua}),\s_0,\s_1)$
with a realization in $S$,
but in this case the realization is not proper and
so it does not determine a cellular structure in a natural way.
The fact is that some components of
$S\setminus |G^*_{\ua}|$ are not disks or pointed disks
but may have positive genus or may contain many
marked points, thus no explicit link with Poincar\'e duality
is any longer possible.

Nevertheless, we can find an open subsurface $N\subset S$
(up to isotopy) containing $|G^*_{\ua}|$ such that
$\pa N$ is a disjoint union of circles and
$N\setminus |G^*_{\ua}|$
is a collection of disks, pointed disks
and cylinders that touch the boundary $\pa N$.
More concretely, $N$ is obtained as the union of
a fattening of $|G^*_{\ua}|$ inside $S$
and of those components of $S\setminus |G^*_{\ua}|$ that are
disks or pointed disks; thus, $N$ can well be disconnected
(for instance, see Fig.~\ref{fig:arc}).

If $N$ is not the whole $S$, then $N$ is not compact.
Call $\hat{N}$ the unique compactification of $N$ to
a surface obtained by adding a finite set $Q$ of ``special points'',
and notice that $G^*_{\ua}$ has a canonical proper realization
in $\hat{N}$.

Now we can construct the ``dual'' ribbon graph
$G_{\ua}=(G^*_{\ua})^*$ by simply setting
$G_{\ua}:=(X(G^*_{\ua}),\s_{\infty}^{-1},\s_1)$
as before.
Notice that the holes of $G_{\ua}$ are naturally marked
by some elements of $P$ and that $Q$ marks some vertices
in such a way that $G_{\ua}$ admits a canonical realization
in $\hat{N}$ up to isotopy which is compatible with the markings
(Fig.~\ref{fig:arc2}).

As a consequence, for every point $a\in\ua$
we can build a Riemann surface $S(G_{\ua},a)$ as explained before,
which comes endowed with an isotopy class of oriented
diffeomorphisms $\hat{N}\arr{\sim}{\lra} S(G_{\ua},a)$
and a Jenkins-Strebel
differential, whose critical graph realizes $G_{\ua}$.

The complement $N^c:=S\setminus N$ is just a topological subsurface and
it cannot be given a complex structure in a natural way.
\begin{figure}[h]
\resizebox{10cm}{!}{\input{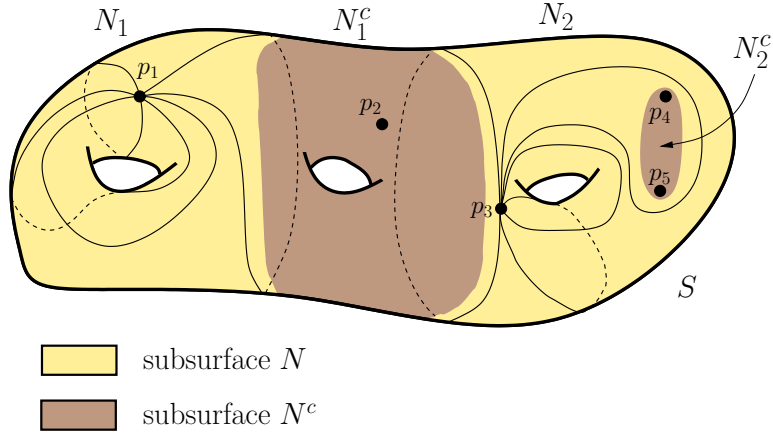}}
\caption{An arc system corresponding to an improper simplex}
\label{fig:arc}
\end{figure}
\begin{figure}[h]
\resizebox{10cm}{!}{\input{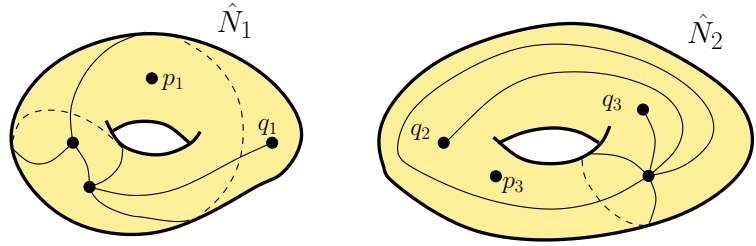}}
\caption{Surface $\hat{N}=\hat{N}_1\cup\hat{N}_2$
with a proper realization of $G_{\ua}$}
\label{fig:arc2}
\end{figure}
Hence, by means of Theorem \ref{th:js} applied to $N$ componentwise,
a point $a\in|A(S,P)|$ corresponds to an equivalence class
of couples $(N,l)$, where
\begin{itemize}
\item
$N\subset S$ is a nonempty open subsurface with complex structure,
such that $\pa N$ consists of disjoint circles, and $l$ is a point of $\D_P$
\item
every connected component $N_i$ of $N$
contains at least one point $p$ of $P$ such that $l_p>0$ and
the Euler characteristic of
$N_i\setminus P$ is negative
\item
no connected component of
the complement $N^c$
is a disk or a pointed disk
and $l_p=0$ for every $p \in P \cap N^c$
\end{itemize}
and two such couples $(N_1,l_1)$ and $(N_2,l_2)$ are to be considered
equivalent if $l_1=l_2$ and
there is a biholomorphism $N_1 \arr{\sim}{\lra} N_2$
such that the inclusion $N_1 \hra S$ and the composition
$N_1 \arr{\sim}{\lra} N_2 \hra S$ are isotopic
relatively to the markings.
\newline
\begin{figure}[h]
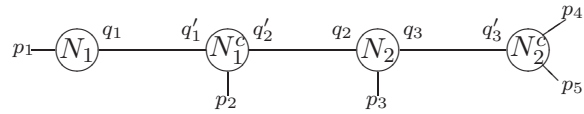

\[
\xy
\ar@{-}^{q_1\qquad q'_1}(3,0);(17,0)
\ar@{-}^{q'_2\qquad q_2}(23,0);(37,0)
\ar@{-}^{q_3\qquad q'_3}(43,0);(57,0)
\ar@{-}(-3,0);(-6,0)
\ar@{-}(20,-3);(20,-6)
\ar@{-}(40,-3);(40,-6)
\ar@{-}(62,2);(65,4)
\ar@{-}(64,-4);(66,-5)
,(-7,0)*{_{p_1}}
,(20,-7)*{_{p_2}}
,(40,-7)*{_{p_3}}
,(66,5)*{_{p_4}}
,(66,-5)*{_{p_5}}
,(0,0)*+{N_1}*\cir<9pt>{}
,(20,0)*+{N^c_1}*\cir<9pt>{}
,(40,0)*+{N_2}*\cir<9pt>{}
,(60,0)*+{N^c_2}*\cir<9pt>{}
\endxy
\]
\caption{Dual graph of $(N,l)$ as in Figures ~\ref{fig:arc} and \ref{fig:arc2}}
\label{fig:dual_graph}
\end{figure}
\newline
As a consequence, we can define a map
\[
\ol{\Phi}:|A(S,P)|/\G_{S,P} \lra \Mbartri_{g,P}
\]
in the following way.
Pick a point $a\in|A(S,P)|$ and consider the associated
$N\subset S$. Then consider the surface $\ol{S}$ obtained from $S$
collapsing every circle of $\pa N$ to a point (Fig.~\ref{fig:arc3}).
Moreover the positive components of $\ol{S}$
(which correspond to the connected components of $\hat{N}$)
are given a complex
structure thanks to the diffeomorphism $S(G_{\ua},a)\cong\hat{N}$.
Hence, we can set $\ol{\Phi}(a):=[\ol{S},\l(a)]$.
\begin{figure}[h]
\resizebox{11cm}{!}{\input{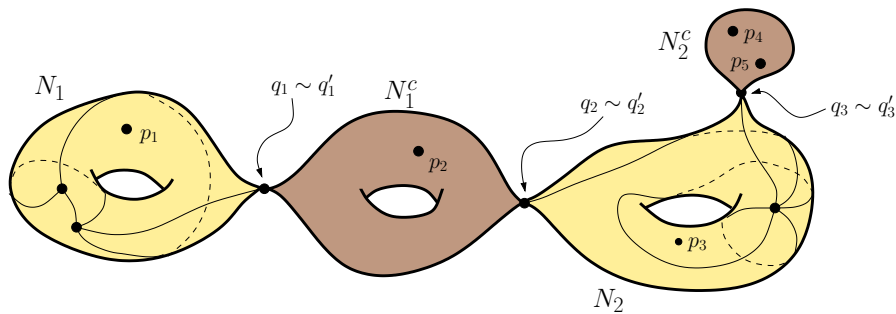}}
\caption{Surface $\ol{S}$ associated to $(N,l)$ of Fig.~\ref{fig:arc}}
\label{fig:arc3}
\end{figure}
It is easy to see that $\ol{\Phi}$ is well-defined and bijective.
One can also prove that $\ol{\Phi}$ is continuous
(see \cite{looijenga:cellular} and the Appendix),
and so it is a homeomorphism.

\begin{remark}
Unfortunately, the arc complex does not give a cellularization of the
Deligne-Mumford compactification of $\M_{g,P}$. On the contrary,
we have just seen that there is a proper surjective map
\[
\Mbar_{g,P}\times\D_P \lra |A(S,P)|/\G_{S,P}
\]
which is in fact a quotient.
In the Appendix, we sketch the construction
of Looijenga's modification of the arc complex $\wh{A}(S,P)$
(see \cite{looijenga:cellular}) which ``nearly'' does the job.
\end{remark}
\end{subsection}
\begin{subsection}{The ribbon graph complex}
\label{subsec:ribbon}
Here we introduce the second complex we are interested in,
for which we follow Kontsevich \cite{kontsevich:intersection}.
The point of view is reversed: the central object is
the ribbon graph and no longer the arc system
(a point of view which was already present in
\cite{penner:euler}).

Form the category $\mathcal{RG}_{g,P}$ of
$P$-marked ribbon graphs
of genus $g$ as follows. Its objects are
the ribbon graphs $G_{\ua}$ with $\ua$ in
$\Ao(S,P)$, and its morphisms are compositions of
isomorphisms of pointed ribbon graphs and contractions of one edge.
Denote by $\M$ (resp. $\Mbar$) the functor
$\mathcal{RG}_{g,P} \lra \mathrm{(Top. \,\, spaces)}$
that associates $(|\ua|\cap |\Ao|)\times\R_+$
(resp. $|\ua|\times\R_+$) to every $G_{\ua}$
and by $\Mcomb_{g,P}$ (resp. $\Mbarcomb_{g,P}$) its limit
functor. Remark that both functors are represented by
orbicellular complexes and that
$\Mcomb_{g,P} \subset \Mbarcomb_{g,P}$
embeds as a dense open subspace. Moreover we can define a
circumference function
$\bar{\l}:\Mbarcomb_{g,P} \rar \R_{\geq 0}^P\setminus\{0\}$
as in the case of the arc complex.
\begin{remark}
Notice that our definition of $\Mcomb_{g,P}$ and $\Mbarcomb_{g,P}$
slightly differs from that of Kontsevich. In fact, we allow some
holes to have perimeter zero (i.e. we admit marked vertices)
while Kontsevich does not. Briefly, Kontsevich's spaces are obtained
from ours by intersecting $\Mcomb_{g,P}$ and $\Mbarcomb_{g,P}$
with $\bar{\l}^{-1}(\R_+^P)$.
\end{remark}
Observe that the points of $\Mbarcomb_{g,P}$ correspond
to the following data:
\begin{enumerate}
\item[-]
a connected graph $\g$ (the ``dual graph of the pointed surface'')
with vertices $V$ labelled
by pairs $(g_v,P_v)$ such that $\sqcup_{v\in V} P_v=P$
and $\sum_{v\in V} g_v+\mathrm{dim}\, H^1(|\g|)=g$
(where $|\g|$ is a topological
realization of $\g$)
\item[-]
a subset $V_+\subset V$ of vertices of $\g$ (the ``positive vertices
of the dual graph'')
\item[-]
for every vertex $v\in V_+$
an ordinary $P_v\cup Q_v$-marked ribbon graph $(G_v,x_v)$ of genus $g_v$
with positive metric (possibly of total length different from $1$)
such that $Q_v$ marks only vertices of $G_v$, where $Q_v$
bijectively correspond to the set of
half-edges of $\g$ incident on $v$.
\end{enumerate}
We require moreover that no edge of $\g$ joins two nonpositive
vertices (i.e. that the $\g$ is reduced).

It is easy to see that the natural map
\[
F:\Mbarcomb_{g,P} \arr{\sim}{\lra}
\Mbartri_{g,P}\times\R_+
\]
is bijective and proper (see \cite{kontsevich:intersection}).
One can prove that $F$ is also continuous (see \cite{looijenga:cellular}
and the Appendix), so $F$ is a homeomorphism.

We summarize the preceding observations in the following commutative diagram
\[
\xymatrix@!C=2cm{
&&& \Mbar_{g,P}\times (\R_{\geq 0}^P\setminus\{0\}) \ar@{->>}[d]^{\xi} \\
& \M_{g,P}\times (\R_{\geq 0}^P\setminus\{0\})
\ar@{^(->}[urr] \ar@{^(->}[rr]
&& \Mbartri_{g,P}\times \R_+ \\
(|\Ao(S,P)|/\G_{S,P})\times\R_+ \ar[ru]_{\cong}^{\Phi} \ar@{^(->}[rr]
\ar[dr]_{\cong}
& \ar[u]
& (|A(S,P)|/\G_{S,P})\times\R_+ \ar[ru]^{\ol{\Phi}}_{\cong}
\ar[dr]_{\cong} \\
& \Mcomb_{g,P} \ar@{^(->}[rr] \ar@{-}[u]_{\cong}
&& \Mbarcomb_{g,P} \ar[uu]^{F}_{\cong}
}
\]
and we recall that $\xi$ is the map that collapses nonpositive
components so that its fibers are isomorphic to products
of smaller moduli spaces.
Notice that the homeomorphism $F^{-1}\ol{\Phi}$ is induced by
the $\G_{S,P}$-equivariant map
\[
|A(S,P)|\times\R_+ \lra \Mbarcomb_{g,P}
\]
that sends a point $a\in |A(S,P)|$ corresponding to $(N,l)$
to the point of $\Mbarcomb_{g,P}$ associated to
$(\g_{\ol{S},l},G_a)$, where $\g_{\ol{S},l}$ is the (reduced) dual
graph of $\ol{S}$ (which is obtained from $S$ by collapsing every
component of $\pa N$ to a point) and $G_a$ is the (possibly disconnected)
metrized ribbon graph, whose components correspond to
positive vertices of $\g_{\ol{S},l}$.

In what follows we will always identify $\Mbarcomb_{g,P}$,
$|A(S,P)|\times\R_+$ and $\Mbartri_{g,P}\times\R_+$,
so that we will regard $\xi$
as a map from $\Mbar_{g,P}\times(\R_{\geq 0}\setminus\{0\})$
to $\Mbarcomb_{g,P}$.
\begin{remark}
The orbispace $\Mbarcomb_{g,P}$
has an orbi-piecewise-linear structure;
so de Rham theorem holds and gives an isomorphism between
rational singular cohomology and rational PL de Rham
cohomology.
In what follows all cohomology groups
will be considered with rational coefficients,
even though tautological and combinatorial classes are defined
over $\Z$, so that all results still hold in integral cohomology
modulo torsion.
\end{remark}

Now we define
\[
\Scomb_{g,P}(l):=
\left\{ (N,l)\in\Mbarcomb_{g,P}\,\Big|\,
\begin{matrix}
\text{a component of $N^c$ is not}\\
\text{a sphere with $2$ or $3$ punctures}
\end{matrix}
\right\}
\]
and similarly
\[
\S_{g,P}(l):=
\left\{ (C,l)\in\Mbar_{g,P}\times\{l\} \,\Big|\,
\begin{matrix}
\text{$(C,l)$ has a nonpositive component}\\
\text{which is not a $3$-punctured sphere}
\end{matrix}
\right\}
\]
in such a way that that $\xi^{-1}(\Scomb_{g,P}(l))=\S_{g,P}(l)$.
One can check that the restriction of $\xi$ to
$\Mbar_{g,P}\times\{l\}\setminus\S_{g,P}(l)$ is a
homeomorphism onto $\Mbarcomb_{g,P}(l)\setminus\Scomb_{g,P}(l)$,
so that $H_*(\Mbar_{g,P}\times\{l\},\S_{g,P}(l))\cong
H_*(\Mbarcomb_{g,P}(l),\Scomb_{g,P}(l))$.

Moreover we set
\[
\S^Q_{g,P}:=
\left\{ C \in\Mbar_{g,P\cup Q} \,\Big|\,
\begin{matrix}
\text{$C$ has a component without $P$-markings}\\
\text{which is not a $3$-punctured sphere}
\end{matrix}
\right\}
\]
so that, for every $\tilde{l}=(l,0)\in\R_+^P\times\{0\}^Q$,
the obvious identification
$\Mbar_{g,P\cup Q}\cong\Mbar_{g,P\cup Q}\times\{\tilde{l}\}$
restricts to $\S^Q_{g,P}\cong\S_{g,P\cup Q}(\tilde{l})$.

These relative homeomorphisms show that one can lift
homology classes from $\Mbarcomb_{g,P}$ to $\Mbar_{g,P}$,
but one has to pay the price of some ambiguity
coming from the homology of $\S_{g,P}$.

\end{subsection}
\end{section}
%
%
\begin{section}{Tautological classes}\label{section:tautological}
Let $g,n$ be nonnegative integers such that $2g-2+n>0$
and let $P:=\{p_1,\dots,p_n\}$.
The moduli spaces $\Mbar_{g,P}$ of complex projective stable curves
have also the structure of smooth proper Deligne-Mumford
stacks over $\C$, so it is possible to define the Chow intersection
ring with rational coefficients $CH^*(\Mbar_{g,P})_{\Q}$ 
(as the $\Mbar_{g,P}$ are global quotients
of smooth projective varieties by finite groups
it is also possible to define integral Chow rings and
$CH^*(\Mbar_{g,P})_{\Q}=
CH^*(\Mbar_{g,P})\otimes\Q$).

There are three kinds of natural maps:
\begin{enumerate}
\item
the proper and flat map $\pi_q: \Mbar_{g,P\cup\{q\}} \rar \Mbar_{g,P}$
that forgets the point $q$ and stabilizes the curve
(i.e. contracts two-pointed spheres), which is isomorphic to
the universal family and so is endowed with
natural sections $\vartheta_{0,\{p_i,q\}}:
\Mbar_{g,P} \rar \Mbar_{g,P\cup\{q\}}$
for all $p_i\in P$
\item
the proper map $\vartheta_{irr}: \Mbar_{g-1,P\cup\{p',p''\}} \rar \Mbar_{g,P}$
(defined for $g>0$) that glues $p'$ and $p''$ together
\item
the proper map $\vartheta_{g',I}:
\Mbar_{g',I\cup\{p'\}} \times \Mbar_{g-g',I^c\cup\{p''\}}
\rar \Mbar_{g,P}$ (defined for every $0\leq g'\leq g$ and $I\subseteq P$
such that the spaces involved are nonempty)
that takes two curves and glues them together identifying $p'$ and $p''$.
\end{enumerate}
Call $\d_{irr}\subset \Mbar_{g,P}$ and
$\d_{g',I}\subset\Mbar_{g,P}$ the divisors supported
on the image of $\vartheta_{irr}$ and $\vartheta_{g',I}$.
When no confusion can occur, we will denote by
the same symbol the divisors and the associated classes
in the Chow ring or in homology.

If $\omega_{\pi_q}$ denote the relative dualizing sheaf
and $\mathcal{L}_{p_i}:=\vartheta^*_{0,\{p_i,q\}}(\omega_{\pi_q})$,
then define the Miller classes as
\[
\psi_{p_i}:=c_1(\mathcal{L}_{p_i})
\in CH^1(\Mbar_{g,P})_{\Q}
\]
and the modified (by Arbarello-Cornalba)
Mumford-Morita classes as
\[
\k_j:=(\pi_q)_*(c_1(\omega_{\pi_q}(\sum_i \d_{0,\{p_i,q\}} ))^{j+1})
\in CH^j(\Mbar_{g,P})_{\Q}.
\]
When there is no risk of ambiguity,
we will denote in the same way the line bundles
$\mathcal{L}$ and the classes $\psi$ and $\k$
belonging to different $\Mbar_{g,P}$'s.
Because of the natural definition of $\k$ and $\psi$ classes,
the subring $R^*(\M_{g,P})$ of $CH^*(\M_{g,P})_{\Q}$
they generate is called {\it tautological ring}
of $\M_{g,P}$. Its image $RH^*(\M_{g,P})$ through the cycle class map
is called cohomology tautological ring.

The {\it system of tautological rings}
$(R^*(\Mbar_{g,P}) \subset
CH^*(\Mbar_{g,P})_{\Q} \ )_{2g-2+n>0}$ is
the minimal system of subrings that contain the
classes $\k$ and $\psi$ which is
closed under the push-forward maps $\pi_*$, $(\vartheta_{irr})_*$
and $(\vartheta_{g',I})_*$. The definition is the same
for the rational cohomology.

The psi classes and the kappa classes are related in a very nice way.
In fact
\begin{align*}
(\pi_q)_*(\psi_{p_1}^{r_1}\cdots\psi_{p_n}^{r_n})=
\sum_{\{i| r_i>0\}} \psi_{p_1}^{r_1}\cdots\psi_{p_i}^{r_i-1}
\cdots\psi_{p_n}^{r_n} \\
(\pi_q)_*(\psi_{p_1}^{r_1}\cdots\psi_{p_n}^{r_n}\psi_q^{b+1})=
\psi_{p_1}^{r_1}\cdots\psi_{p_n}^{r_n}\k_b
\end{align*}
where the first one is the so-called {\it string equation} and
the second one for $b=0$ is the {\it dilaton equation}.

In \cite{arbarello-cornalba:combinatorial}
Arbarello and Cornalba noticed that the previous formulas,
together with the relation
\[
\pi_q^*(\k_b)=\k_b+\psi_q^b
\]
and a repeated use of the push-pull formula,
make it possible to compute the push-forward of
any polynomial in the $\psi$'s through the forgetful maps.

The following explicit formula for maps that forget more than one point
was found by Faber.
Let $Q:=\{q_1,\dots,q_m\}$
and let $\pi_Q:\Mbar_{g,P\cup Q}\rar\Mbar_{g,P}$ be the forgetful map.
Then
\[
(\pi_Q)_*(\psi_{p_1}^{r_1}\cdots\psi_{p_n}^{r_n}
\psi_{q_1}^{b_1+1}\cdots\psi_{q_m}^{b_m+1})=
\psi_{q_1}^{r_1}\cdots\psi_{q_n}^{r_n}
K_{b_1\cdots b_m}
\]
where $K_{b_1\cdots b_m}=\sum_{\s\in\mathfrak{S}_m} \k_{b(\s)}$ and
$\k_{b(\s)}$ is defined in the following way.
If $\g=(c_1,\dots,c_l)$ is a cycle,
then set $b(\g):=\sum_{j=1}^l b_{c_j}$.
If $\s=\g_1\cdots\g_{\nu}$ is the decomposition in
disjoint cycles (including 1-cycles), then we define
$k_{b(\s)}:=\prod_{i=1}^{\nu} \k_{b(\g_i)}$.

Another property we will use is the following.
Consider the map
\[
\vartheta_{0,\{p_i,q\}}:\Mbar_{g,P}\lra\Mbar_{g,P\cup\{q\}}.
\]
One can easily check that the line bundle
$\vartheta_{0,\{p_i,q\}}^*(\mathcal{L}_{q})$ is trivial,
so that $\vartheta_{0,\{p_i,q\}}^*(\psi_q)=0$
in $H^2(\Mbar_{g,P})$.
If we call $D_q$ the (disjoint) union $\cup_i \d_{0,\{q,p_i\}}$,
then the exact sequence of the couple
\[
0=H^1(D_q) \lra H^2(\Mbar_{g,P\cup\{q\}},D_q)
\lra H^2(\Mbar_{g,P\cup\{q\}})
\lra H^2(D_q) \lra
\]
shows that $\psi_q$ lifts to a class in $H^2(\Mbar_{g,P\cup\{q\}},D_q)$.

Similarly, consider the following situation.
Let $Q=\{q_1,\dots,q_m\}$ and for every $Q'\subset Q$
define $D_{Q',P}\subset\Mbar_{g,P\cup Q}$
to be the union of all divisors of the type
$\d_{0,\{q_{j_1},\dots,q_{j_h},p_i\}}$
with $p_i\in P$ and $\{q_{j_1},\dots,q_{j_h}\}\subset Q'$.
Then, for every $b_1,\dots,b_m\geq 1$, the class
$\psi_{q_1}^{b_1}\cdots\psi_{q_m}^{b_m}$ admits a lift
to $H^*(\Mbar_{g,P\cup Q},D_{Q',P})$.
In fact, as $\pi_q^*(\psi_{p_i})=\psi_{p_i}-\d_{0,\{p_i,q\}}$,
the class we are interested in coincides with
\[
\pi_{q_m,\dots,q_2}^*(\psi_{q_1}^{b_1})   \cdot
\pi_{q_m,\dots,q_3}^*(\psi_{q_2}^{b_2})   \cdots
\pi_{q_m}^*(\psi_{q_{m-1}}^{b_{m-1}})     \cdot
\psi_{q_m}^{b_m}
\]
which lies in the image of the homomorphism
\[
\bigotimes_{j=1}^m
H^*(\Mbar_{g,P\cup \{q_1,\dots,q_j\}},D_{\{q_j\},P})
\lra
H^*(\Mbar_{g,P\cup Q},D_{Q',P})
\]
which pulls classes back via appropriate forgetful
morphisms and then multiplies them together.\\

We refer to \cite{kmz} for some remarks
on Faber's formula, its inversion and for other useful formulas
on Weil-Petersson volumes.

Moreover, we refer \cite{arbarello-cornalba:combinatorial} for
a proof of the above relations between $\kappa$ classes
and $\psi$ classes, and to \cite{faber:conjectures} for
a conjectural description of the tautological rings.
\end{section}
%
%
\begin{section}{Combinatorial classes}\label{section:combinatorial}
Fix a compact oriented surface $S$ of genus $g$ and a nonempty subset
$P=\{p_1,\dots,p_n\}$ of $n$ distinct points of $S$ such that $2g-2+n>0$.

Let $m_*=(m_{-1},m_0,m_1,\dots)$ be a sequence of nonnegative integers
such that
\[
\sum_{i\geq -1} (2i+1)m_i=4g-4+2n
\]
and define $(m_*)!:=\prod_{i\geq -1}m_i!$
and $r:=\sum_{i\geq -1} i\, m_i$.
\begin{warning}
Even if one could still define combinatorial classes with $m_{-1}>0$,
much more care is needed to handle them. Moreover, in this case the argument
of Lemma~\ref{lemma:enum} would not be sufficient to complete
the proof of Theorem~\ref{th:second}. Therefore
we assume $m_{-1}=0$ in what follows.
\end{warning}
Reasoning as in Section~\ref{subsec:ribbon}, it is possible to construct
an orbispace $\Mbarcomb_{m_*,P}$ whose cells of maximal dimension
are indexed by isomorphism classes of ordinary ribbon graphs that have
exactly $m_i$ vertices of valency $2i+3$. Analogously it is
possible to define an arc complex $A(S,P)_{m_*}$
as the smallest subcomplex of $A(S,P)$
that contains all
simplices $\ua$ such that $G_{\ua}$ is an ordinary ribbon
graph with exactly $m_i$ vertices
of valency $2i+3$. Notice that both these complexes are acted
on by $\G_{S,P}$ and so is
$\Ao(S,P)_{m_*}:=A(S,P)_{m_*}\cap\Ao(S,P)$.
Thus, the following diagram is obviously commutative.
\[
\xymatrix{
(|\Ao(S,P)_{m_*}|/\G_{S,P})\times\R_+ \ar[r]^{\qquad\qquad\cong} \ar@{^(->}[d]
       &  \Mcomb_{m_*,P} \ar@{^(->}[d]  \\
(|A(S,P)_{m_*}|/\G_{S,P})\times\R_+ \ar[r]^{\qquad\qquad\cong} &
               \Mbarcomb_{m_*,P} 
}
\]
In what follows, we will naturally identify arc complexes with
combinatorial moduli spaces.
\begin{remark}
In the case $m_{-1}>0$ it is still possible to define
$A(S,P)_{m_*}$ (and $\Ao(S,P)_{m_*}$)
as a subcomplex of an extended arc complex
$\widetilde{A}(S,P)$, obtained adding
to $\mathcal{A}$ contractible loops (i.e. unmarked tails
in the corresponding ribbon graph picture).
However $A(S,P)_{m_*}$ is no longer a subcomplex of $A(S,P)$, so
we can only define
$\Mcomb_{m_*,P}\rar \M_{g,P}\times (\R_{\geq 0}^P\setminus\{0\})$
as a classifying map by constructing a family of Riemann surfaces
over $\Mcomb_{m_*,P}$.
\end{remark}
Now consider the following commutative diagram.
\[
\xymatrix{
 \Mcomb_{m_*,P} \ar@{^(->}[d] \ar@{^(->}[r] & \Mcomb_{g,P} \ar@{^(->}[d]
       & \M_{g,P}\times(\R^P_{\geq 0}\setminus\{0\}) \ar[l]_{\cong\qquad} \ar@{^(->}[d] \\
\Mbarcomb_{m_*,P} \ar@{^(->}[r] & \Mbarcomb_{g,P}
       & \Mbar_{g,P}\times(\R^P_{\geq 0}\setminus\{0\}) \ar@{->>}[l]_{\xi\qquad}
}
\]
For every $l \in \R_{\geq 0}^P\setminus\{0\}$
call $\Mbarcomb_{g,P}(l)$ the {\it slice}
$\bar{\l}^{-1}(l) \subset \Mbarcomb_{g,P}$.
In the same way we can define $\Mbarcomb_{m_*,P}(l)$
and the restriction
\[
\xymatrix{
\xi_l: \Mbar_{g,P} \ar@{->>}[r] & \Mbarcomb_{g,P}(l)
}
\]
of $\xi$.
Notice that the dimensions of the slices are
the expected ones, because in every cell they are described
by a system of $n$ independent linear equations.
\begin{lemma}[\cite{kontsevich:intersection}]
Fix $p$ in $P$ and $l\in\R_{\geq 0}^P$ such that $l_p>0$.
Then, on every simplex $\ua \in \Mbarcomb_{g,P}(l)$, define
\[
\ol{\omega}_p|_{\ua}:=\sum_{1\leq s<t\leq k}
d \left( \frac{e_s(a)}{2l_p} \right)
\wedge
d \left( \frac{e_t(a)}{2l_p} \right)
\]
where $x(p)$ is a hole with cyclically ordered sides
$(e_1,\dots,e_k)$. These $2$-forms glue to
give a piecewise-linear $2$-form
$\ol{\omega}_p$ on $\Mbarcomb_{g,P}(l)$
and the class
$\xi_l^*[\ol{\omega}_p]$ pulls back to
$\psi_p$ in $H^2(\Mbar_{g,P})$.
\end{lemma}
\begin{lemma}[\cite{kontsevich:intersection}]
For every $l\in\R_{\geq 0}^P\setminus\{0\}$
the restriction of $\ol{\Omega}:=\sum_{p\in P} l_p^2 \ol{\omega}_p$
to $\Mbarcomb_{g,P}(l)$ is a nondegenerate symplectic form, so
$\ol{\Omega}^r$ defines an orientation on $\Mbarcomb_{m_*,P}(l)$.
Hence, $\ol{\Omega}^r\wedge\bar{\l}^*d\vol_{\R_+^P}$
is an orientation on $\Mbarcomb_{m_*,P}(\R_+^P)$.
\end{lemma}
\begin{lemma}[\cite{kontsevich:intersection}]
With the given orientation $\Mbarcomb_{m_*,P}(l)$ is a cycle
for all $l\in\R_{\geq 0}^P\setminus\{0\}$
and $\Mbarcomb_{m_*,P}(\R_+^P)$ is
a cycle with noncompact support.
\end{lemma}

Notice that the space $\Mbarcomb_{m_*,P}$
reduces to $\Mcomb_{m_*,P}$ when restricted to the locus of
ordinary ribbon graphs, and it coincides with the closure
of $\Mcomb_{m_*,P}$
in $\Mbarcomb_{g,P}$.

Define the {\it combinatorial classes}
$W_{m_*,P}(l):=[\Mcomb_{m_*,P}(l)]$ and
$\Wbar_{m_*,P}(l):=[\Mbarcomb_{m_*,P}(l)]$, and observe that
K\"unneth formula and $H^{BM}_*(\R_+^n)\cong\Q[-n]$ give the
following isomorphism
\[
\xymatrix@R=0.1cm{
H^{BM}_{6g-6+3n-2r}(\Mcomb_{g,P}(\R_+^P)) \ar[r]^{\sim} &
H^{BM}_{6g-6+2n-2r}(\Mcomb_{g,P}(l)) \\
W_{m_*,P}(\R_+^P) \ar@{|->}[r] & W_{m_*,P}(l)
}
\]
and analogously
\[
\xymatrix@R=0.1cm{
H^{BM}_{6g-6+3n-2r}(\Mbarcomb_{g,P}(\R_+^P)) \ar[r]^{\sim} &
H_{6g-6+2n-2r}(\Mbarcomb_{g,P}(l)) \\
\Wbar_{m_*,P}(\R_+^P) \ar@{|->}[r] & \Wbar_{m_*,P}(l)
}
\]
for every $l \in \R_+^P$, naturally with respect the inclusion
$\Mcomb_{g,P}\hra\Mbarcomb_{g,P}$.
Therefore,
we will write $W_{m_*,P}$ and $\Wbar_{m_*,P}$
instead of $W_{m_*,P}(l)$ and $\Wbar_{m_*,P}(l)$
for the homology classes they define on
$\M_{g,P}(l)\cong\M_{g,P}$
and $\Mbarcomb_{g,P}(l)$ respectively,
for any $l\in\R_+^P$.
Moreover, we will also identify $W_{m_*,P}$ with its
Poincar\'e dual in $H^{2r}(\M_{g,P})$. \\

It is possible to define a slight generalization of
the previous classes, prescribing that some markings
hit vertices with assigned valency.

Given a finite set $Q=\{q_1,\dots,q_h\}$ and a map
$\rho:Q\rar \Z_{\geq -1}$, we define
$m^{\rho}_*=(m^{\rho}_{-1},m^{\rho}_0,\dots)$
as $m^{\rho}_i:=|\rho^{-1}(i)|$.
Consider now an $m_*$ and a $\rho$ such that $m^{\rho}_{-1}=m_{-1}$,
$m^{\rho}_*\leq m_*$ and 
$\sum_{i\geq -1}(2i+1)m_i=4g-4+2|P|$, and
call $\Mbarcomb_{m_*,\rho,P}$ the subcomplex of
$\Mbarcomb_{m_*,P\cup Q}$ whose
simplices of maximal dimension are ordinary ribbon
graphs in which $q_j$ marks a vertex of valency $2\rho(q_j)+3$
for every $j=1,\dots,h$.
Then, denote by
$\Wbar_{m_*,\rho,P}$
its homology class inside some slice
of $\Mbarcomb_{g,P\cup Q}(\tilde{l})$
with $\tilde{l}\in \R_+^P\times\{0\}^Q$
(as before the orientation is determined
by $\sum_{p\in P}l_{p_i}^2\ol{\omega}_p$).
Define analogously
$\Mcomb_{m_*,\rho,P}=\Mbarcomb_{m_*,\rho,P}\cap\Mcomb_{g,P\cup Q}$
and let $W_{m_*,\rho,P}$ be its Borel-Moore homology
class inside the same slice of
$\Mcomb_{g,P\cup Q}$. \\

Notice that the combinatorial classes $\Wbar_{m_*,P}$
are not defined on the Deligne-Mumford compactification but
only on the combinatorial compactification
(analogously for $\Wbar_{m_*,\rho,P}$). This last space
is often singular, so we cannot use Poincar\'e duality in order
to obtain a cohomology class to lift to the Deligne-Mumford
compactification.

However the singular locus
of $\Mbarcomb_{g,P}(l)$
is contained inside $\Scomb_{g,P}(l)$.
Thus, the following diagram
\[
\xymatrix@R=0.3cm{
H_{6g-6+2n-2r}(\Mbarcomb_{g,P}(l)) \ar[r] &
H_{6g-6+2n-2r}(\Mbarcomb_{g,P}(l),\Scomb_{g,P}(l))
\ar@{=}[d]^{\wr} \\
& H_{6g-6+2n-2r}(\Mbar_{g,P}(l),\S_{g,P}(l))
}
\]
shows that we can lift the combinatorial classes on the
Deligne-Mumford compactification up to some ambiguities.

\begin{remark}
One can try solve these ambiguities using a different combinatorial
compactification $\Mhat_{g,P}$, which maps to $\Mbar_{g,P}$.
Looijenga's modification of the arc complex (see \cite{looijenga:cellular})
offers a natural candidate, but it seems much more difficult to
prove that the combinatorial chains are in fact cycles
(see the Appendix).
\end{remark}

As it is evident from the definition, the cycles $\Wbar_{m_*,\rho,P}$ and
$\Wbar_{m_*,P}$ are closely related. It seems that one can obtain
a multiple of $\Wbar_{m_*,P}$ from $\Wbar_{m_*,\rho,P}$ by
simply ``forgetting'' the $Q$-marking, so that one
would like to deduce a kind of push-forward formula
for combinatorial classes which looks like
\begin{center}
``$(\pi_Q)_*(\Wbar_{m_*,\rho,P})=
N(m_*,\rho) \Wbar_{m_*,P}$''
\end{center}
where the integer coefficient $N(m_*,\rho)$ counts in how many ways
one can obtain a top-dimensional cell of $\Wbar_{m_*,\rho,P}$
by $Q$-marking vertices of a fixed top-dimensional cell
of $\Wbar_{m_*,P}$.

\begin{example}
In general, $\pi_Q(\S^Q_{g,P})$
is not contained inside $\S_{g,P}$
or even in the boundary $\pa\M_{g,P}$.
For instance, in the case $Q=\{q\}$, the image $\pi_q(\S^q_{g,P})$
does not lie inside $\S_{g,P}$,
but $\pi_q(\S^q_{g,P})\setminus \S_{g,P}$
contains the locus of surfaces with an unmarked
three-noded sphere. Hence it has complex codimension three.
Nevertheless, in this particular case,
$\pi_q(\S^q_{g,P})$ is contained inside $\pa\M_{g,P}$.
\end{example}
\end{section}
%
%
\begin{section}{Shrinking map and combinatorial forgetful map}
\label{section:pi}
Before introducing some technical tools, we want to describe
the basic ideas involved in the proof of Theorem \ref{th:first}.

The first observation is that $\Mcomb_{g,P\cup\{q\}}(\tilde{l})$ 
is homeomorphic
to $\M_{g,P\cup\{q\}}$ for every
$\tilde{l}\in \R_{\geq 0}^{P\cup\{q\}}\setminus\{0\}$.
The second remark is that the differential form
$\omega_q:=\ol{\omega}_q|_{\Mcomb_{g,P\cup\{q\}} }$ lives
on the slices $\Mcomb_{g,P\cup\{q\}}(\tilde{l})$ such that $\tilde{l}_q>0$,
while the generalized combinatorial class
(which we will briefly denote by $W^q_{2r+3}$),
defined prescribing that
$q$ marks a vertex of valency (at least) $2r+3$, lives on
the slices that have $\tilde{l}_q=0$.

Consequently, a deformation retraction $\H^q_0$ of $\Mcomb_{g,P\cup\{q\}}$
onto the slice defined by $\tilde{l}_q=0$ would help us to compare
$\omega_q^{r+1}$ and the combinatorial class $W^q_{2r+3}$
as functionals on the cohomology of $\Mcomb_{g,P\cup\{q\}}(\tilde{l})$.

However, the deformation retraction $\H^q_0$ we will construct
does not preserve the locus of smooth curves, but it
retracts $\Mbarcomb_{g,P\cup\{q\}}$ onto
the slice $\Mbarcomb_{g,P\cup\{q\}}(\{\tilde{l}_q=0\})$.
In fact, $\H^q_0$ is defined sending
all the edges bordering the \mbox{$q$-th} hole to zero
(and it is defined only when $\tilde{l}_q$ is ``small'', because
we must avoid the situation in which $\H^q_0$ would squeeze
another hole beside $q$).
Thus, it shrinks a circular \mbox{$q$-th} hole (i.e. such
that $\bar{T}_q$ is a disk) to a $q$-marked vertex;
while it produces a surface with a $q$-marked
nonpositive component, if the topology
of the \mbox{$q$-th} hole is more complicated.
Anyway, if we subdivide the complex $\Mbarcomb_{g,P\cup\{q\}}$
into subcomplexes $\ol{Y}^{\bullet}_*$
according to the topology of the \mbox{$q$-th} hole,
then the restriction of $\H^q_0$ to each subcomplex becomes
a simplicial fibration.

Finally, we consider a differential form $\eta$ on
$\Mbarcomb_{g,P\cup\{q\}}(\{\tilde{l}_q=0\})$ and we compare
the integral of $\eta$ on $\Wbar^q_{2r+3}(\tilde{l})$
(the closure of $W^q_{2r+3}(\tilde{l})$) for an $\tilde{l}$ such
that $\tilde{l}_q=0$ with the
integral of $\ol{\omega}_q^{r+1}\wedge(\H^q_0)^*\eta$
on $\Mbarcomb_{g,P\cup\{q\}}(\tilde{l})$ for an $\tilde{l}$
such that $\tilde{l}_q>0$.
Here we notice that the restriction of
the form $\ol{\omega}_q^{r+1}\wedge(\H^q_0)^*\eta$
to $\Mbarcomb_{g,P\cup\{q\}}(\tilde{l})$
has support on the top-dimensional simplices
whose \mbox{$q$-th} hole has exactly $2r+3$
distinct edges;
then the integral of $\ol{\omega}_q^{r+1}\wedge(\H^q_0)^*\eta$ is
evaluated by calculating
the integral of $\ol{\omega}_q^{r+1}$
on the fibers of the restriction of $\H^q_0$
to $\ol{Y}^{\bullet}_*$, for each subcomplex $\ol{Y}^{\bullet}_*$.
In the case of a
circular \mbox{$q$-th} hole with $2r+3$ edges,
from this integration
we obtain the coefficient $2^{r+1}(2r+1)!!$.

The analogous result for the ordinary combinatorial
class $\Wbar_{2r+3}$
(supported on the smallest subcomplex containing all
ribbon graphs with a $(2r+3)$-valent vertex)
and $\k_r$ is derived from the
previous one by noticing that the forgetful morphism
$\pi_q$ has a combinatorial analogue $\pi_q^{comb}$ on
the combinatorial spaces (another little technical problem
is due to the fact that $\pi_q^{comb}$ is not defined
on the whole $\Mbarcomb_{g,P\cup\{q\}}$, but this is also
superable).
As a consequence, $(\pi_q)_*(\psi_q^{r+1})=\k_r$
and $(\pi_q^{comb})_*$ sends $\Wbar_{2r+3}^q$ to $\Wbar_{2r+3}$.
Hence we obtain our result for the kappa classes too.

As the shrinking map $\H^q_0$ and the combinatorial forgetful map
$\pi_q^{comb}$ play a main role, this section will deal with their
definition and their properties.\\

%
%
Fix $g \geq 0$ and a nonempty set of markings $P=\{p_1,\dots,p_n\}$
such that $2g-2+n>0$ and define
\[
C_{P,q}:=\{\tilde{l}\in \R_{\geq 0}^{P\cup\{q\}}|\, \tilde{l}_q<\tilde{l}_p
\quad\text{for all $p\in P$}\}.
\]
Denote by $\pi_q:\Mbar_{g,P\cup\{q\}}\times C_{P,q}
\rar \Mbar_{g,P}\times\R_+^P$ the map that forgets $q$
and the $q$-th coordinate. We can define $\pi_q^{comb}$
forcing the commutativity of the following diagram
\begin{equation*}
\xymatrix{
(\Mbar_{g,P\cup\{q\}}\setminus D_q)\times C_{P,q}
\ar[r]^{\tilde{\xi}\qquad} \ar[d]_{\pi_q} &
\Mbarcomb_{g,P\cup\{q\}}(C_{P,q})\setminus D_q^{comb}(C_{P,q})
\ar[d]^{\pi_q^{comb}} \\
\Mbar_{g,P}\times\R_+^P \ar[r]^{\xi} & \Mbarcomb_{g,P}(\R_+^P)
}
\end{equation*}
where $D_q=\cup_{p\in P}\d_{0,\{q,p\}}$ and
$D_q^{comb}=\tilde{\xi}(D_q)$.
We remark that $\xi\pi_q$ does not factorize through
$\tilde{\xi}:\Mbar_{g,P\cup\{q\}}\times C_{P,q}\rar
\Mbarcomb_{g,P\cup\{q\}}(C_{P,q})$.
\begin{figure}[h]
\resizebox{7cm}{!}{\input{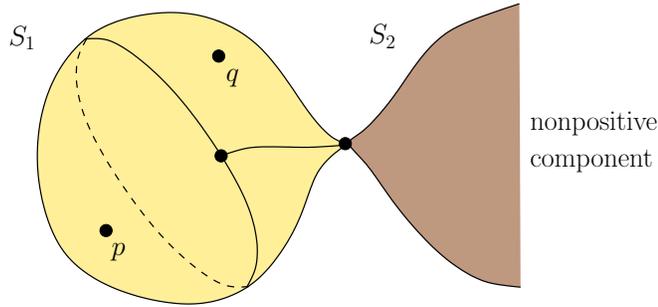}}
\caption{$\pi_q^{comb}$ is not defined in this case}
\label{fig:nonpositive}
\end{figure}
In fact pick a point $(S,\tilde{l})$ in $\Mbar_{g,P\cup\{q\}}\times C_{P,q}$
such that $q$ and $p$ lie on a two-pointed component
$S_1$ of $S$ of genus zero which has only one singular point and
suppose that the adjacent component $S_2$ is nonpositive
(see Fig.~\ref{fig:nonpositive}).
Then $\tilde{\xi}(S,\tilde{l})$
does not ``remember'' the analytic type of $S_2$
but $\xi\pi_q(S,\tilde{l})$
does (if $\tilde{l}_p>0$) because the $p$-marking now hits $S_2$
after forgetting $q$ and stabilizing. However this is
the only bad case, so it sufficient to excise $D_q$ and $D_q^{comb}$.
\begin{remark}
The behaviour of the map $\pi_q^{comb}$ is really mysterious as
we do not know in general how Jenkins-Strebel's differential changes when we
delete the marked point $q$ and consequently how the critical graph modifies.
However, we know that if $q$ marks a vertex
of valency at least two,
then the new critical
graph is obtained simply forgetting the marking.
On the contrary, we get in troubles if we try to
forget the $q$-marking and $q$ marks a univalent vertex.
This explains why it is more difficult to deal with
univalent vertices.
\end{remark}
\begin{notation}
Call $\ol{Y}_h\subset\Mbarcomb_{g,P\cup\{q\}}$ the closure
of the locus of graphs where
the $q$-marked hole has positive perimeter and
consists exactly of $h$ distinct (unoriented) edges.
Moreover, set $\ol{Y}_{\geq h}:=\cup_{i\geq h}\ol{Y}_i$
and $\ol{Y}_{\leq h}:=\cup_{i\leq h}\ol{Y}_i$.
\end{notation}
Clearly, the topological boundary $\pa\ol{Y}_{\geq h}$ is
contained inside $\ol{Y}_{\leq h-1}$.
Moreover, $\ol{Y}_{\geq 2}(C_{P,q})$ is contained
inside $\Mbarcomb_{g,P\cup\{q\}}(C_{P,q})\setminus D_q^{comb}(C_{P,q})$.
In fact, taking perimeters in $C_{P,q}$ is essential to our purposes;
a consequence of this choice is that
$\ol{Y}_1(C_{P,q})$ is a closed neighbourhood of $D_q^{comb}(C_{P,q})$
inside $\Mbarcomb_{g,P\cup\{q\}}(C_{P,q})$, which is false
for $\ol{Y}_1(\tilde{l})$ with $\tilde{l}\notin C_{P,q}$.

Now, since we will work with the differential form
$\ol{\omega}_q$, which
is defined only where $\tilde{l}_q>0$, then we
let the perimeters vary in the subset
$C_{P,q}^+:=C_{P,q}\cap\{\tilde{l}_q>0\}$ only.

Remark that
$\ol{\omega}_q|_{\ol{Y}_1(\tilde{l})}=0$ with $\tilde{l}\in C_{P,q}^+$,
because the hole $q$ does not contain enough edges.
Thus, if $\eta$ is a piecewise-linear differential form
on $\ol{Y}_{\geq 2}(C^+_{P,q})\subset \Mbarcomb_{g,P\cup\{q\}}$,
then $\eta\wedge\ol{\omega}_q^{r+1}$ regularly extends by
zero to the whole $\Mbarcomb_{g,P\cup\{q\}}(C^+_{P,q})$
for every $r\geq 0$.

The previous observations and the commutativity of the
following diagram
\[
\xymatrix{
(\Mbar_{g,P\cup\{q\}}\setminus D_q)\times C^+_{P,q}
\ar[r]^{\tilde{\xi}\qquad} \ar[d]_{\pi_q} &
\Mbarcomb_{g,P\cup\{q\}}(C^+_{P,q})\setminus D_q^{comb}(C^+_{P,q})
\ar[d]^{\pi_q^{comb}} \\
\Mbar_{g,P}\times\R_+^P \ar[r]^{\xi} & \Mbarcomb_{g,P}(\R_+^P)
}
\]
tell us how to use the combinatorial forgetful map.
Namely, for every $\tilde{l}\in C_{P,q}^+$
we can associate to $[\ol{\omega}_q]^{r+1}\in
H^{2r+2}(\Mbarcomb_{g,P\cup\{q\}}(\tilde{l}), D_q^{comb}(\tilde{l}))$
a homology class $[\ol{\omega}_q^{r+1}]^*$ in
$H_{6g-6+2n-2r}(\Mbarcomb_{g,P\cup\{q\}}(\tilde{l})
\setminus D_q^{comb}(\tilde{l}))$,
which is defined as the functional on
$H^{6g-6+2n-2r}(\Mbarcomb_{g,P\cup\{q\}}(\tilde{l})
\setminus D_q^{comb}(\tilde{l}))$
given by
\[
\eta \lmt
\int_{\Mbarcomb_{g,P\cup\{q\}}(\tilde{l})}\eta\wedge\ol{\omega}^{r+1}.
\]
Thus, $\tilde{\xi}_*(\psi_q^{r+1})^*=[\ol{\omega}_q^{r+1}]^*$ and
$\pi_{q,*}(\psi_q^{r+1})^*=\k_r^*$. Consequently,
$\pi^{comb}_{q,*}[\ol{\omega}_q^{r+1}]^*=\xi_*\k_r^*$.
\begin{proposition}
There is a deformation retraction
\[
\tilde{\H}^q:\Mbarcomb_{g,P\cup\{q\}}(C_{P,q})\times[0,1]\lra
\Mbarcomb_{g,P\cup\{q\}}(C_{P,q})
\]
such that $\tilde{\H}^q_1$ is
the identity and $\tilde{\H}^q_0$ is the piecewise-linear retraction
onto $\Mbarcomb_{g,P\cup\{q\}}(\R_+^P\times\{0\})$.
Moreover, $\tilde{\H}^q_t(\ol{Y}_h)\subset\ol{Y}_h$ and
$\tilde{\H}^q_t(D^{comb}_{q})\subset D^{comb}_q$
for all $t\in[0,1]$.
\end{proposition}
\begin{proof}
Consider a cell $\bar{\l}^{-1}(C_{P,q})\cap (|\ua|\times\R_+)$ inside
$\Mbarcomb_{g,P\cup\{q\}}(C_{P,q})$.
Denote by $e_1,\dots,e_h$ the coordinates of $|\ua|\times\R_+$
corresponding to the unoriented edges of $G_{\ua}$
that border the hole $q$ and by $f_1,\dots,f_k$ the remaining ones.
Then it is sufficient to define $\tilde{\H}^q_t$ as the map that sends
$e_i \lmt t\cdot e_i$ and $f_j \lmt f_j$ and to observe that
all these deformation retractions glue to give a global $\tilde{\H}^q$.
By definition, it is immediate that $\tilde{\H}^q_t(\ol{Y}_h)\subset \ol{Y}_h$
and one can also easily check that the locus $D^{comb}_q$ is
preserved.
\end{proof}
\begin{figure}[h]
\resizebox{12cm}{!}{\input{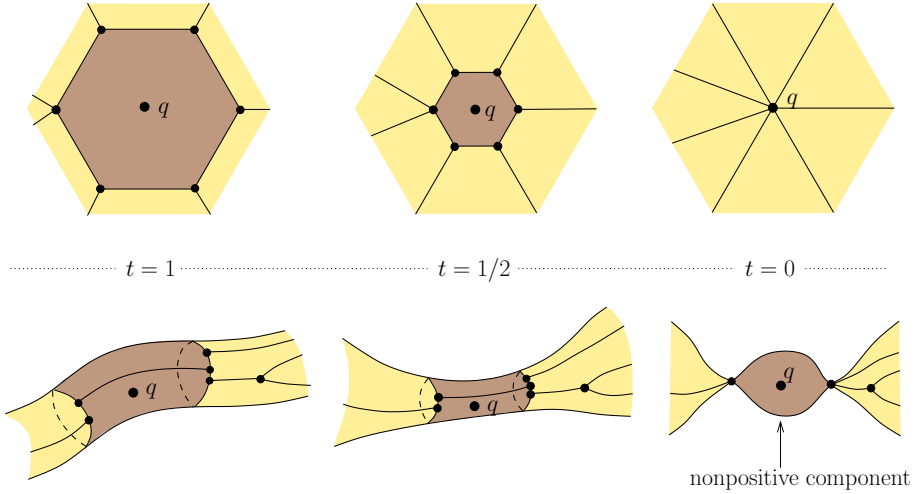}}
\caption{The deformation retraction $\tilde{\H}^q$}
\label{fig:retraction}
\end{figure}
\begin{definition}
The map $\H^q_0:\Mbarcomb_{g,P\cup\{q\}}(C_{P,q})\twoheadrightarrow
\Mbarcomb_{g,P\cup\{q\}}(\R_+^P\times\{0\})$
obtained from $\tilde{\H}^q_0$
by restriction of the codomain
is called {\it shrinking} of the $q$-th hole.
\end{definition}
A consequence of the previous proposition is that
we can compare $[\ol{\omega}_q^{r+1}]^*$ with
the combinatorial classes that live in
$\Mbarcomb_{g,P\cup\{q\}}(\R_+^P\times\{0\})$
using the following diagram.
\[
\xymatrix@C=3cm{
\Mbarcomb_{g,P\cup\{q\}}(\R^P_+\times\{0\})
\ar@{^(->}[r]^{\simeq} &
\Mbarcomb_{g,P\cup\{q\}}(C_{P,q})
\ar@/^2pc/[l]^{\H^q_0} \\
& \Mbarcomb_{g,P\cup\{q\}}(C^+_{P,q})
\ar@{^(->}[u]
}
\]
In particular, we can consider
the functional $\ol{\omega}_q^{r+1}|^*_{\tilde{l}}$
associated to the restriction
of $\ol{\omega}_q^{r+1}$ to
$\Mbarcomb_{g,P\cup\{q\}}(\tilde{l})$
with $\tilde{l}\in C^+_{P,q}$
as a homology class in
$H_*(\Mbarcomb_{g,P\cup\{q\}}(\tilde{l})\setminus D_q^{comb})$
and then we can look at its image in
$H_*(\Mbarcomb_{g,P\cup\{q\}}(\R_+^P\times\{0\})\setminus D_q^{comb} )$
through $(\H^q_0)_*$.

Henceforth, with the aid of this last diagram
\[
\xymatrix@!C=2.6cm{
H_*(\Mbar_{g,P\cup\{q\}}\setminus D_q)
\ar[rr]^{\tilde{\xi}_{\tilde{l}} \qquad} \ar[d] & &
H_*(\Mbarcomb_{g,P\cup\{q\}}(C_{P,q})\setminus D^{comb}_q))
\ar[d]^{(\H^q_0)_*}_{\cong} \\
H_*(\Mbar_{g,P\cup\{q\}}\setminus D_q,
\S^q_{g,P}\setminus D_q)
\ar[rd]^{\cong} & &
H_*(\Mbarcomb_{g,P\cup\{q\}}(\R_+^P\times\{0\})\setminus D^{comb}_q)
\ar[ld] \\
& H_*(\Mbarcomb_{g,P\cup\{q\}}(\R_+^P\times\{0\})\setminus D^{comb}_q,
\Scomb_{g,P\cup\{q\}}(\R^P_+\times\{0\})\setminus D^{comb}_q)
}
\]
we obtain that, if we prove that
$(\H^q_0)_*[\ol{\omega}_q^{r+1}|^*_{\tilde{l}}]=X$
in
$H_*(\Mbarcomb_{g,P\cup\{q\}}(\R^P_+\times\{0\})\setminus D^{comb}_q,
\Scomb_{g,P\cup\{q\}}(\R^P_+\times\{0\})\setminus D^{comb}_q)$
with $X$ a combinatorial cycle, then
we have that $(\psi_q^{r+1})^*=X$ in
$H_*(\Mbar_{g,P\cup\{q\}}\setminus D_q,
\S^q_{g,P}\setminus D_q)$.
\end{section}
%
\begin{section}{Classes with one nontrivalent vertex}
\setcounter{mytheorem}{0}
\setcounter{mycorollary}{0}
Let $P=\{p_1,\dots,p_n\}$ be a nonempty set of markings
and $g\geq 0$ an integer such that $2g-2+n>0$.
For every integer $r \geq -1$
denote by $\Wbar^q_{2r+3}$ the generalized combinatorial class on
$\Mbarcomb_{g,P\cup\{q\}}$ corresponding to ribbon graphs whose
vertices are all trivalent except one
which has valency $2r+3$ and is marked by $q$. Analogously,
for every $r\geq 0$ call
$\Wbar_{2r+3}$ the combinatorial class on
$\Mbarcomb_{g,P}$ correspondent to ribbon graphs
whose vertices are all trivalent except one
which has valency $2r+3$ (in the case $r=0$ all the vertices are
trivalent). In what follows, when there is no risk of
ambiguity, we will use the same symbols
$\Wbar^q_{2r+3}$ and $\Wbar_{2r+3}$ for the combinatorial classes
and for the subcomplexes of the combinatorial moduli spaces
the classes are supported on.

Before stating the first result, we need to introduce
some auxiliary combinatorial classes.
Call $\ol{N}^q_{v_1,v_2}$
the subcomplex of $\Mbarcomb_{g,P\cup\{q\}}(\{l_q=0\})$
that parametrizes surfaces with a nonpositive
component $S_0$ of genus $0$,
such that: $S_0$ is $q$-marked, it has two singular points
and the two nodes
identify points of $S_0$ to vertices of
valencies $v_1$ and $v_2$.
Call $\ol{N}_{v_1,v_2}$ the subcomplex of
$\Mbarcomb_{g,P}$ which parametrizes surfaces
with a node that identifies a vertex of valency $v_1$
and a vertex of valency $v_2$.

It is not difficult to check that,
for odd $v_1$ and $v_2$, these subcomplexes
define cycles, which in fact could also be obtained
as images of combinatorial classes of the type
$\Wbar$ through ``combinatorial boundary maps''.
\begin{mytheorem}
For any $g$ and $n\geq 1$, the equality
\[
\Wbar^q_{2r+3}+
\sum_{\substack{i,j\geq 0 \\ i+j=r-1}} (2i+1)(2j+1)
\ol{N}^q_{2i+1,2j+1}
=\frac{(2r+2)!}{(r+1)!}(\psi_q^{r+1})^*
\]
holds in $H_{2s}(\Mbar_{g,P\cup\{q\}},
\S^q_{g,P})$ for every $r\geq -1$,
where $s=3g-3+n-r$.
As a consequence, for $r\geq 1$ the equality
\[
\Wbar_{2r+3}+
\sum_{\substack{i,j\geq 0 \\ i+j=r-1}} (2i+1)(2j+1)
\ol{N}_{2i+1,2j+1}
=2^{r+1}(2r+1)!!\,\k^*_r
\]
holds in $H_{2s}(\Mbar_{g,P},\pi_q(\S^q_{g,P}))$.
\end{mytheorem}
\begin{strategy}
First, notice that for $r\geq 0$ and $\tilde{l}=(l,0)\in\R^P_+\times\{0\}$
the cycles $\Wbar^q_{2r+3}(\tilde{l})$
and $\ol{N}^q_{2i+1,2j+1}(\tilde{l})$ live also in
$H_{2s}(\Mbarcomb_{g,P\cup\{q\}}(\tilde{l})\setminus D^{comb}_q(\tilde{l}))$,
because $\Wbar^q_{2r+3}(\tilde{l})\cap D^{comb}_q(\tilde{l})=\emptyset$
and $\ol{N}^q_{2i+1,2j+1}(\tilde{l})\cap D^{comb}_q(\tilde{l})=\emptyset$.

Then, considerations developed in Section \ref{section:pi}
show that it is sufficient to compare
\[
\frac{(2r+2)!}{(r+1)!}(\psi_q^{r+1})^* \qquad \text{and} \qquad
\Wbar^q_{2r+3}+\sum_{\substack{i,j\geq 0 \\ i+j=r-1}} (2i+1)(2j+1)
\ol{N}^q_{2i+1,2j+1}
\]
as elements of
$H_{2s}(\Mbarcomb_{g,P\cup\{q\}}(\R^P_+\times\{0\})\setminus D^{comb}_q,
\Scomb_{g,P\cup\{q\}}(\R^P_+\times\{0\})\setminus D^{comb}_q)$;
so we couple them with classes
\[
[\eta] \in
H^{2s}(\Mbarcomb_{g,P\cup\{q\}}(\R^P_+\times\{0\})\setminus D^{comb}_q,
\Scomb_{g,P\cup\{q\}}(\R^P_+\times\{0\})\setminus D^{comb}_q)
\]
and a simple computation proves the first relation.

For the second claim, it is sufficient to recall that
we can push the first relation forward, using the following
diagram
\[
\footnotesize
\xymatrix{
H_{2s}(\Mbar_{g,P\cup\{q\}}\setminus D_q,\S^q_{g,P}\setminus D_q)
\ar[d]^{\pi_{q,*}} \ar[r]^{\cong\qquad\qquad} &
H_{2s}(\Mbarcomb_{g,P\cup\{q\}}(\tilde{l})\setminus D^{comb}_q,
\Scomb_{g,P\cup\{q\}}(\tilde{l})\setminus D^{comb}_q)
\ar[d]^{\pi_{q,*}^{comb}} \\
H_{2s}(\Mbar_{g,P},\pi_q(\S^q_{g,P})) \ar[r]^{\cong\qquad} &
H_{2s}(\Mbarcomb_{g,P}(l),\pi_q^{comb}(\Scomb_{g,P})(l))
}
\]
in order to obtain the second relation
\[
\pi^{comb}_{q,*}\Big(\Wbar^q_{2r+3}
+  \!\!\! \sum_{\substack{i,j\geq 0 \\ i+j=r-1}} (2i+1)(2j+1)
\ol{N}^q_{2i+1,2j+1}\Big)
=\frac{(2r+2)!}{(r+1)!}\pi_{q,*}(\psi_q^{r+1})^*
\]
thus concluding the argument.
\end{strategy}
\begin{proof}[Proof of Theorem \ref{th:first}]
Consider a closed PL differential form $\eta$ of degree $2s$ on
$(\Mbarcomb_{g,P\cup\{q\}}\setminus D^{comb}_q)(\R_+^P\times\{0\})$
that vanishes on
$(\Scomb_{g,P\cup\{q\}}\setminus D^{comb}_q)(\R_+^P\times\{0\})$.

The form
$(\H^q_0)^*\eta\wedge\ol{\omega}_q^{r+1}$
extends by zero to $\Mbarcomb_{g,P\cup\{q\}}(C^+_{P,q})$,
because
\[
(\H^q_0)^{-1}(D^{comb}_q) \cap \Mbarcomb_{g,P\cup\{q\}}(C^+_{P,q})
\subset \ol{Y}_1(C^+_{P,q})
\]
is an inclusion in a closed neighbourhood
and $\ol{\omega}_q^{r+1}$ vanishes on $\ol{Y}_1$.
Moreover, for simple reasons of degree,
the restriction of
$(\H^q_0)^*\eta\wedge\ol{\omega}_q^{r+1}$
to $\Mbarcomb_{g,P\cup\{q\}}(\tilde{l})$
has support contained inside $\ol{Y}_{2r+3}(\tilde{l})$
for every $\tilde{l}\in C_{P,q}^+$.
In fact, the restriction of
$\ol{\omega}_q^{r+1}$ has support inside
$\ol{Y}_{\geq 2r+3}(\tilde{l})$, while
the restriction of $(\H^q_0)^*\eta$ has support
inside $\ol{Y}_{\leq 2r+3}(\tilde{l})$.
\begin{figure}[h]
\resizebox{10cm}{!}{\input{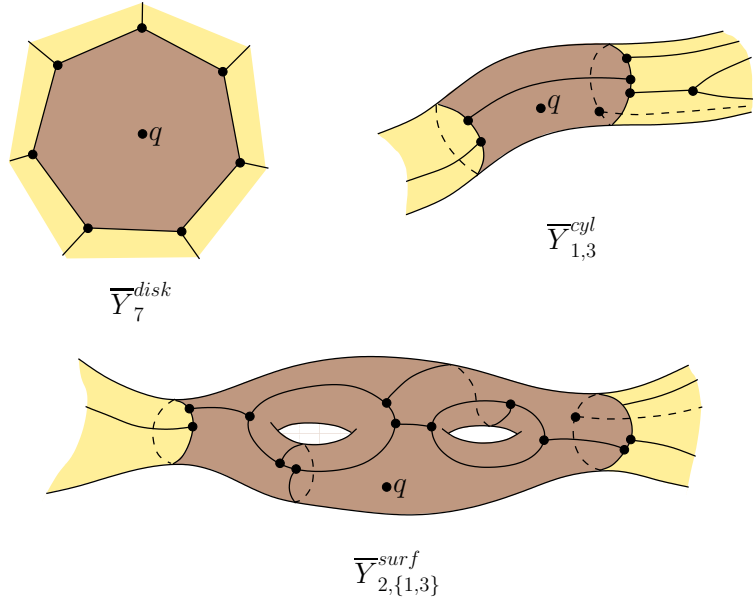}}
\caption{Three examples of loci $\ol{Y}$}
\label{fig:Yloci}
\end{figure}
Now decompose $\ol{Y}_{2r+3}(C_{P,q}^+)$ into three families of subsets:
\begin{enumerate}
\item
the closure $\ol{Y}^{disk}_{2r+3}(C_{P,q}^+)$ of the locus of graphs where
the surface $\bar{T}_q$ is a disk (see Subsection \ref{subsec:arc});
in this case $\H^q_0(\ol{Y}^{disk}_{2r+3}(C_{P,q}^+))$ is exactly
$\Wbar^q_{2r+3}(\R_+^P\times\{0\})$
\item
the closure $\ol{Y}^{cyl}_{v_1,v_2}(C_{P,q}^+)$ of the locus of
graphs where $\bar{T}_q$ is a cylinder
with exactly one internal edge $e$, which divides the other
edges of $x(q)$ into
two subsets of cardinality $v_1+1$ and $v_2+1=2r-v_1+1$;
its image via $\H^q_0$
is exactly $\ol{N}^q_{v_1,v_2}(\R_+^P\times\{0\})$
\item
the closure $\ol{Y}^{surf}_{h,\{v_1,\dots,v_{\nu}\}}(C_{P,q}^+)$ of the locus
of graphs where $\bar{T}_q$ is a surface of genus $h$ with $\nu\geq 1$
(if $h=0$, then $\nu\geq 3$)
boundary components
which touch $v_1,\dots,v_{\nu}$ external edges
(i.e. not in $\bar{T}_q$) respectively,
where $6h-6+\sum_{j=1}^{\nu}(v_j+3)=2r$;
its image via $\H^q_0$ is the locus
$\ol{Z}^q_{h,\{v_1,\dots,v_{\nu}\}}(\R_+^P\times\{0\})$ of graphs
with one nonpositive component of genus $h$ which
has the $q$-marking and $\nu$ nodes corresponding to vertices
of valencies $v_1,\dots,v_{\nu}$.
\end{enumerate}
Notice that all the subcomplexes $\ol{Z}^q_{h,v_*}(\R_+^P\times\{0\})$
are contained inside $\Scomb_{g,P\cup\{q\}}(\R_+^P\times\{0\})$,
hence $(\H^q_0)^*\eta$ vanishes on
$\ol{Y}^{surf}_{h,v_*}(C_{P,q}^+)$.

Remark that $\bar{\l}_{p_i}(\H^q_0(\ol{Y}_{2r+3}(\tilde{l})))$ takes
values between $\tilde{l}_{p_i}-\tilde{l}_q$ and $\tilde{l}_{p_i}$.
Then, choose $0<\e<L''\ll L'$
and notice that $\H^q_0(\ol{Y}_{2r+3}([L'',L']^n\times\{\e\}))$ contains
\[
\left(
\cup_{h,v_*}\ol{Z}^q_{h,v_*}\cup\Wbar^q_{2r+3}
\cup_{v_1,v_2} \ol{N}^q_{v_1,v_2}
\right)
([L'',L'-\e]^n\times\{0\})
\]
and is contained inside
\[
\left(
\cup_{h,v_*}\ol{Z}^q_{h,v_*}\cup\Wbar^q_{2r+3}
\cup_{v_1,v_2} \ol{N}^q_{v_1,v_2}
\right)
([L''-\e,L']^n\times\{0\}).
\]
Since the volume of the difference $[L''-\e,L']^n\setminus[L'',L'-\e]^n$
goes to zero as $\e$ decreases, we have
\begin{align*}
& \int_{[L'',L']^n} d\tilde{l}_{p_1} \wedge\dots\wedge d\tilde{l}_{p_n}
\int_{\Mbar_{g,P\cup\{q\} }} \psi_q^{r+1}\smallsmile
\tilde{\xi}^*(\H^q_0)^*[\eta] = \\
& =\lim_{\e \rar 0}
\int_{\ol{Y}_{2r+3}([L'',L']^n\times\{\e\})}
\bar{\l}^*(d\tilde{l}_{p_1}\wedge\dots \wedge d\tilde{l}_{p_n})\wedge
\ol{\omega}_q^{r+1}\wedge (\H^q_0)^*\eta = \\
& =\lim_{\e\rar 0} \Big(
\int_{\ol{Y}^{disk}_{2r+3}([L'',L']^n\times\{\e\})}
\bar{\l}^*(d\tilde{l}_{p_1}\wedge\dots d\tilde{l}_{p_n})\wedge
\ol{\omega}_q^{r+1}\wedge(\H^q_0)^*\eta+ \\
&\quad +\sum_{v_1+v_2=2r}\int_{\ol{Y}^{cyl}_{v_1,v_2}([L'',L']^n\times\{\e\})}
\bar{\l}^*(d\tilde{l}_{p_1}\wedge\dots d\tilde{l}_{p_n})\wedge
\ol{\omega}_q^{r+1}\wedge(\H^q_0)^*\eta+ \\
&\quad +\sum_{h,v_*}\int_{\ol{Y}^{surf}_{h,v_*}([L'',L']^n\times\{\e\})}
\bar{\l}^*(d\tilde{l}_{p_1}\wedge\dots d\tilde{l}_{p_n})\wedge
\ol{\omega}_q^{r+1}\wedge(\H^q_0)^*\eta  \Big) = \\
& =\lim_{\e \rar 0} \Big(
\int_{\Wbar^q_{2r+3}([L'',L'-\e]^n\times\{0\})}
\bar{\l}^*(d\tilde{l}_{p_1}\wedge\dots d\tilde{l}_{p_n})\wedge\eta
\int_{F_{2r+3}^{disk}(\e)} \ol{\omega}_q^{r+1}+ \\
&\quad +\sum_{v_1+v_2=2r}\int_{\ol{N}^q_{v_1,v_2}([L'',L'-\e]^n\times\{0\})}
\bar{\l}^*(d\tilde{l}_{p_1}\wedge\dots d\tilde{l}_{p_n})\wedge\eta
\int_{F^{cyl}_{v_1,v_2}(\e)} \ol{\omega}_q^{r+1} \Big) = \\
& =\int_{[L'',L']^n}d\tilde{l}_{p_1}\wedge\dots\wedge d\tilde{l}_{p_n} \Big(
\int_{\Wbar^q_{2r+3}(\tilde{l})}\eta
\int_{F_{2r+3}^{disk}(\e)}\ol{\omega}_q^{r+1} + \\
& \quad + \sum_{v_1+v_2=2r} \int_{\ol{N}^q_{v_1,v_2}(\tilde{l})}\eta
\int_{F^{cyl}_{v_1,v_2}(\e)}\ol{\omega}_q^{r+1} \Big)
\end{align*}
where $\tilde{l}$ belongs to $\R_+^P\times\{0\}$ and
$F_{2r+3}^{disk}(\e)$ is the intersection of the generic fiber
of $\H^q_0$ over $\Wbar^q_{2r+3}([L'',L'-\e]^n\times\{0\})$
with $\ol{Y}_{2r+3}([L'',L']^n\times\{\e\})$ and similarly
for $F^{cyl}_{v_1,v_2}$.
\begin{remark}
In the first equality,
we used the push-forward
through the map
\[
\tilde{\xi}_\e:\Mbar_{g,P\cup\{q\}}\times\R_+^P\times\{\e\}
\lra \Mbarcomb_{g,P\cup\{q\}}(\R_+^P\times\{\e\}).
\]
In the third equality, we used that $\H^q_0$ restricts to
\[
\xymatrix{
\ol{Y}_{2r+3}^{disk}([L'',L']^n\times\{\e\}) \ar[r]^{\H^q_0} &
\Wbar^q_{2r+3}([L''-\e,L']^n\times\{0\}) \\
(\H^q_0)^{-1}(\Wbar^q_{2r+3}([L'',L'-\e]^n\times\{0\})) \ar@{^(->}[u]
\ar@{->>}[r] & \Wbar^q_{2r+3}([L'',L'-\e]^n\times\{0\}) \ar@{^(->}[u]
}
\]
where the lower map is a fibration with fiber $F_{2r+3}^{disk}(\e)$,
and the fact that the volumes of
the differences $\ol{Y}_{2r+3}^{disk}([L'',L']^n\times\{\e\})\setminus
(\H^q_0)^{-1}(\Wbar^q_{2r+3}([L'',L'-\e]^n\times\{0\}))$
and $\Wbar^q_{2r+3}([L''-\e,L']^n\times\{0\})\setminus
\Wbar^q_{2r+3}([L'',L'-\e]^n\times\{0\})$ tend to zero with $\e$.
Analogous considerations hold for $\ol{Y}_{v_1,v_2}^{cyl}$
and $\ol{N}^q_{v_1,v_2}$.
\end{remark}
It is easy to see that $F_{2r+3}^{disk}(\e)$
is a simplex of dimension $2r+2$ with affine
coordinates $e_0,\dots,e_{2r+2}$
subject to the constraint
$\sum_{j=0}^{2r+2} e_j=2\e$,
where $e_j$ are the (unoriented) edges of the $q$-marked hole.
It is also immediate to see that $\ol{\omega}_q^{r+1}$ is equal to
$(r+1)!d(\frac{e_1}{2\e})\wedge\dots\wedge d(\frac{e_{2r+2}}{2\e})$
on $F_{2r+3}^{disk}(\e)$,
so that
\[
\int_{F_{2r+3}^{disk}(\e)} \ol{\omega}_q^{r+1}=
(r+1)!\mathrm{vol}(\D_{2r+2})=
\frac{(r+1)!}{(2r+2)!} \, .
\]

A simple computation shows that $\ol{\omega}_q^{r+1}$ vanishes
on $F^{cyl}_{v_1,v_2}$ if $v_1$ and $v_2$ are even; while
$\ol{\omega}_q^{r+1}=(r+1)!d(\frac{e_1}{2\e})\wedge\dots
\wedge d(\frac{e_{2r+2}}{2\e})$ if $v_1$ and $v_2$ are odd,
where $2e_0+\sum_{j=1}^{2r+2} e_j=2\e$ and $e_0$ is the ``separating''
edge of the cylinder. We deduce that, for $v_1$ and $v_2$ odd,
\[
\int_{F^{cyl}_{v_1,v_2}(\e)} \ol{\omega}_q^{r+1}=v_1 v_2\frac{(r+1)!}{(2r+2)!}
\]
because $F^{cyl}_{v_1,v_2}$ contains $v_1 v_2$ top-dimensional simplices.

Hence, we conclude that
\begin{multline*}
\frac{(2r+2)!}{(r+1)!}
\int_{\Mbar_{g,P\cup\{q\} }} \psi_q^{r+1}\smallsmile
\tilde{\xi}^*(\H^q_0)^*[\eta]=
\int_{\Wbar^q_{2r+3}(\tilde{l})}\eta+ \\
+\sum_{\substack{i,j\geq 0 \\ i+j=r-1}} (2i+1)(2j+1)
\int_{\ol{N}^q_{2i+1,2j+1}(\tilde{l})}\eta
\end{multline*}
and the proof is complete.
\end{proof}
\begin{mycorollary}
For every $g$ and $|P|=n\geq 1$ such that $2g-2+n>0$
the following equalities hold
\[
\begin{array}{ll}
\displaystyle{
\Wbar^q_5+\d^q_{irr}+
\sum_{\substack{g',I\neq \emptyset,P}}  \d^q_{g',I}
=12(\psi^2_q)^* }
& \quad\text{in $H_{6g-8+2n}(\Mbar_{g,P\cup\{q\}},\S^q_{g,P})$} \\
\displaystyle{
\Wbar_5+\d_{irr}+
\sum_{\substack{g',I\neq\emptyset,P}}  \d_{g',I}
=12 \k^*_1 }
& \quad\text{in $H_{6g-8+2n}(\Mbar_{g,P},\S_{g,P})$}
\end{array}
\]
where $\d^q_{g',I}$ is the image of the morphism
\[
\Mbar_{g',I\cup\{p'\}} \times\Mbar_{0,\{q,q',q''\}}
\times\Mbar_{g-g',I^c\cup\{p''\}} \rar \Mbar_{g,P\cup\{q\}}
\]
that glues $p'$ with $q'$ and $p''$ with $q''$
(analogously for $\d^q_{irr}$).
\end{mycorollary}
The second equality of the
previous corollary has been proven first by
Arbarello and Cornalba \cite{arbarello-cornalba:combinatorial}
in a very different manner.
Here it is a consequence
of the proof of Theorem \ref{th:first}
and of the remark
at the end of Section \ref{section:combinatorial}.
In fact, the difference between $\pi_q(\S^q_{g,P})$
and $\S_{g,P}$ has complex codimension three.
Moreover, for $r=1$ all simplices of top dimension in
$\ol{Y}^{cyl}$ have $v_1=v_2=1$.
\end{section}
\begin{section}{More on the shrinking map and the forgetful map}
\label{section:pi2}
We now want to examine the case of an arbitrary class
$\Wbar_{m_*,\rho,P}$ on $\Mbarcomb_{g,P\cup Q}$
for some $\rho:Q\rar\N$.

The ideas involved in the proof are the same as before,
but more care is needed.
The main difference is that $\pi_Q(\S^Q_{g,P})$
is not contained in $\pa\M_{g,P}$, as in the case $Q=\{q\}$.

However, we do not need to consider
homology classes relative to $\S^Q_{g,P}$. In fact,
the singular locus of $\Mbarcomb_{g,P\cup Q}(\R_+^P\times\{0\}^Q)$ is
smaller than $\Scomb_{g,P\cup Q}(\R_+^P\times\{0\}^Q)$.
In this way, we can adapt the combinatorial forgetful map
and the shrinking map to the case $|Q|>1$.\\

%
Fix $P=\{p_1,\dots,p_n\}$,
$Q'=\{q_1,\dots,q_s\}$ and $Q''=\{q_{s+1},\dots,q_{s+u}\}$,
and let $Q=Q'\cup Q''$.
Consider the forgetful map
$\pi_{Q''}:\Mbar_{g,P\cup Q}\lra\Mbar_{g,P\cup Q'}$
and call $\Cc_{g,P\cup Q'}^{Q''}$ the inverse
image $\pi_{Q''}^{-1}(\M_{g,P\cup Q'})$,
which is the locus of curves with one component of
geometric genus $g$ plus some rational tails,
each one containing at most one point of $P\cup Q'$;
moreover, call $\pa\Cc_{g,P\cup Q'}^{Q''}$ its boundary
$\Mbar_{g,P\cup Q}\setminus\Cc_{g,P\cup Q'}^{Q''}$.

On the combinatorial side, consider the map
\[
\tilde{\xi}: \Mbar_{g,P\cup Q}\times\R_+^P\times\{0\}^Q
\lra \Mbarcomb_{g,P\cup Q}(\R_+^P\times\{0\}^Q)
\]
and call $\Cc^{Q'',comb}_{g,P\cup Q'}(\R_+^P\times\{0\}^Q)$
the image of $\Cc^{Q''}_{g,P\cup Q'}\times\R_+^P\times\{0\}^Q$,
and analogously $\pa\Cc^{Q'',comb}_{g,P\cup Q'}(\R_+^P\times\{0\}^Q)$
its boundary
$(\Mbarcomb_{g,P\cup Q}\setminus\Cc^{Q'',comb}_{g,P\cup Q'})
(\R_+^P\times\{0\}^Q)$.

Remember that $D_{Q'',P}\subset\Mbar_{g,P\cup Q}$
is the union of all divisors of the type
$\d_{0,\{ q_{j_1},\dots,q_{j_h},p_i \} }$
for $p_i\in P$ and $\{ q_{j_1},\dots,q_{j_h} \} \subset Q''$.
Call $D^{comb}_{Q'',P}$ its image through $\tilde{\xi}$ and
define $\pi^{comb}_{Q''}$ as the composition
$\pi^{comb}_{q_{s+1}}\cdots\pi^{comb}_{q_{s+u}}$, so that
the following diagram
\[
\xymatrix{
(\Mbar_{g,P\cup Q}\setminus D_{Q'',P})\times \R_+^P\times\{0\}^Q
\ar[d]^{\pi_{Q''}} \ar[r]^{\tilde{\xi}} &
\Mbarcomb_{g,P\cup Q}(\R_+^P\times\{0\}^Q)\setminus D^{comb}_{Q'',P}
\ar[d]^{\pi^{comb}_{Q''}} \\
\Mbar_{g,P\cup Q'}\times \R_+^P\times\{0\}^{Q'}
\ar[r]^{\xi} &
\Mbarcomb_{g,P\cup Q'}(\R_+^P\times\{0\}^{Q'})
}
\]
is commutative.
\begin{lemma}
The locus $(\Cc^{Q'',comb}_{g,P\cup Q'}
\setminus D^{comb}_{Q'',P})(\R^P_+\times\{0\})$
is smooth.
\end{lemma}
\begin{proof}
Notice that the restriction of $\tilde{\xi}$ factorizes as follows
\[
\xymatrix@!C=2.8cm{
(\Cc^{Q''}_{g,P\cup Q'}\setminus D_{Q'',P})\times\R^P_+\times\{0\}^Q
\ar[rr]^{\tilde{\xi}} \ar[rd] &&
(\Cc^{Q'',comb}_{g,P\cup Q'}\setminus D^{comb}_{Q'',P})
(\R_+^P\times\{0\}^Q) \\
& (\Cc^{q_{s+1}}_{g,P\cup Q'}\setminus D_{\{q_{s+1}\},P})
\times_{_{\M_{g,P\cup Q'}}} \dots_{_{\M_{g,P\cup Q'}}} \times
(\Cc^{q_{s+u}}_{g,P\cup Q'}\setminus D_{\{q_{s+u}\},P})
\times\R^P_+\times\{0\}^Q
\ar[ru]_{\cong}
}
\]
and that the locus we are interested in is homeomorphic to
a fiber product of $u$ copies of
$\Cc^q_{g,P\cup Q'}\setminus D_{\{q\},P}$
over $\M_{g,P\cup Q}$, which is smooth.
\end{proof}
We now want to show that for every $\tilde{l}\in\R_+^P\times\{0\}^Q$
we can lift homology classes via $\tilde{\xi}$ from
\[
H_*((\Mbarcomb_{g,P\cup Q}\setminus D^{comb}_{Q'',P})(\tilde{l}),
(\pa\Cc^{Q'',comb}_{g,P\cup Q'}\setminus D^{comb}_{Q'',P})(\tilde{l}))
\]
to
\[
H_*((\Mbar_{g,P\cup Q}\setminus D_{Q'',P})\times\{\tilde{l}\},
(\pa\Cc^{Q''}_{g,P\cup Q'}\setminus D_{Q'',P})\times\{\tilde{l}\})
\]
using a sort of Lefschetz duality.
\begin{remark}
Let be given a compact connected triangulated space $K$
and two nonempty proper subcomplexes $F$ and $G$, such that 
$K\setminus (F\cup G)$
is a connected oriented topological manifold of dimension $d$.
Suppose that $F$ (resp. $G$) has a closed neighbourhood  
$N_F$ (resp. $N_G$) that retracts on $F$ (resp. on $G$)
by deformation. Moreover, suppose that $\pa N_F\setminus N_G$ is
a topological manifold of dimension $d-1$.
Then we can define a duality homomorphism
\[
H_k(K\setminus N_F, N_G\setminus N_F)
\lra
H^{d-k}(K\setminus(N_F^o\cup N_G),\pa N_F\setminus N_G) 
\]
noticing that every dual cocycle
vanishes on $\pa N_F\setminus N_G$.
Moreover, this map fits in the following diagram
\[
\scriptsize
\xymatrix@C=0.2cm{
H^{d-k-1}(\pa N_F\setminus N_G) \ar[r] &
H^{d-k}(K\setminus(N_F^o\cup N_G),\pa N_F\setminus N_G)
\ar[r] &
H^{d-k}(K\setminus(N_F\cup N_G)) \\
H_k(\pa N_F,\pa N_F\cap N_G) \ar[u]^{\cong} \ar[r] &
H_k(K\setminus N_F^o, N_G\setminus N_F^o) \ar[u] \ar[r] &
H_k(K\setminus N_F^o,\pa N_F\cup (N_G\setminus N_F^o))
\ar[u]^{\cong}
}
\]
where the first row is the exact sequence of the couple
\[
(K\setminus(N_F^o\cup N_G),\pa N_F\setminus N_G)
\]
and the second row is the exact sequence of the triple
\[
(K\setminus N_F^o,\pa N_F\cup (N_G\setminus N_F^o),
N_G\setminus N_F^o),
\]
and we are using the homotopy equivalences
\begin{align*}
K\setminus (N_F^o\cup N_G) & \simeq
K\setminus (N_F\cup N_G) \\
(K\setminus N_F, N_G\setminus N_F) & \simeq
(K\setminus N_F^o, N_G\setminus N_F^o).
\end{align*}
Notice that we are identifying
\[
H_k(\pa N_F,\pa N_F\cap N_G) \qquad\text{and}\qquad
H_k(\pa N_F\cup (N_G\setminus N_F^o),N_G\setminus N_F^o)
\]
in the following way.
First we have that
\[
H_k(\pa N_F,\pa N_F \cap N_G) \cong
H_k(\pa N_F\setminus G, \pa N_F \cap (N_G \setminus G))
\]
by excision of $\pa N_F \cap G$.
Then we use that $\pa N_F \cap (N_G\setminus G)$
has a tubular neighbourhood inside $N_G\setminus G$
and we get
\[
H_k(\pa N_F\setminus G, \pa N_F \cap (N_G \setminus G))
\cong
H_k(\pa N_F \cup (N_G \setminus N_F^o)\setminus G,
N_G\setminus (N_F^o\cup G)).
\]
Finally, the excision of $G\setminus N_F^o$ gives
\[
H_k(\pa N_F \cup (N_G \setminus N_F^o)\setminus G,
N_G\setminus (N_F^o\cup G))
\cong
H_k(\pa N_F \cup (N_G \setminus N_F^o),
N_G\setminus N_F^o)
\]
which establishes the desired identification.
Hence, we can conclude that the vertical arrow
in the middle is an isomorphism by the Five Lemma.
\end{remark}
Now we want to show that the hypotheses of
the previous remark are satisfied in two particular
situations.
In the first case, $F=D_{Q'',P}$, $G=\pa\mathcal{C}^{Q''}_{g,P\cup Q'}$
and $K=\Mbar_{g,P\cup Q}$, so we can immediately conclude
because $K$ is smooth and
$F$ is a divisor with normal crossings
(as $\Q$-varieties).
In the second case,
$F=D^{comb}_{Q'',P}(\tilde{l})$,
$G=\pa\mathcal{C}^{Q'',comb}_{g,P\cup Q'}(\tilde{l})$
and $K=\Mbarcomb_{g,P\cup Q}(\tilde{l})$, and
we need to prove that there exist ``good''
neighbourhoods $N_F$ and $N_G$.
Call $\widetilde{\mathcal{C}}^{Q''}_{g,P\cup Q'}$
the fiber product of all
$\Mbar_{g,P\cup Q'\cup\{q_{s+i}\}}$
for $i=1,\dots,u$ over $\Mbar_{g,P\cup Q'}$.
The continuous map
\[
\prod^{i=1,\dots,u}_{\M_{g,P\cup Q'}}
\mathcal{C}^{q_{s+i}}_{g,P\cup Q'}\setminus
D_{\{q_{s+i}\},P}
\longrightarrow
\Mbarcomb_{g,P\cup Q}(\tilde{l})
\]
does not extend to
$\widetilde{\mathcal{C}}^{Q''}_{g,P\cup Q'}$.
However, there is a well-defined map
$\Mbar_{g,P\cup Q} \twoheadrightarrow \Mbarcomb_{g,P\cup Q}(\tilde{l})$
and an algebraic morphism
$\Mbar_{g,P\cup Q} \twoheadrightarrow
\widetilde{\mathcal{C}}^{Q''}_{g,P\cup Q'}$.
As a consequence,
even though $\Mbarcomb_{g,P\cup Q}(l)$
may not be an algebraic space,
there is a resolution
$\widehat{\mathcal{C}}^{Q''}_{g,P\cup Q'}$,
obtained from
$\widetilde{\mathcal{C}}^{Q''}_{g,P\cup Q'}$
by performing iterated blow-ups
with centers outside
\[
\prod^{i=1,\dots,u}_{\M_{g,P\cup Q'}}
\mathcal{C}^{q_{s+i}}_{g,P\cup Q'}\setminus
D_{\{q_{s+i}\},P} \,\, ,
\]
such that the following diagram commutes
\[
\xymatrix{
& \Mbar_{g,P\cup Q} \ar@{->>}[r] \ar@{->>}[d]
& \Mbarcomb_{g,P\cup Q}(\tilde{l}) \\
& \widehat{\mathcal{C}}^{Q''}_{g,P\cup Q'} \ar@{->>}[d] \ar@{->>}[ur] \\
\prod\limits^{i=1,\dots,u}_{\M_{g,P\cup Q'}}
\mathcal{C}^{q_{s+i}}_{g,P\cup Q'}\setminus
D_{\{q_{s+i}\},P} \ar@{^(->}[r] \ar@{^(->}[ur]
& \widetilde{\mathcal{C}}^{Q''}_{g,P\cup Q'} \ar@{-->}[uur]
}
\]
and the inverse image
$\widehat{D}_{Q'',P}\subset \widehat{\mathcal{C}}^{Q''}_{g,P\cup Q'}$
of $D^{comb}_{Q'',P}$
is a divisor with normal crossings
inside a smooth variety. Hence,
$\widehat{D}_{Q'',P}$ has
a closed neighbourhood whose
boundary is a closed
manifold of real codimension $1$ and
we can define $N_F$ to be its
image in $\Mbarcomb_{g,P\cup Q}(\tilde{l})$.
Also the inverse image of
$\pa\mathcal{C}^{Q'',comb}_{g,P\cup Q'}(\tilde{l})$
through last map
is a divisor with normal crossings,
so we can choose
the image of an adequate neighbourhood
of this divisor as $N_G$.

Thus, we can consider a cycle
on $(\Mbarcomb_{g,P\cup Q}\setminus D^{comb}_{Q'',P})(\tilde{l})$
with boundary contained in
\hbox{$(\pa\Cc^{Q'',comb}_{g,P\cup Q'}\setminus D^{comb}_{Q'',P})(\tilde{l})$},
dualize it, pull it back and dualize it again, in order
to get a cycle on
$(\Mbar_{g,P\cup Q}\setminus D_{Q'',P})\times\{\tilde{l}\}$
relative to
$(\pa\Cc^{Q''}_{g,P\cup Q'}\setminus D_{Q'',P})\times\{\tilde{l}\}$.
This defines the desired pull-back map in homology.
In what follows we will denote
by the same symbols combinatorial classes and their lifts.\\

Remember that,
for every $a_1,\dots,a_s,b_1,\dots,b_u\geq 1$, the class
\[
\psi_{q_1}^{a_1}\cdots\psi_{q_s}^{a_s}\cdot
\psi_{q_{s+1}}^{b_1}\cdots\psi_{q_{s+u}}^{b_u}
\]
which is equal to
\[
\pi_{q_{s+u},\dots,q_2}^*(\psi_{q_1}^{a_1})     \cdot
\pi_{q_{s+u},\dots,q_3}^*(\psi_{q_2}^{a_2})     \cdots
\pi_{q_{s+u}}^*(\psi_{q_{s+u-1}}^{b_{u-1}})     \cdot
\psi_{q_{s+u}}^{b_u}
\]
in $H^*(\Mbar_{g,P\cup Q})$, admits a natural lift
to $H^*(\Mbar_{g,P\cup Q},D_{Q'',P})$ (see
Section \ref{section:tautological}).
Thus, we can lift the Poincar\'e dual of this class
from $H_*(\Mbar_{g,P\cup Q})$ to
$H_*(\Mbar_{g,P\cup Q}\setminus D_{Q'',P})$, and so
it makes sense to compare in
$H_*(\Mbar_{g,P\cup Q}\setminus D_{Q'',P},
\pa\Cc^{Q''}_{g,P\cup Q'}\setminus D_{Q'',P})$
the tautological class above
and a combinatorial class that lives in
$H_*(\Mbarcomb_{g,P\cup Q}(\tilde{l})\setminus D^{comb}_{Q'',P}(\tilde{l}))$
with $\tilde{l}\in\R_+^P\times\{0\}^Q$.

In order to exploit the representatives $\ol{\omega}$ of
the $\psi$ classes, we now define a new
shrinking map $\H^{Q}_0$, which is in fact the composition
of simple shrinking maps $\H^q_0$.

For every $k=1,\dots,s+u$, call $Q_k:=\{q_1,\dots,q_k\}$
and, analogously to Section \ref{section:pi},
call $C_{P,k}$ the subset
of $\tilde{l}\in\R_{\geq 0}^{P\cup Q}$ defined by
\[
\begin{cases}
l_{q_j}=0 & \text{for all $j=k+1,\dots,s+u$} \\
\sum_{i=j+1}^k l_{q_i}<l_{q_j} & \text{for all $j=1,\dots,k-1$} \\
\sum_{i=1}^k l_{q_i}<l_{p_j} & \text{for all $j=1,\dots,n$}
\end{cases}
\]
and set $C_{P,k}^+:=C_{P,k}\cap\{l_{q_k}>0\}$.
Define coherently $C_{P,0}:=\R_+^P\times\{0\}^Q$.

Call $\H^Q_0:\Mbarcomb_{g,P\cup Q}(C_{P,s+u})\twoheadrightarrow
\Mbarcomb_{g,P\cup Q}(C_{P,0})$
the composition $\H^{q_1}_0 \H^{q_2}_0\cdots \H^{q_{s+u}}_0$ of all
the retractions
\[
\H^{q_i}_0:\Mbarcomb_{g,P\cup Q}(C_{P,i})
\rar \Mbarcomb_{g,P\cup Q}(C_{P,i-1})
\]
and remark that
$\ol{\omega}_{q_1}^{a_1}\wedge\dots\wedge\ol{\omega}_{q_s}^{a_s}\wedge
\ol{\omega}_{q_{s+1}}^{b_1}\wedge\dots\wedge\ol{\omega}_{q_{s+u}}^{b_u}$
is cohomologous to
\[
(\H_0^{q_{s+u}}\cdots\H_0^{q_2})^*(\ol{\omega}_{q_1}^{a_1})
\wedge \dots \wedge
(\H_0^{q_{s+u}})^*       (\ol{\omega}_{q_{s+u-1}}^{b_{u-1}})  \wedge
                          \ol{\omega}_{q_{s+u}}^{b_u}
\]
on $\Mbarcomb_{g,P\cup Q}(C^+_{P,s+u})$.
As a consequence, if we call $\ol{Y}_{t_1,\dots,t_{s+u}}\subset
\Mbarcomb_{g,P\cup Q}$ the closure of the locus of graphs
where
\begin{itemize}
\item[-]
the $Q$-marked holes have positive perimeters
\item[-]
the hole marked by $q_{s+u}$ consists of $t_{s+u}$
distinct edges
\item[-]
for every $i=1,\dots,s+u-1$
the hole marked by $q_i$ consists of $t_i$ distinct
edges beside those which border the holes marked by
$q_{i+1},\dots,q_{s+u}$,
\end{itemize}
then we notice that the restriction of this last PL differential form
to $\Mbarcomb_{g,P\cup Q}(\tilde{l})$ with $\tilde{l}\in C_{P,s+u}^+$
has support contained in
$\ol{Y}_{\geq 2,\dots,\geq 2}(C_{P,s+u}^+)$.
By a simple computation of Euler characteristic,
we immediately see that
$D^{comb}_{Q'',P}(\tilde{l})\cap\ol{Y}_{\geq 2,\dots,\geq 2}(\tilde{l})
=\emptyset$ and so the previous PL differential form vanishes on
$D^{comb}_{Q'',P}(\tilde{l})$.
Thus, it produces a class in
$H_*((\Mbarcomb_{g,P\cup Q}\setminus D^{comb}_{Q'',P})(\tilde{l}))$,
which we can push forward to
$H_*((\Mbarcomb_{g,P\cup Q}\setminus D^{comb}_{Q'',P})(\R_+^P\times\{0\}^Q))$
through $\H^Q_0$.
In this way, we obtain the required representative
for a tautological class to compare with a combinatorial class.

\vspace{0.5cm}
\hspace{-2cm}
\mbox{
\scriptsize
\xymatrix@C=0.3cm{
\\
& \Mbarcomb_{g,P\cup Q}(\tilde{l})\setminus D^{comb}_{Q'',P}
\ar[ld] \ar[d]
%
%
\POS(40,-18);(35,-25)**\dir{-}
\POS(26,-37);(22,-42)**\dir{-}
\POS(18,-47);(0,-70)**\dir{-}?>*\dir{>}
\POS(6,-58)*{\pi^{comb}_{Q''}}
%
& \Mbar_{g,P\cup Q}\setminus D_{Q'',P}
\ar[l] \ar[d] \ar[dr] \\
\Mbarcomb_{g,P\cup Q}(\tilde{l})
\ar[d]
& (\Mbarcomb_{g,P\cup Q},
\pa\Cc^{Q'',comb}_{g,P\cup Q'})(\tilde{l})\setminus D^{comb}_{Q'',P}
\ar[ld] \ar@{-}[d] 
& (\Mbar_{g,P\cup Q},
\pa\Cc^{Q''}_{g,P\cup Q'})\setminus D_{Q'',P}
\ar[l] \ar[d]
& \Mbar_{g,P\cup Q}
\ar[ld] \ar[ddd]^{\pi_{Q''}} \ar@/_8pc/[lll]^{\tilde{\xi}_l} \\
(\Mbarcomb_{g,P\cup Q},
\pa\Cc^{Q'',comb}_{g,P\cup Q'})(\tilde{l})
& \ar[d]
& (\Mbar_{g,P\cup Q},
\pa\Cc^{Q''}_{g,P\cup Q'})
\ar[ll] \ar[d] \\
& (\Mbarcomb_{g,P\cup Q'},\pa\Mcomb_{g,P\cup Q'})(l)
& (\Mbar_{g,P\cup Q'},\pa\M_{g,P\cup Q'})
\ar[l]_{\quad\cong} \\
\Mbarcomb_{g,P\cup Q'}(l)
\ar[ru]
&&&  \Mbar_{g,P\cup Q'}
\ar[lu] \ar[lll]_{\qquad\xi_l}
}
}
\newline

Now, look at the big diagram above with $l\in\R_+^P\times\{0\}^{Q'}$
and $\tilde{l}=(l,0)\in\R_+^P\times\{0\}^Q$, and suppose
we have already proven an equality of the type
\[
X=(\H^Q_0)_*\left[(\H_0^{q_{s+u}}\cdots\H_0^{q_2})^*
(\ol{\omega}_{q_1}^{a_1}) \wedge\dots\wedge
(\H_0^{q_{s+u}})^*(\ol{\omega}_{q_{s+u-1}}^{b_{u-1}}) \wedge
\ol{\omega}_{q_{s+u}}^{b_u}\right]^*
\]
in $H_*((\Mbarcomb_{g,P\cup Q},
\pa\Cc^{Q'',comb}_{g,P\cup Q'})(\R_+^P\times\{0\}^{Q})
\setminus D^{comb}_{Q'',P})$,
where $X$ is a combinatorial class.
Then this equality lifts to
\[
X=(\psi_{q_1}^{a_1}\cdots\psi_{q_{s+u}}^{b_u})^*
\]
in $H_*(\Mbar_{g,P\cup Q},\pa\Cc^{Q''}_{g,P\cup Q'})$.
Moreover, we can push it forward to get
\[
(\pi_{Q''}^{comb})_*(X)=
(\pi_{Q''})_*(\psi_{q_1}^{a_1}\cdots\psi_{q_{s+u}}^{b_u})^*
\]
in $H_*(\Mbar_{g,P\cup Q'},\pa\M_{g,P\cup Q'})$.
\end{section}
%
%
\begin{section}{Classes with many nontrivalent vertices}\label{sec:second}
\setcounter{mycorollary}{0}
As in the previous section,
fix $P=\{p_1,\dots,p_n\}$,
$Q'=\{q_1,\dots,q_s\}$ and $Q''=\{q_{s+1},\dots,q_{s+u}\}$,
and let $Q=Q'\cup Q''$.
Moreover, let be given $m_*=(0,m_0,m_1,\dots)$ and
$\rho:Q\lra \N$ in such a way that $m_i=m_i^{\rho}$
for every $i>0$ and $m_0 \geq m_0^{\rho}$.
Clearly, one must have $4g-4+2|P|=\sum_{i\geq 0} (2i+1) m_i$.
In what follows we only consider combinatorial classes
$\Wbar_{m_*,\rho,P}$ with $\rho$ taking nonnegative values, because proofs
and notations are simpler in this case.

As we are interested in cycles
in $(\Mbarcomb_{g,P\cup Q}\setminus D^{comb}_{Q'',P},
\pa\Cc^{Q'',comb}_{g,P\cup Q'}\setminus D_{Q'',P}^{comb})$
that are sent through $\pi_{Q''}^{comb}$
to combinatorial cycles
in $(\Mbarcomb_{g,P\cup Q'},\pa\Mcomb_{g,P\cup Q'})$,
now we introduce the last type of combinatorial classes
$\Wbar^{rt}_{M,\rho,P}$ we will have to deal with.

\begin{notation}
We denote by $\mathfrak{P}_{Q}$ the
set of partitions of $Q$.
We denote by $\mathfrak{P}_{Q,Q'}$ the subset of
$\mathfrak{P}_{Q}$ consisting of all
$M=\{\mu_1,\dots,\mu_k\}$
such that the restriction $M\cap Q':=\{\mu_1\cap Q',\dots,\mu_k\cap Q'\}$
is the discrete partition of $Q'$.
For every $q\in Q'$ we will denote by
$\mu_q$ the element of $M$ that contains $q$.
\end{notation}

Let $M=\{\mu_1,\dots,\mu_k\}\in\mathfrak{P}_{Q}$
be a partition of $Q$, and set $I_M:=\{\i_{\mu_1},\dots,\i_{\mu_k}\}$.
For $l\in\R_{\geq 0}^P\setminus\{0\}$
and $\tilde{l}=(l,0)\in \R_{\geq 0}^P\times\{0\}^Q$, consider the
following commutative diagram
\[
\xymatrix{
\Mbar_{g,P\cup I_M}\times
\prod_{\mu\in M}\ol{M}_{0,\mu\cup\{\i'_\mu\}}
\ar[r]^{\qquad\qquad\vartheta^{rt}_{M,P}} \ar[d] &
\Mbar_{g,P\cup Q} \ar[d] \\
\Mbarcomb_{g,P\cup I_M}(l,0)
\ar[r]^{\qquad\vartheta^{rt,comb}_{M,P}} &
\Mbarcomb_{g,P\cup Q}(\tilde{l})
}
\]
where the map $\vartheta^{rt}_{M,P}$
glues the points $\i_\mu$ with $\i'_{\mu}$ together
for every $\mu\in M$ and, if necessary, stabilizes the curve;
while $\vartheta^{rt,comb}_{M,P}$
glues a nonpositive $\mu$-marked sphere on $\i_{\mu}$
and, if necessary, stabilizes the surface and
reduces the resulting dual graph.
Both the horizontal maps are closed immersions.
Remark that in the diagram we have considered $\ol{M}_{0,2}$
as a point.

Call $\d^{rt}_{M,P}$ the image of $\vartheta^{rt}_{M,P}$
and $\d^{rt,comb}_{M,P}(\tilde{l})$ the image of $\vartheta^{rt,comb}_{M,P}$.
Notice that, if $M$ does not belong to $\mathfrak{P}_{Q,Q'}$, then
$\d^{rt}_{M,P}$ is contained inside $\pa\Cc^{Q''}_{g,P\cup Q'}$
and $\d^{rt,comb}_{M,P}(\tilde{l})$ is contained inside
$\pa\Cc^{Q'',comb}_{g,P\cup Q'}(\tilde{l})$.

For every partition $M$, define the function $\rho|_M$ as
\[
\begin{matrix}
\rho|_M: &  I_M     & \lra & \N \\
         & \i_{\mu} & \lmt & \rho_{\mu} & :=\sum_{q\in\mu}\rho(q)
\end{matrix}
\]
and set $m_i(M):=m_i^{\rho|_M}$ if $i\neq 0$, and $m_0(M)$
in such a way that $\sum_{i\geq 0} (2i+1)m_i(M)=4g-4+2n$.

If $M\in\mathfrak{P}_{Q,Q'}$, then we define another function
$\tau_M$ as
\[
\begin{matrix}
\tau_M: & Q' & \lra & \N \\
        & q  & \lra & \rho_{\mu_q}
\end{matrix}
\]

\begin{definition}
The combinatorial class
{\it with rational tails} $\Wbar^{rt}_{M,\rho,P}$
on $\Mbarcomb_{g,P\cup Q}$
is the image through $\vartheta^{rt,comb}_{M,P}$
of the combinatorial class
$\Wbar_{m_*(M),\rho|_M,P}$ (which lives
on $\Mbarcomb_{g,P\cup I_M}$). Clearly
it is contained inside $\d^{rt,comb}_{M,P}$.
\end{definition}
\begin{lemma}\label{lemma:forget}
Let the notation be as before and
let $M$ be a partition in $\mathfrak{P}_{Q,Q'}$.
Consider the forgetful map
\[
\pi_{Q''}^{comb}:
\Mbarcomb_{g,P\cup Q}(\R_+^P\times\{0\}^Q)\setminus D^{comb}_{Q'',P} \lra
\Mbarcomb_{g,P\cup Q'}(\R_+^P\times\{0\}^{Q'})
\]
Then $\pi_{Q''}^{comb}(\Wbar^{rt}_{M,\rho,P})=\Wbar_{m_*(M),\tau_M,P}$
and the restriction of $\pi_{Q''}^{comb}$ to
$\Wbar^{rt}_{M,\rho,P}$ has degree
\[
\frac{(m_*(M)-m_*^{\tau_M})!}{(m_0(M)-m_0^{\rho|_M})!}
\]
onto its image.
\end{lemma}
\begin{proof}
Consider the following commutative diagram,
where restrictions are understood
and $\pi^{comb}$ forgets those markings $\i_\mu$
such that $\mu\cap Q'=\emptyset$.
\[
\xymatrix@C=2cm{
& \Wbar^{rt}_{M,\rho,P}
\ar[d]^{\pi_{Q''}^{comb}} \\
\Wbar_{m_*(M),\rho|_M,P}
\ar[r]_{\quad\pi^{comb}}
\ar[ru]^{\vartheta^{rt,comb}_{M,P}}
& \Wbar_{m_*(M),\tau_M,P}
}
\]
As the restriction of
$\vartheta^{rt,comb}_{M,P}$
has degree $1$, we only need to compute
the degree of the restriction
of $\pi^{comb}$.

Consider vertices of valency $2i+3$ in
the general simplex of $\Wbar_{m_*(M),\tau_M,P}$:
among them, $m_i^{\tau_M}$ are marked by some
elements of $Q'$ and
$(m_i(M)-m_i^{\tau_M})$ are not.
On the other side, the general simplex of $\Wbar_{m_*(M),\rho|_M,P}$
has $m_i^{\rho|_M}$ vertices marked by some elements of $I_M$.
Thus, the cardinality of the fiber of $\pi^{comb}$ over
a general simplex of $\Wbar_{m_*(M),\tau_M,P}$ is given by the
product for every $i\geq 0$ of
number of ordered $(m_i^{\rho|_M}-m_i^{\tau_M})$-uples
in a set of $(m_i(M)-m_i^{\tau_M})$ elements.
The claim follows because $m_i^\rho=m_i$ for all $i>0$.
\end{proof}

Now we can state the main result.
\begin{mytheorem}
Let $P=\{p_1,\dots,p_n\}$ and $Q=\{q_1,\dots,q_{s+u}\}$ and
let be given $\rho:Q\lra\N$ and $m_*=(0,m_0,m_1,\dots)$
such that $m_i=m^{\rho}_i$ for $i>0$ and $m_0 \geq m_0^{\rho}$.
Then the following equation holds in
$H_{6g-6+2n-2r}(\Mbar_{g,P\cup Q},\pa\Cc_P^Q)$
\begin{multline*}
\prod_{q\in Q}
2^{\rho(q)+1}
(2\rho(q)+1)!!
(\psi_q^{\rho(q)+1})^* =
\Wbar_{m_*,\rho,P}+ \!\!\!\!
\sum_{\substack{M \in \mathfrak{P}_Q \\ \text{$M$ not discrete}}}
\!\!\!\!
c_{_M} \Wbar^{rt}_{M,\rho,P}
\end{multline*}
where $r=\sum_{i\geq 1}i\, m_i$ and the coefficient $c_{_M}$ is
defined as
\[
\begin{array}{ccc}
\displaystyle{
c_{_M}:=\prod_{\mu\in M} c_{\mu}
} & \text{and} &
\displaystyle{
c_{\mu}=\frac{(2\rho_{\mu}+2|\mu|-1)!!}{(2\rho_{\mu} +1)!!} \, .
}
\end{array}
\]
\end{mytheorem}
\begin{mycorollary}\label{co:second2}
With the same notation as in the theorem,
the following relation holds
\begin{multline*}
\left( \prod_{q\in Q}
2^{\rho(q)+1} (2\rho(q)+1)!! \right) 
\prod_{q\in Q'}(\psi_q^{\rho(q)+1})^*
\sum_{\s\in\mathfrak{S}_{Q''}} \k^*_{\rho_\s} =\\
=\sum_{M \in \mathfrak{P}_{Q,Q'}}
\frac{(m_*(M)-m_*^{\tau_M})!}{(m_0(M)-m_0^{\rho|_M})!} c_{_M}
\Wbar_{m_*(M),\tau_M,P}
\end{multline*}
in $H_{6g-6+2n-2r}(\Mbar_{g,P\cup Q'},\pa\M_{g,P\cup Q'})$.
\end{mycorollary}
We remark that Theorem \ref{th:first} and Corollary \ref{co:second2} give an
inductive method to express all $\Wbar_{m_*,\rho,P}$ in terms of the
tautological classes and vice versa.
In fact it is sufficient to isolate the term on the right hand side
which corresponds to the discrete partition to obtain
the recursion or to isolate the term on the left hand side that
corresponds to $\s=e$.
\begin{proof}[Proof of Corollary \ref{co:second2}]
The discussion developed in Section \ref{section:pi2}
and an inspection of the proof of the previous theorem
show that the equality of Theorem \ref{th:second}
lifts to
$H_{6g-6+2n-2r}(\Mbar_{g,P\cup Q}\setminus D_{Q'',P},
\pa\Cc^{Q''}_{g,P\cup Q'}\setminus D_{Q'',P})$,
if we choose an adequate representative of
the dual of
the monomial in the $\psi$'s in the left hand side.

Then we can push the relation forward through
the forgetful map (look at the big diagram in Section \ref{section:pi2})
and we get our result by using Faber's formula
(see Section \ref{section:tautological})
and Lemma \ref{lemma:forget}.
\end{proof}
\begin{strategy}
The proof of Theorem \ref{th:second}
basically follows that of Theorem \ref{th:first}.
To establish the asserted equality, we evaluate
the integral
\[
\int_{\Mbarcomb_{g,P\cup Q}(\tilde{l}')}
(\H_0^{Q})^*(\eta)\wedge
(\H_0^{q_{s+u}}\cdots\H_0^{q_2})^*(\ol{\omega}_{q_1}^{\rho(q_1)+1})
\wedge\cdots\wedge
 \ol{\omega}_{q_{s+u}}^{\rho(q_{s+u})+1}
\]
where $[\eta]\in H^{6g-6+2n-2r}(\Mbarcomb_{g,P\cup Q}(\tilde{l})
\setminus D^{comb}_{Q,P},
\pa\Cc^{Q,comb}_{g,P}(\tilde{l})\setminus D^{comb}_{Q,P})$
with $\tilde{l}\in C_{P,0}$,
and $\tilde{l}'$ belongs to $C^+_{P,s+u}$.

As in the proof of Theorem \ref{th:first}, the integral
is restricted to a certain subcomplex and decomposed into
a sum of integrals on fibrations (which are
restrictions of $\H^Q_0$), so that we
can integrate the $\ol{\omega}$'s over the fibers of $\H^Q_0$.
Each fiber consists of a cellular complex: the volume
of a single cell (which is a product of simplices)
is responsible of the factor in
the left hand side, while $c_{_M}$ counts the number
of maximal cells in the fiber.
\end{strategy}
\begin{proof}[Proof of Theorem \ref{th:second}]
Pick a closed PL differential form $\eta$
of degree $6g-6+2n-2r$ on
$\Mbarcomb_{g,P\cup Q}(C_{P,0})\setminus D^{comb}_{Q,P}$,
that vanishes on $\pa\Cc_{g,P}^{Q,comb}(C_{P,0})\setminus D^{comb}_{Q,P}$.

Consider then the forms $\beta_i(\eta)$
on $\Mbarcomb_{g,P\cup Q}(C_{P,i}^+)\setminus D^{comb}_{Q,P}$
defined as
\[
(\H_0^{Q_i})^*(\eta)\wedge
(\H_0^{q_i}\cdots\H_0^{q_2})^*(\ol{\omega}_{q_1}^{\rho(q_1)+1})
\wedge\cdots\wedge
 \ol{\omega}_{q_i}^{\rho(q_i)+1}
\]
for $i=1,\dots,s+u$
and remember that $\beta_{s+u}(\eta)$
extends by zero outside
$\Mbarcomb_{g,P\cup Q}(C_{P,s+u}^+) \setminus D^{comb}_{Q,P}$.
As already mentioned in Section \ref{section:pi}, the
form $\beta_{s+u}(1)$ has the property that its cohomology class
pulls back via $\tilde{\xi}$ to
\[
\psi_{q_1}^{\rho(q_1)+1}
\psi_{q_2}^{\rho(q_2)+1}
\cdots
\psi_{q_{s+u}}^{\rho(q_{s+u})+1}
\]
on $\Mbar_{g,P\cup Q}\times C_{P,s+u}^+$.

As in the proof of Theorem \ref{th:first},
it is easy to see that the restriction of $\beta_{s+u}(\eta)$
to $\Mbarcomb_{g,P\cup Q}(\tilde{l})$
for some $\tilde{l}\in C_{P,s+u}^+$ 
has support
contained inside the locus
$\ol{Y}_{2\rho(q_1)+3,\dots,2\rho(q_{s+u})+3}(\tilde{l})$
by reasons of degree.
Now we want to analyze its image through $\H^Q_0$
which consists of several components.
\begin{definition}
Given a ribbon graph, we say that two holes
are {\it adjacent} if they have at least one vertex in common.
Moreover, we say that a subset $\mu$ of markings forms
a {\it cluster} if
\begin{enumerate}
\item[-]
every vertex of $x(\mu)$ contains an edge that belongs to a hole
in $x(\mu)$
\item[-]
every two distinct holes in $x(\mu)$
are joined by a chain of pairwise
adjacent holes belonging to $x(\mu)$.
\end{enumerate}
Two clusters $\mu$ and $\mu'$ are {\it disjoint} if $\mu\cup\mu'$ is
not a cluster (in particular $\mu$ and $\mu'$ are disjoint as sets).
\end{definition}
We associate to every partition $M=\{\mu_1,\dots,\mu_k\}$ in
$\mathfrak{P}_Q$ the closure $\ol{Y}_M(C_{P,s+u}^+)$
of the locus of top-dimensional simplices
of $\ol{Y}_{2\rho(q_1)+3,\dots,2\rho(q_{s+u})+3}(C_{P,s+u}^+)$
such that $\mu_1,\dots,\mu_k$ form disjoint
clusters.
It is obvious that $\cup_M \ol{Y}_M(C_{P,s+u}^+)$ is a dissection of
$\ol{Y}_{2\rho(q_1)+3,\dots,2\rho(q_{s+u})+3}(C_{P,s+u}^+)$.
Really, these subcomplexes overlap on simplices
of nonmaximal dimension, but this fact is not important for what follows.
Strictly speaking, we need a refinement of this dissection:
for every tripartition $M^{\bullet}=(M^{disk},M^{cyl},M^{surf})$
of $M$, we denote
by $\ol{Y}_{M^{\bullet}}(C_{P,s+u}^+)$ the closure of the locus
in $\ol{Y}_M(C_{P,s+u}^+)$ where the subsurface
$\cup_{q\in\mu} \bar{T}_q$ form a disk
(resp. a cylinder or a surface with negative Euler characteristic)
for every cluster $\mu$ in $M^{disk}$
(resp. in $M^{cyl}$ or in $M^{surf}$).

Then $\H^Q_0(\ol{Y}_{M^{\bullet}}(C_{P,s+u}^+))$
is the smallest subcomplex of $\Mbarcomb_{g,P\cup Q}(C_{P,0})$
containing all the simplices
indexed by ribbon graphs $G$ such that:
\begin{enumerate}
\item
every $\mu \in M^{disk}$ marks a nonpositive
sphere that intersects only one positive component
in a vertex of $G$ of valency $2\rho_{\mu}+3$
(if $\mu=\{q_j\}$ we should simply say: $q_j$ marks
a vertex of $G$ of valency $2\rho_{\mu}+3$),
while all the other vertices
in the smooth locus of $G$ are trivalent
\item
every $\mu \in M^{cyl}$ marks a nonpositive sphere
that intersects the positive components in two
vertices of $G$ of valencies $v_1$ and $v_2$,
with $v_1+v_2=2\rho_{\mu}$
\item
every $\mu \in M^{surf}$ marks a nonpositive component
of some genus $h$ which intersects the positive components
in $\nu$ vertices of $G$ of valencies
$v_1,\dots,v_{\nu}$, provided
$6h-6+\sum_{j=1}^{\nu}(v_j+3)=2\rho_{\mu}$.
\end{enumerate}

Notice that $\H^Q_0(\ol{Y}_{M^{\bullet}}(C_{P,s+u}^+))$ is contained
inside $\pa\Cc^{Q,comb}_{g,P}(C_{P,0})$ if
$M^{cyl}$ or $M^{surf}$ are nonempty. Hence, we can
restrict to $M^{cyl}=M^{surf}=\emptyset$,
in which case
\begin{align*}
\H^Q_0(\ol{Y}_{(M,\emptyset,\emptyset)}(C_{P,s+u}^+))=
\begin{cases}
\Wbar_{m_*,\rho,P} & \text{if $M$ is discrete} \\
\Wbar^{rt}_{M,\rho,P} & \text{if $M$ is not discrete.}
\end{cases}
\end{align*}

The configuration of the $Q$-marked holes in the ribbon
graphs $G$ that correspond to simplices of top dimension in
$\ol{Y}_{(M,\emptyset,\emptyset)}$
is not so difficult to describe.
We can restrict our analysis to a single cluster $\mu\in M$,
because in $\ol{Y}_{(M,\emptyset,\emptyset)}$ they are
disjoint.

We already know that the cluster $\mu$ has circular shape,
that is $\cup_{q\in\mu} \bar{T}_q\subset S(G)$
is a disk. The fact that $\rho$ takes nonnegative
values leads to a stronger conclusion.
\begin{lemma}\label{lemma:disk}
In the previous hypotheses, every $\bar{T}_q$ is a disk.
\end{lemma}
\begin{proof}[Proof of Lemma \ref{lemma:disk}]
By contradiction, suppose $h$ to be the minimum integer such
that $q_h\in\mu$ and $\bar{T}_{q_h}$ is not a disk.
Then it is a surface of
genus zero with at least two boundary components. Moreover,
if we remove $\bar{T}_{q_h}$, we disconnect the surface $S(G)$.
Call $S_1,\dots,S_f$ the connected components
of $S(G)\setminus \bar{T}_{q_h}$
(which are necessarily disks)
that do not contain $P$-markings,
Call, for instance,
$q_{i_1},\dots,q_{i_j}$ the marked points of $S_1$, with $i_1<\dots<i_j$.

One can see that $i_2<h$ is impossible, because
$\H^{q_h}_0\cdots\H^{q_{s+u}}_0$ sends $S_1$ to a genus zero component
with a node and at least two marked points. A simple
Euler characteristic computation shows that the holes in $S_1$
cannot have the right number of edges. In algebro-geometric terms,
the equivalent assertion is that on $\Mbar_{0,m+1}$ the class
$\psi_1^{t_1}\cdots\psi_m^{t_m}$ vanishes
if $t_1,\dots,r_m\geq 1$, because $\mathrm{dim}_{\C}\Mbar_{0,m+1}=m-2$.

One can also exclude the case $i_1<h<i_2$. In fact, the
hole $\bar{T}_{q_{i_1}}$ is a disk (because $i_1<h$)
and $G$ has no univalent vertices,
hence $\bar{T}_{q_{i_1}}$ does not contain internal edges.
As a consequence,
the hole $q_{i_1}$ has all its edges in common with
$q_h,q_{i_2},\dots,q_{i_j}$, which is impossible.

A similar conclusion holds if we replace $S_1$ with any
other component $S_2,\dots,S_f$. Thus, the hole $q_h$ is collapsed
after all the components $S_1,\dots,S_f$.
Hence,
$\H^{q_{h+1}}_0\dots\H^{q_{s+u}}_0$ shrinks $S_1,\dots,S_f$
to nonpositive spheres
that touch the rest of the surface in vertices of valency
at least $3$, because $\rho$ takes nonnegative values.
Moreover, these vertices sit in the internal part of
the hole $q_h$. On the other side, after the previous shrinking,
the hole $\bar{T}_{q_h}$ has necessarily become a disk, and
a disk that has internal edges must contain also
an internal univalent vertex.
Thus, we have found a contradiction. Hence, every $\bar{T}_{q_h}$
is a disk.
\end{proof}
As it is evident, we have used in an essential way the hypothesis
that $\rho$ takes nonnegative values.
\begin{remark}
A deeper inspection of the previous proof reveals more.
Not only the map $\H^Q_0$ shrinks the cluster associated to $\mu$
to a vertex (if we do not care about nonpositive components),
but also every partial shrinking $\H^{q_t}_0\cdots\H^{q_{s+u}}_0$
really sends the cluster associated to $\mu$ to a cluster
with circular shape.
\end{remark}

To integrate our differential form $\beta_{s+u}(\eta)$,
we proceed as in the
proof of Theorem \ref{th:first}, shrinking one hole at
a time. We choose positive lengths $L''<L'$, $\e''_i<\e'_i$
for all $i=1,\dots,s+u-1$ and $\e'_{s+u}$ in such a way that
\[
\sum_{j=i}^{s+u}
\e'_j \ll
\mathrm{min}\{\e''_{i-1},\dots,\e''_1,L'',
\e'_{i-1}-\e''_{i-1},\dots,\e'_1-\e''_1,L'-L''\}
\]
and we evaluate
{\footnotesize{
\begin{align*}
& \int_{\Mbar_{g,P\cup Q}}\psi_{q_1}^{\rho(q_1)+1}\smallsmile
\dots\smallsmile\psi_{q_{s+u}}^{\rho(q_{s+u})+1}\smallsmile
(\H^{Q}_0\tilde{\xi}_{\tilde{l}})^*[\eta] = \\
= & \lim_{\e'_{s+u} \rar 0}
\sum_{M\in\mathfrak{P}_Q}
\int_{\ol{Y}_{(M,\emptyset,\emptyset)}([L'',L']^n
\times\dots\times\{\e'_{s+u}\})}
\frac{
\beta_{s+u}(\eta) \wedge
\bar{\l}^*(d\tilde{l}_{p_1}\wedge\dots\wedge d\tilde{l}_{q_{s+u-1}})
}
{ (L'-L'')^n (\e'_1-\e''_1)\cdots(\e'_{s+u-1}-\e''_{s+u-1}) } = \\
= & \sum_{M\in\mathfrak{P}_Q}
\int_{
\H_0^{q_{s+u}}(\ol{Y}_{(M,\emptyset,\emptyset)})([L'',L']^n
\times\dots\times[\e''_{s+u-1},\e'_{s+u-1}]\times\{0\})
}
\beta_{s+u-1}(\eta) \wedge \\
& \quad \wedge \frac{
\bar{\l}^*(d\tilde{l}_{p_1}\wedge\dots\wedge d\tilde{l}_{q_{s+u-1}})
} { (L'-L'')^n (\e'_1-\e''_1)\cdots(\e'_{s+u-1}-\e''_{s+u-1}) }
\int_{F_M^{s+u}(\e'_{s+u})}
\ol{\omega}_{s+u}^{\rho(q_{s+u})+1}
\end{align*} }}
where $F_M^{s+u}(\e'_{s+u})$ is the intersection of
$\ol{Y}_{(M,\emptyset,\emptyset)}([L'',L']^n
\times\dots\times\{\e'_{s+u}\})$ with the generic fiber
of $\H_0^{q_{s+u}}$ over
$\H_0^{q_{s+u}}(\ol{Y}_{(M,\emptyset,\emptyset)})(L'',\dots,\e''_{s+u-1},0)$.
Moreover,
\[
\int_{F_M^{s+u}(\e'_{s+u})}
\ol{\omega}_{s+u}^{\rho(q_{s+u})+1}=
\frac{(\rho(q_{s+u})+1)!}{(2\rho(q_{s+u})+2)!} N_M^{s+u}
\]
where $N_M^{s+u}$ is the number of simplices of top dimension
in $F_M^{s+u}(\e'_{s+u})$.

Then we let $\e''_{s+u-1}$ and $\e'_{s+u-1}$ go to zero,
keeping their difference $\e'_{s+u-1}-\e''_{s+u-1}$ constant,
so that the integral above is equal to
{\footnotesize{
\begin{align*}
& \sum_{M\in\mathfrak{P}_Q}
\int_{
\H_0^{\{q_{s+u-1},q_{s+u}\}}(\ol{Y}_{(M,\emptyset,\emptyset)})([L'',L']^n
\times\dots\times[\e''_{s+u-2},\e'_{s+u-2}]\times\{0\})
}
\beta_{s+u-2}(\eta) \wedge \\
& \quad \wedge \frac{
\bar{\l}^*(d\tilde{l}_{p_1}\wedge\dots\wedge d\tilde{l}_{q_{s+u-2}})
} { (L'-L'')^n (\e'_1-\e''_1)\cdots(\e'_{s+u-2}-\e''_{s+u-2}) }
\frac{(\rho(q_{s+u})+1)! (\rho(q_{s+u-1})+1)!}
{(2\rho(q_{s+u})+2)! (2\rho(q_{s+u-1})+2)!} N_M^{s+u-1}
\end{align*} }}
where $N_M^{s+u-1}$ is the number of simplices of top dimension
contained in the intersection of
$\ol{Y}_{(M,\emptyset,\emptyset)}([L'',L']^n
\times\dots\times\{\e'_{s+u}\})$ with
the generic fiber of $\H_0^{\{q_{s+u-1},q_{s+u}\}}$ over
$\H_0^{\{q_{s+u-1},q_{s+u}\}}(\ol{Y}_{(M,\emptyset,\emptyset)})
(L'',\dots,\e''_{s+u-2},0,0)$.

At the end, we obtain that the first integral is equal to
{\footnotesize{
\begin{align*}
& \sum_{M\in\mathfrak{P}_Q}
\int_{
\H_0^{Q}(\ol{Y}_{(M,\emptyset,\emptyset)})([L'',L']^n
\times\{0\})
}
\eta \wedge
\frac{
\bar{\l}^*(d\tilde{l}_{p_1}\wedge\dots\wedge d\tilde{l}_{p_n})
} { (L'-L'')^n } \\
& \hspace{5cm}
\frac{(\rho(q_{s+u})+1)!\dots (\rho(q_{1})+1)!}
{(2\rho(q_{s+u})+2)! \dots (2\rho(q_{1})+2)!} N_M^{1} =\\
= &
\Big(
\int_{\Wbar_{m_*,\rho,P}(l,0) }
\eta +
\sum_{\substack{M\in\mathfrak{P}_Q \\ M\,\, \text{not discrete}}}
c_{_M}
\int_{\Wbar^{rt}_{M,\rho,P}(l,0) }
\eta \Big) \\
& \hspace{5cm}
\frac{(\rho(q_{s+u})+1)!\dots (\rho(q_{1})+1)!}
{(2\rho(q_{s+u})+2)! \dots (2\rho(q_{1})+2)!}
\end{align*} }}
where $l$ belongs to $[L'',L']^n$ and
$c_{_M}:=N_M^1$ is the number of top-dimensional simplices
contained in the intersection $F_M^1$ of
$\ol{Y}_{(M,\emptyset,\emptyset)}([L'',L']^n
\times\dots\times\{\e'_{s+u}\})$ with
the generic fiber of $\H_0^Q$ over
$\H_0^{Q}(\ol{Y}_{(M,\emptyset,\emptyset)})(L'',0)$.
\begin{figure}[h]
\resizebox{7cm}{!}{\input{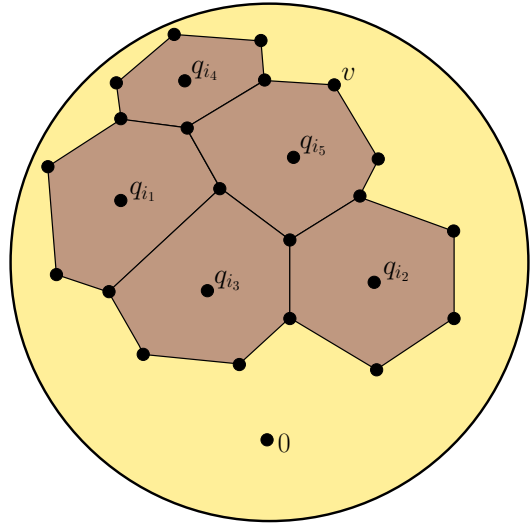}}
\caption{An example of admissible cluster}
\label{fig:admissible}
\end{figure}

Now we are left to determine $c_{_M}$.
In order to reconstruct a ribbon graph $G$ parametrizing
a maximal simplex in $F^1_M$
from its image through $\H_0^Q$, it is sufficient
to know the ``configuration'' of the clusters
associated to $\mu_1,\dots,\mu_k$.
The key point is that this enumeration is ``local'': we mean that
we only use that each cluster is a disk and it is
made of holes that keep their shape circular after
each shrinking, and that these holes
have the ``right'' number of edges.
Briefly, we do not need to know about the surface
apart from the clusters.
Hence, if we denote by
$c_{\mu_i}$ the number of possible configurations for
the cluster $\mu_i$,
then we can conclude that $F^1_M$
contains $c_{_M}=c_{\mu_1}\cdots c_{\mu_k}$ maximal simplices.

Consequently, we need to determine the number $c_{\mu}$
of all possible ``isomorphism types'' of these clusters associated to $\mu$.
More formally, by ``isomorphism classes'' of {\it admissible clusters}
we mean isomorphism classes of ribbon graphs $G'$ such that:
\begin{itemize}
\item[-]
$G'$ is a connected ordinary ribbon graph marked by the set $\mu\cup\{0,v\}$
\item[-]
$S(G')$ is a sphere and $\mu$ forms a cluster
\item[-]
$\bar{T}_0$ is a disk and for every $q_j\in\mu$
the subsurface $\bar{T}_{q_j}$ is a disk
\item[-]
for every $j=2,\dots,h$, the shrinking
of the holes $q_{i_h},q_{i_{h-1}},\dots,q_{i_j}$
produces a ribbon graph with only one positive component
\item[-]
the vertices of $G'$ have valency two or three;
the bivalent ones always border the hole $0$
and one of them is marked by $v$
\item[-]
if $\mu=\{q_{i_1},\dots,q_{i_h}\}$ with $i_1<i_2<\dots<i_h$,
then the hole $q_{i_h}$ has $2\rho(q_{i_h})+3$ sides,
and for all $j=1,\dots,h-1$ the hole $q_{i_j}$ has $2\rho(q_{i_j})+3$ sides
beside those which border the holes $q_{i_h},\dots,q_{i_{j+1}}$.
\end{itemize}
(Although we are only interested in the case $\rho(q_{i_j})\geq 0$,
we remark that the above definition and part of the subsequent
calculation
makes sense also for $\rho(q_{i_1})=-1$.)

Given $G'_{\mu_1},\dots,G'_{\mu_k}$, 
we can reconstruct $G$ from $\H^Q_0(G)$
in the following way.
For every $\mu_i\in M$, we discard the nonpositive
component associated to $\mu_i$ and we expand the
corresponding $(2\rho_{\mu_i}+3)$-valent vertex in order to obtain
a polygon with $2\rho_{\mu_i}+3$ edges bordering a new hole.
Then we glue $G'_{\mu_i}$ to $G$ in such
a way that the new hole of $G$ corresponds
to the $0$-th hole of $G'_{\mu_i}$ and the
bivalent vertices of $G'_{\mu_i}$ are glued
to the vertices of the polygon.
The different $v$-markings on $G'_{\mu_i}$ count
in how many ways we can perform this gluing.

So we are left to prove the following lemma.
\end{proof}
\begin{lemma}\label{lemma:enum}
The number of isomorphism classes of admissible
clusters associated to $\mu$ is
\[
c_{\mu}=
\frac{(2\rho_{\mu}+2|\mu|-1)!!}{(2\rho_{\mu}+1)!!}.
\]
where we have conventionally set $(-1)!!=1$.
\end{lemma}
\begin{proof}[Proof of Lemma \ref{lemma:enum}]
Remark that the calculation has a clear geometrical
meaning even if we allow $\rho(q_{i_1})$ to assume the value $-1$.

Clearly, if $h=1$ then $c_{\mu}=1$.
If $h>1$, then the cluster has no rotational symmetries and
so the possible $v$-markings are exactly
$2\rho_{\mu}+3$.

In particular, if $h=2$
then the cluster consists of the holes $q_{i_1}$
with $2\rho(q_{i_1})+4$ sides
and $q_{i_2}$ with $2\rho(q_{i_2})+3$ sides
that have exactly one edge in common.
Thus, in this case $c_{\mu}=2\rho_{\mu}+3$.

Now we deal with the case $h>2$, so we suppose that the formula
holds for all $|\mu|<h$, with $\rho(q_{i_1})\geq -1$ and
$\rho(q_{i_2}),\dots,\rho(q_{i_{|\mu|}})\geq 0$. We want
to prove that the formula holds for $|\mu|=h$, with
$\rho(q_{i_1})\geq -1$ and $\rho(q_{i_2}),\dots,\rho(q_{i_h})\geq 0$.

If $\rho(q_{i_1})=-1$, then we look at the situation just
before shrinking $q_{i_1}$. We have a loop surrounding $q_{i_1}$ and
its vertex has valency $2\rho_{\mu}+5$. So this vertex is obtained
collapsing the subcluster $\mu'=\mu\setminus\{q_{i_1}\}$.
By induction hypothesis, $c_{\mu}$ is
$(2\rho_{\mu}+3)\rho_{\mu'}=(2\rho_{\mu}+3)(2(\rho_{\mu}+1)+2(|\mu|-1)-1)
\cdots (2(\rho_{\mu}+1)+3)=
(2\rho_{\mu}+3)(2\rho_{\mu}+2|\mu|-1)\cdots(2\rho_{\mu}+5)$.
 
If $\rho(q_{i_1})\geq 0$, then
we look at the situation before collapsing $q_{i_2}$ and $q_{i_1}$.
There are two possibilities: the holes may touch each other in one edge
(case (a)) or in one vertex (case (b)).
We want to show that in both cases the number
of configurations (which we denote respectively
by $c^a_{\mu}$ and $c^b_{\mu}$)
depends only on $\rho(q_{i_1})+\rho(q_{i_2})$ and not on $\rho(q_{i_1})$
and $\rho(q_{i_2})$ separately. Hence $c_{\mu}$ depends
only on $\rho(q_{i_1})+\rho(q_{i_2})$ too, because
$c_{\mu}=c^a_{\mu}+c^b_{\mu}$.
Henceforth, we can use the computation made for $\rho(q_{i_1})=-1$
and we are done.
\begin{figure}[h]
\resizebox{8cm}{!}{\input{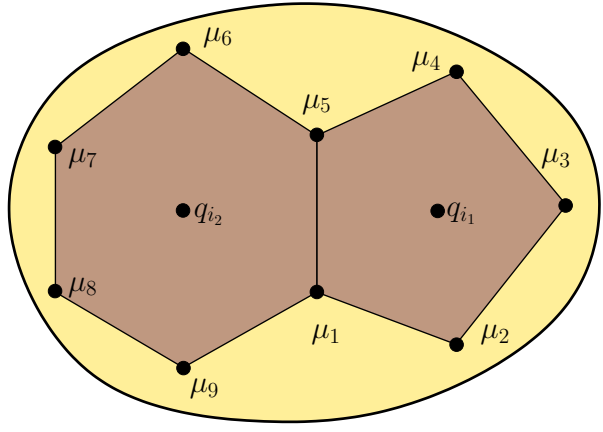}}
\caption{An example of case (a) with $\rho(q_{i_1})=\rho(q_{i_2})=1$}
\label{fig:adm1}
\end{figure}

In case (a),
the two holes $q_{i_1}$ and $q_{i_2}$ touch in one edge.
Hence, the holes $q_{i_3},\dots,q_{i_h}$ are distributed in
$t=(2\rho(q_{i_1})+3)+(2\rho(q_{i_2})+4)-2=2(\rho(q_{i_1})+\rho(q_{i_2}))+5$
subclusters $\mu_1,\dots,\mu_t$.
Then we obtain
\[
c^{a}_{\mu}=(2\rho_{\mu}+3)\sum_{j\in J}
\prod_{k=1}^t c_{j^{-1}(k)}
\]
where $J=\{j:\mu\setminus\{q_{i_1},q_{i_2}\}\rar \{1,\dots,t\} \, \}$
and conventionally $c_\emptyset=1$,
which depends only on $\rho(q_{i_1})+\rho(q_{i_2})$
and not on $\rho(q_{i_1})$ and $\rho(q_{i_2})$ separately.
\begin{figure}[h]
\resizebox{8cm}{!}{\input{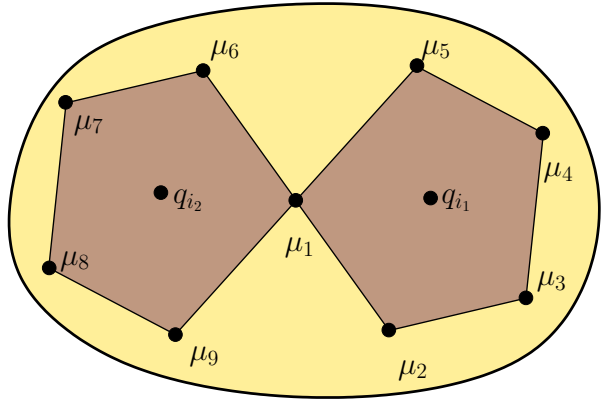}}
\caption{An example of case (b) with $\rho(q_{i_1})=\rho(q_{i_2})=1$}
\label{fig:adm2}
\end{figure}
In case (b),
the two holes $q_{i_1}$ and $q_{i_2}$ touch in a vertex.
Hence, we have
$t=(2\rho(q_{i_1})+3)+(2\rho(q_{i_2})+3)-1=
2(\rho(q_{i_1})+\rho(q_{i_2}))+5$ subclusters; but
the cluster $\mu_1$,
which corresponds to the common vertex, must
be at least $4$-valent.
Moreover, the two holes can touch
the cluster $\mu_1$
in $2\rho_{\mu_1}$ ways,
hence we obtain
\[
c^{b}_{\mu}=(2\rho_{\mu}+3)(2\rho_{\mu_1})
\sum_{j\in J}  \prod_{k=1}^t c_{j^{-1}(k)}
\]
where  $J=\{j:\mu\setminus\{q_{i_1},q_{i_2}\}\rar \{1,\dots,t\}| 
\rho_{j^{-1}(1)}\geq 1 \}$,
which depends only on $\rho(q_{i_1})+\rho(q_{i_2})$
and not on $\rho(q_{i_1})$ and $\rho(q_{i_2})$ separately.
\end{proof}
%
%
\begin{example}
In the case $\mu=\{q_{i_1},q_{i_2},q_{i_3}\}$, the induction
works as follows.
\begin{figure}[h]
\resizebox{10cm}{!}{\input{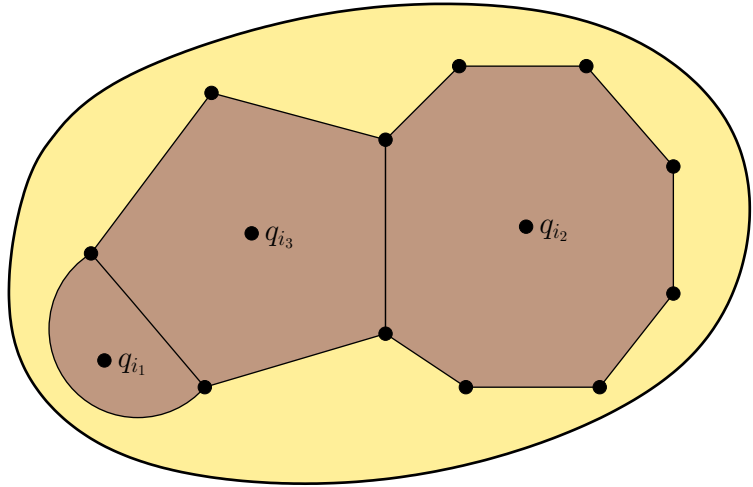}}
\caption{Example with $h=3$, $\rho(q_{i_1})=-1$,
$\rho(q_{i_2})=2$ and $\rho(q_{i_3})=1$}
\label{fig:admex}
\end{figure}
If $\rho(q_{i_1})=-1$, we have that $q_{i_2}$ and $q_{i_3}$
touch in an edge. Then the hole $q_{i_1}$ can be ``attached''
in $2\rho(q_{i_2})+2\rho(q_{i_3})+3$ ways on the rest
of the cluster (see Fig. \ref{fig:admex}).
At the end, there are $2\rho(q_{i_2})+2\rho(q_{i_3})+1$
bivalent vertices that can be marked by $v$. Hence, there are
$(2\rho_{\mu}+3)(2\rho_{\mu}+5)$ total possibilities, and
so in this case $c_{\mu}=(2\rho_{\mu}+5)(2\rho_{\mu}+3)$.

If $\rho(q_{i_1})\geq 0$, then we look at the cluster after
the collapsing of $q_{i_3}$.

In case (a) the two holes $q_{i_1}$ and $q_{i_2}$ touch in an edge
and there are $2\rho(q_{i_1})+2\rho(q_{i_2})+5$ vertices where
the $q_{i_3}$-marking could be. So the number of total possibilities in
case (a) is $c^a_{\mu}=(2\rho_{\mu}+3)(2\rho(q_{i_1})+2\rho(q_{i_2})+5)$,
which correctly does not depend on $\rho(q_{i_1})$ and $\rho(q_{i_2})$
separately but only on $\rho(q_{i_1})+\rho(q_{i_2})$.

In case (b) the two holes $q_{i_1}$ and $q_{i_2}$ touch in a vertex
of valency at least $4$,
which is necessarily marked by $q_{i_3}$ (and so this case may
occur only if $\rho(q_{i_3})\geq 1$).
It is easy to see that, before the shrinking,
there were exactly $2\rho(q_{i_3})$ distinct
configurations in which $q_{i_3}$ had one edge
in common with $q_{i_1}$ and one edge in common with $q_{i_2}$,
but $q_{i_1}$ and $q_{i_2}$ did not have edges in common
(here we use that at the beginning vertices have valencies
at most $3$). Thus, in case (b) the number of possibilities is
$c^b_{\mu}=(2\rho_{\mu}+3)(2\rho(q_{i_3}))$, which correctly vanishes
if $\rho(q_{i_3})=0$ and which depends on $\rho(q_{i_1})+\rho(q_{i_2})$
and not on $\rho(q_{i_1})$ and $\rho(q_{i_2})$ separately.
Finally, there are $c^a_{\mu}+c^b_{\mu}=(2\rho_{\mu}+3)(2\rho_{\mu}+5)$
total configurations.
\end{example}
As an application, we compute the case $\Wbar_{2a+3,2b+3}$
of graphs with two vertices of valencies $2a+3$ and $2b+3$.
\begin{mycorollary}
For every $g\geq 0$ and $P\neq\emptyset$ such that $2g-2+|P|>0$
and for every $a,b\geq 1$, the
following identity
\begin{multline*}
2^{\d_{a,b}}\Wbar_{2a+3,2b+3}
=2^{a+b+2}(2a+1)!!(2b+1)!!(\k^*_a\k^*_b+\k^*_{a+b}) \\
-2^{a+b+1}(2a+2b+3)!!\k^*_{a+b}
\end{multline*}
holds in $H_{6g-6+2n-2a-2b}(\Mbar_{g,P},\pa\M_{g,P})$.
\end{mycorollary}
\begin{proof}
It is a simple application of Theorem \ref{th:second}.
In fact, from the recursive formula one obtains
\begin{align*}
& 2^{a+b+2}(2a+1)!!(2b+1)!!(\k^*_a\k^*_b+\k^*_{a+b})= \\
& =2^{\d_{a,b}}\Wbar_{2a+3,2b+3}+(2a+2b+3)\Wbar_{2a+2b+3}= \\
& =2^{\d_{a,b}}\Wbar_{2a+3,2b+3}+
(2a+2b+3)\left[2^{a+b+1}(2a+2b+1)!!\, \k^*_{a+b}\right].
\end{align*}
\end{proof}
\end{section}
%
%
\appendix
\section*[]{Appendix. On Looijenga's modification
 of the arc complex}
It is evident that a lot of technical problems
come from the fact that we have not found a
geometrical way to lift the combinatorial classes
to the Deligne-Mumford compactification of the moduli
space. In other words, we have not a triangulation
of $\Mbar_{g,P}$ that supports combinatorial classes.

In \cite{looijenga:cellular}, Looijenga defined a
modification $\Ahat(S,P)$ of the arc complex that maps to $\Mbar_{g,P}$.
If $\Ahat(S,P)$ supported combinatorial cycles, then
we could push their class forward to $\Mbar_{g,P}$.
The critical point is to prove that the combinatorial
subcomplexes are really cycles. We are unable to do that
in general, but we can handle the case of one nontrivalent vertex.

Here we sketch the construction of $\Ahat(S,P)$, but
we refer to Looijenga's paper for a more detailed treatment.
\subsection*{The modified arc complex}
The notation being as in Section \ref{sec:combinatorial},
let $Z \subset X_1(G)$ be a nonempty subset of edges of an ordinary
connected ribbon graph $G$.
We can construct two new ribbon graphs.

The {\it subgraph} $G_Z=(X(G_Z),\s'_0,\s'_1)$
has $X(G_Z)$ equal
to the set of orientations of edges in $Z$, its
$\s'_1$ is the natural restriction of $\s_1$, and its $\s'_0$ sends
an oriented edge
to the next one belonging to $X(G_Z)$ with respect to the
cyclic order induced by $\s_0$. If $Z$ does not coincide with $X_1(G)$, then
$G_Z$ has some new {\it exceptional holes} corresponding
to orbits in $X(G_Z)\subset X(G)$ under $\s'_{\infty}$
which are not orbits under the action of $\s_{\infty}$.

Consider now a proper subset $Z$ of $X_1(G)$. Then
the {\it quotient graph} $G/G_Z$ has $X(G/G_Z)$ equal to
$X(G)\setminus X(G_Z)$,
its $\s'_1$ is the restriction and its $\s'_\infty$ sends
an oriented edge to the next one of $X(G/G_Z)$ with respect to
the cyclic order induced by $\s_\infty$. If $Z$ is nonempty,
then $G/G_Z$ has {\it exceptional vertices}
corresponding to orbits in $X(G_Z)\subset X(G)$ under $\s'_0$
that are not orbits under the action of $\s_0$.

\begin{figure}[h]
\resizebox{12cm}{!}{\input{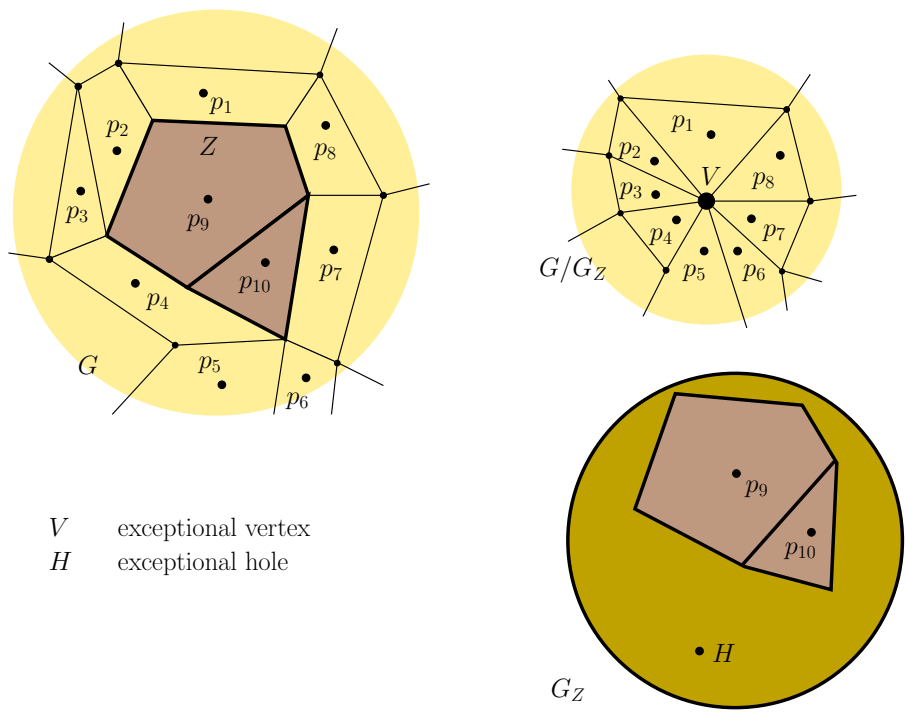}}
\caption{Example of correspondence between
exceptional holes and exceptional vertices}
\label{fig:exceptional}
\end{figure}

Notice that there is a canonical correspondence
between exceptional vertices of $G/G_Z$ and exceptional
holes of $G_Z$ (see Fig.~\ref{fig:exceptional}).
In fact consider an exceptional hole $H$
of $G_Z$. For every (oriented) edge $e\in H$
call $c_e>0$ the minimum integer such that
$\s_1\s_0^{c_e}(e)$ belongs to $H$. Then the
subset $\{\s_0^i(e)|\,e\in H \,\text{and}\,0<i<c_e\}$
is the corresponding exceptional vertex in $G/G_Z$.
Conversely, given an exceptional vertex $V$ of $G/G_Z$
and an $e\in V$ call $c_e>0$ the minimum integer such that
$\s_{\infty}^{c_e}\s_1(e)$ belongs to $V$. Then
$\{\s_{\infty}^i\s_1(e)|\,e\in H \, \text{and}\,0<i<c_e\}$
is the corresponding exceptional hole in $G_Z$.

To introduce the definition of stable $P$-marked ribbon graph
consider how an ordinary
metrized $P$-marked ribbon graph $G$
can degenerate: it happens when the lengths
of a subset $Z$ of edges go to zero.
As we can work componentwise, we suppose $Z$ connected.
Then various cases can occur:
\begin{enumerate}
\item
$Z$ is a tree and contains at most one marked point,
so it is {\it contractible}:
then we can collapse it to a vertex and put, if necessary,
the marking on it, that is we simply obtain $G/G_Z$
\item
$Z$ is homotopic to a circle and no vertex of $Z$ is marked,
so we call it {\it semistable}.
If $Z$ surrounds a single hole, then it shrinks to
a vertex which inherits the marking in $G/G_Z$;
otherwise $G/G_Z$ contains two exceptional vertices
\item
if $Z$ is neither contractible nor semistable, then its
collapsing
gives rise to a new irreducible component. If $Z$ contains
no unmarked tails, then we call $Z$ a {\it stable} subset. Notice that
every $Z$ which is not contractible or semistable
contains a nonempty maximal stable subset $Z^{st}$.
We denote by $G_{\ol{Z^{st}}}$ the reduced ribbon graph
obtained from $G_{Z^{st}}$ by ``deleting'' the bivalent unmarked
vertices and we set $\ol{Z^{st}}=X_1(G_{\ol{Z^{st}}})$.
\end{enumerate}
Now we want to produce a stable version of ribbon graphs
by successive collapsing of semistable or stable subsets of edges.

Given an ordinary connected $P$-marked ribbon graph $(G,x)$,
we call $Z_\bullet=(Z_0,Z_1,\dots,Z_k)$ a
{\it permissible sequence} for $(G,x)$
if $Z_0=X_1(G)$ and $Z_{j+1}\subset \ol{Z_j^{st}}$
is a nonempty subset not containing a whole component of
$\ol{Z_j^{st}}$ for every $j=0,\dots,k-1$.
Given such a $Z_\bullet$ we can produce a (quasi)stable
$P$-marked ribbon graph taking the triple $(G(Z_\bullet),\xb,\i)$
where
\[
G(Z_{\bullet}):=(G_{Z_0}/G_{Z_1}) \sqcup (G_{\ol{Z_1^{st}}}/G_{Z_2})
\sqcup \cdots \sqcup (G_{\ol{Z_{k-1}^{st}}}/G_{Z_k})
\sqcup G_{\ol{Z_k^{st}}},
\]
$\xb:P\hra X_{\infty}(G(Z_{\bullet}))\cup X_0(G(Z_{\bullet}))$
is induced by $x$
and $\i$ is a fixed-point-free involution
that exchanges
an exceptional hole and its corresponding exceptional vertex.
The ``stabilized'' $P$-marked ribbon graph is simply
obtained discarding possible unstable components, namely unmarked
spheres with two exceptional holes, and making $\i$ exchange
the two corresponding exceptional vertices. In any case, $\i$ never
exchanges two holes.

We say that the (stable) components of $G_{\ol{Z_i^{st}}}/G_{Z_{i+1}}$
have {\it order} $i$ and we define
$H_i$ as the set of holes belonging to components of
order $i$ and $V_i$ as the set of marked or exceptional vertices
belonging to components of order $i$. Finally we say that
$\Sigma:=\cup_i (H_i\cup V_i)$ is the set of {\it special points}.
\begin{definition*}
A {\it stable metric} with respect to $Z_\bullet$ is
a sequence of metrics $(a_i)_{i=0}^k$ where
$a_0 \in \Dc_{Z_0\setminus Z_1}$ and
$a_i$ is a metric on $G_{\ol{Z_i^{st}}}/G_{Z_{i+1}}$, which is positive
of total length $1$ on every irreducible component.
\end{definition*}
Thus, given a stable metric for $Z_\bullet$, we can build a
stable marked Riemann surface $S(G,Z_\bullet,a_\bullet)$.
In fact we first consider the disjoint union of
the surfaces $S(G_{\ol{Z_i^{st}}}/G_{Z_{i+1}},a_i)$
for $i=0,\dots,k$ and then we identify some pairs of
points according to $\i$.
Remark that there is an extended circumference function
\[
\wh{\La}: \{ \text{stable metrics on $S(G,Z_{\bullet})$} \}
\lra \prod_{i=0}^k \D_{H_i}
\]
that restricts to a map
$\hat{\l}:=\wh{\La}_0:
\{\text{stable metrics on $S(G,Z_{\bullet})$}\}\lra \D_P$.

Now we can give the formal definition of stable $P$-marked ribbon graph.
\begin{definition*}
Consider a metrized (possibly disconnected) ribbon graph $G$
with an injection $x: P \hra \Sigma$ in a subset of ``special points''
$\Sigma \subset X_0(G) \sqcup X_\infty(G)$
such that $\Sigma\supset X_\infty(G)$, plus 
a fixed-point-free involution $\i$ acting on the set
of ``exceptional points'' $\Sigma \setminus x(P)$.
We say that an order function that assigns a natural
number to each connected component of $G$ is {\it admissible} if
\begin{enumerate}
\item[-]
components of order $0$ contain at least one $P$-marked hole
\item[-]
when $\i$ exchanges all the elements of $\Sigma\setminus x(P)$
belonging to the component $G_j$
with elements belonging to components of order $\leq k$,
then $G_j$ has order $\leq k+1$
\item[-]
every $h \in X_\infty(G)\setminus x(P)$ belongs to a component
of order $k>0$ and the point $\i(h)$ sits in a component of order
at most $k-1$ (and so is a vertex).
\end{enumerate}
We call $(G,x,\i)$ a $P${\it-marked stable ribbon graph} if
there exists an admissible order function on $G$.
A {\it stable metric} on $(G,x,\i)$ is the datum of a positive
metric $a_j$ of total length $1$
for every connected component $G_j$ of $G$.
\end{definition*}
Now let $\ua$ be a proper simplex of $A$ whose
associated marked ribbon graph is $G_{\ua}$. Consider the
set $\mathcal{Z}(G_{\ua})$ of connected stable subsets of $X_1(G_{\ua})$
and for every $Z \in \mathcal{Z}(G_{\ua})$
let $|\ua|^{\circ}\rar \D_Z \rar \D_{\ol{Z}}$ be the natural
map that projects first to a face and then to the space of metrics
of the reduced ribbon graph $G_{\ol{Z}}$.
Define $\hat{\ua}$ to be the closure of the graph of the map
\[
|\ua|^\circ \hra |\ua| \times
\prod_{Z \in \mathcal{Z}(G_{\ua})} \D_{\ol{Z}}
\]
in $|\ua| \times \prod \D_{\ol{Z}}$.

It can be proven that $\hat{\ua}$ parametrizes
all stable degenerations of the ribbon graph $G_{\ua}$.
Moreover,
all the $\hat{\ua}$'s can be glued to obtain a
modification $\Ahat$ of the arc complex.
Remark that $\Ahat(S,P)$ comes with an obvious cellularization
indexed by permissible sequences: for every
$Z_{\bullet}=(Z_0,\dots,Z_k)$ there is a (closed) cell
isomorphic to $|\ua_0|\times\cdots\times|\ua_k|$
that parametrizes stable metrics on $G(Z_{\bullet})$.
The projections $|\ua_0|\times\cdots\times|\ua_k| \rar |\ua_0|$
glue to give a continuous surjection $\Ahat(S,P) \rar |A(S,P)|$
which is actually a quotient (i.e. $|A(S,P)|$ has the quotient topology).

We may regard the map $\Ahat(S,P)/\G_{S,P}\rar|A(S,P)|/\G_{S,P}$ as a sort
of real oriented blow-up along a subcomplex, which is in fact the locus
on which $\Mbar_{g,P}\times\D_P \lra |A(S,P)|/\G_{S,P}$ is not a
homeomorphism.

\begin{theorem*}[\cite{looijenga:cellular}]
The modular group $\G_{g,P}$ naturally acts on $\Ahat(S,P)$
respecting the cellularization. The product of the classifying map
\[
\Ahat(S,P)/\G_{S,P} \lra \Mbar_{g,P}
\]
with $\hat{\l}$ is a continuous surjection
\[
\wh{\Phi}:
\Ahat(S,P)/\G_{S,P} \lra \Mbar_{g,P}\times \D_P
\]
that extends $\Phi$, so it is
one-to-one
when restricted to the dense open subset $|\Ao(S,P)|/\G_{S,P}$.
\end{theorem*}
We still denote by $\wh{\Phi}$ the map
\[
\wh{\Phi}: \Mhat_{g,P} \lra
\Mbar_{g,P}\times(\R_{\geq 0}^P\setminus\{0\}).
\]
where $\Mhat_{g,P}:=(\Ahat(S,P)/\G_{S,P})\times\R_+$.

Finally, the following commutative diagram
\[
\xymatrix@C=3cm{
\Mhat_{g,P} \ar[r]^{\wh{\Phi}} \ar[d] &
\Mbar_{g,P}\times(\R_{\geq 0}^P\setminus\{0\}) \ar[d]^{\xi} \\
|A(S,P)|/\G_{S,P}\times\R_+ \ar[r]^{\ol{\Phi}} & \Mbartri_{g,P}\times\R_+
}
\]
shows that $\ol{\Phi}$ is a homeomorphism, because
it is bijective and proper and both the vertical arrows are quotients.
%
%
%
\subsection*{Combinatorial classes on $\Mhat_{g,P}$}
As the space $\Mcomb_{g,P}$ includes in $\Mhat_{g,P}$,
we can define $\What_{m_*,P}$ as the closure of
$W_{m_*,P}$ inside $\Mhat_{g,P}$ (similarly
for $\What_{m_*,\rho,P}$ and $W_{m_*,\rho,P}$).

Notice that the map $\Mhat_{g,P} \lra \Mbarcomb_{g,P}$
sends $\What_{m_*,P}$ onto $\Wbar_{m_*,P}$ with
degree $1$. Hence, if $\What_{m_*,P}(l)$ were a cycle,
then its image in $\Mbar_{g,P}$ would be a lift
of $\Wbar_{m_*,P}(l)$. We can weaken this requirement:
if the image at the level of chains $\wh{\Phi}_{l,*}(\What_{m_*,P})$
were a cycle in $\Mbar_{g,P}$, then it would be
a lift of $\Wbar_{m_*,P}$.

Anyway, the situation is not so simple as
for $\Wbar_{m_*,P}$. In fact, if we consider
a simplex $\ua$ in $|A(S,P)|$, then $\pa\ua$
is made of simplices corresponding to ribbon
graphs that are obtained from $G_{\ua}$
by contracting one edge.
On the contrary, if we consider $\hat{\ua}$
in $\Ahat(S,P)$, then we may obtain a simplex
in $\pa\hat{\ua}$ by collapsing a stable
subgraph $Z$ (for instance, look at Fig. \ref{fig:What},
where a cell of real codimension one in the
boundary of $\What^q_5$ is obtained collapsing
the subgraph $Z$ of the edges sitting in the grey zone).
\begin{figure}[h]
\resizebox{12cm}{!}{\input{What-color.pstex_t}}
\caption{A cell $\hat{\ub}$ belonging to $\pa\What^q_5$ is obtained
from $\hat{\ua}\subset \What^q_5$ collapsing the edges in the dark zone}
\label{fig:What}
\end{figure}

As an example, we analyze the case $\What^q_{2r+3}$.
The situation is quite simpler than in the general case,
because the only new cells of real codimension one occurring
in the topological boundary $\pa\What^q_{2r+3}$ are indexed
by ribbon graphs that are very similar to the ribbon
graph of Fig. \ref{fig:What}.
In fact, let $\hat{\ub}$ be a cell in $\pa\What^q_{2r+3}$.
The associated ribbon graph $G_{\hat{\ub}}$
must have one positive component $G_{\hat{\ub}}^+$
and one nonpositive component $G_{\hat{\ub}}^0$.
Moreover, we can consider only the case in which $G_{\hat{\ub}}^0$
has one exceptional hole. In fact, if $G_{\hat{\ub}}^0$ had
more exceptional holes, then $\wh{\Phi}(\hat{\ub})$
would have one dimension less than $\hat{\ub}$
and so it would give no contribution
to $\pa\wh{\Phi}_*(\What^q_{2r+3})$.
Notice that, because of this hypothesis, the $(2r+3)$-valent
vertex has split into the exceptional vertex
$(f_1,\dots,f_{v_{ex}})$ of $G_{\hat{\ub}}^+$ and
the $q$-marked vertex $(e_1,\dots,e_{v_q})$
in $G_{\hat{\ub}}^0$.

Consider a cell $\hat{\ua} \in \What^q_{2r+3}$
such that $\hat{\ub}\in\pa\hat{\ua}$.
Then, the $q$-marked vertex in $G_{\hat{\ua}}$ looks like
\[
(e_i,f_1,\dots,f_{t_1-1},e_{i+1},f_{t_1},\dots,f_{t_2-1},e_{i+2},\dots,
e_{i-1},f_{t_{(v_q-1)}+1},\dots,f_{v_{ex}}).
\]

Now we define two operators that permute cells in $\What^q_{2r+3}$
that have $\hat{\ub}$ in their boundary.
The first operator $R_{ex}$ sends the $\hat{\ua}$ above to
the cell $R_{ex}(\hat{\ua})$ with $q$-marked vertex
\[
(e_i,f'_1,\dots,f'_{t_1-1},e_{i+1},f'_{t_1},\dots,f'_{t_2-1},e_{i+2},\dots,
e_{i-1},f'_{t_{(v_q-1)}+1},\dots,f'_{v_{ex}})
\]
where $f'_j=f_{j+1}$ for $j=1,\dots,v_{ex}-1$ and $f'_{v_{ex}}=f_1$.
The second operator $R_q$ sends the $\hat{\ua}$ above to the
cell $R_q(\hat{\ua})$ with $q$-marked vertex
\[
(e'_i,f_1,\dots,f_{t_1-1},e'_{i+1},f_{t_1},\dots,f_{t_2-1},e'_{i+2},\dots,
e'_{i-1},f_{t_{(v_q-1)}+1},\dots,f_{v_{ex}})
\]
where $e'_i=e_{i+1}$ for $i=1,\dots,v_{q}-1$ and $e'_{v_q}=e_1$.

Now, suppose that $v_q$ is even (otherwise, exchange the role
of $(v_q,R_q)$ and $(v_{ex},R_{ex})$).
If we prove that $\hat{\ua}$ and $R_q(\hat{\ua})$
induce opposite orientations on $\hat{\ub}$, then we obtain
that terms $\hat{\ub}$ in $\pa\What^q_{2r+3}$ sum up to zero.

To do that, we use the formalism of Conant and
Vogtmann (see \cite{conant-vogtmann:kontsevich}).

They proved that, for any connected ribbon graph $G$,
\[
\det(H_1(|G|))\otimes\det(\R X_1(G))\cong
\bigotimes_{V\in X_0(G)}\det(\R V)\otimes\det(\bigoplus_{|V|\,even}\R)
\]
while for a disconnected ribbon graph $G$ one has to choose
an ordering of the connected components.
In our case, we can order the connected components
of the ribbon graph $G_{\hat{\ub}}$ by prescribing that
the even vertex sits in the first component; the sign produced
by this choice cancels up with the standard choice in
$\det(\bigoplus_{|V|\,even}\R)$ given by the unique even vertex.

Moreover,
one can easily verify that
the map $|G^0_{\hat{\ub}}|\cup|G^+_{\hat{\ub}}| \lra |G_{\hat{\ua}}|$
that glues the $q$-marked vertex and the exceptional vertex
induces the following isomorphism
\[
H_1(|G^0_{\hat{\ub}}|)\oplus H_1(|G^+_{\hat{\ub}}|)
\cong
H_1(|G_{\hat{\ua}}|).
\]
Hence, the ratio between the orientations on $\hat{\ub}$
induced by $\hat{\ua}$
and $R_q(\hat{\ua})$ is easily seen to be $(-1)^{v_{ex}}=-1$
and so we have proven the following.
\begin{proposition*}
For every $r\geq 0$ and every $\tilde{l}=(l,0)\in\R_+^P\times\{0\}$,
the subcomplex $\What^q_{2r+3}\subset \Mhat_{g,P\cup\{q\}}$
pushes forward to a cycle
\[
\wh{\Phi}_{\tilde{l},*}(\What^q_{2r+3})
\in H_{6g-6+2n-2r}(\Mbar_{g,P\cup\{q\}})
\]
which is a lifting of $\Wbar^q_{2r+3}$.
As a consequence,
for every $r\geq 0$, the cycle
\[
(\pi_q \wh{\Phi}_{\tilde{l}})_*(\What^q_{2r+3})
\in H_{6g-6+2n-2r}(\Mbar_{g,P})
\]
is a lifting of $\Wbar_{2r+3}$.
\end{proposition*}
%
%
%
%
\bibliographystyle{amsalpha}
\bibliography{bibliografia}
\end{document}